\pdfoutput=1
\documentclass[12pt,oneside,a4paper,reqno,openany]{book}
\usepackage[french]{babel}
\usepackage[utf8]{inputenc}
\usepackage{babelbib}
\usepackage[T1]{fontenc}
\usepackage[left=2.8cm,right=2.7cm,top=2.8cm,bottom=2.9cm]{geometry}
\usepackage{srcltx}
\usepackage{multicol}
\usepackage{pdfpages,aeguill}
\usepackage{epic,eepic,color}
\usepackage{amsmath,mathrsfs,amscd,amssymb,amsfonts,latexsym,amsthm,epsf}
\usepackage{graphicx,textcomp}

\usepackage{tikz}

\usepackage{fancybox}
\usepackage{hyperref}

\usepackage[french,boxed,ruled,lined]{algorithm2e}

\newcommand{\beq}{\begin{equation}}
\newcommand{\eeq}{\end{equation}}
\newcommand{\beqn}{\begin{equation*}}
\newcommand{\eeqn}{\end{equation*}}

\theoremstyle{plain}
\newtheorem{thm}{Théorème}[chapter]
\newtheorem{cor}[thm]{Corollaire}
\newtheorem{prop}[thm]{Proposition}
\newtheorem{lem}[thm]{Lemme}

\theoremstyle{plain}
\newtheorem{déf}[thm]{Définition}
\newtheorem{df}[thm]{Définitions}
\newtheorem{rem}[thm]{Remarque}
\newtheorem{ex}[thm]{Exemple}

\newcommand{\bl}{\textcolor{blue}}
\newcommand{\bex}{\begin{ex}}
\newcommand{\eex}{\end{ex}}
\newcommand{\brem}{\begin{rem}}
\newcommand{\erem}{\end{rem}}

\newcommand{\bthm}{\begin{thm}}
\newcommand{\ethm}{\end{thm}}
\newcommand{\bpro}{\begin{prop}}
\newcommand{\epro}{\end{prop}}
\newcommand{\bco}{\begin{cor}}
\newcommand{\eco}{\end{cor}}
\newcommand{\blem}{\begin{lem}}
\newcommand{\elem}{\end{lem}}
\newcommand{\ed}{\end{déf}}
\newcommand{\bd}{\begin{déf}}
\newcommand{\bdf}{\begin{df}}
\newcommand{\edf}{\end{df}}
\addto{\captionsfrench}{ }
\newcommand{\bpr}{\begin{proof}}
\newcommand{\epr}{\end{proof}}

\newcommand{\D}{\mathcal{D}}

\newcommand{\rci}{\rm rci}
\newcommand{\T}{\rm T}

\newcommand{\red}{{\rm red}}

\newcommand{\fp}{{\rm fp}}
\newcommand{\exc}{{\rm exc}}
\newcommand{\inv}{{\rm inv}}
\newcommand{\maj}{{\rm maj}}
\newcommand{\crs}{{\rm cr}}
\newcommand{\nes}{{\rm nes}}
\newcommand{\des}{{\rm des}}

\newcommand{\card}{{\rm card}}
\newcommand{\rt}{{\rm rt}}
\newcommand{\ct}{{\rm ct}}
\newcommand{\RSK}{R\!S\!K}

\newcommand{\ut}{{\rm ut}}
\newcommand{\lt}{{\rm lt}}

\renewcommand{\tablename}{Tableau}

\title{
	\begin{center}
		\Large\textbf{UNIVERSITE D'ANTANANARIVO}\\ °°°°°°°°°°°°°°°°\\ {\upshape Sciences et Technologies\\ °°°°°°°°°°°°°°°°\\
				Mention Mathématiques et Informatique} \\ 	
		\vspace{1cm}
		\textbf{THÈSE DE DOCTORAT\\
				de Mathématiques et Informatique} \\		
		{\normalsize  			Spécialité: \textbf{COMBINATOIRE}} \\
		\vspace{1.2cm} 
		\textbf{ \bl{COMBINATOIRE DES PERMUTATIONS RESTREINTES SELON LE NOMBRE DE CROISEMENTS}}
		\vspace{1.2cm}
		\\		
{\normalsize Soutenu publiquement  le \textit{\textbf{---------------------}} par
 \\
\textbf{ Paul Mazoto RAKOTOMAMONJY}
 \\
 Devant le jury composé de 
}
\end{center}
{ \small\begin{table}[h]
	\begin{center}
		\begin{tabular}{lll}
			Pr  Arthur RANDRIANARIVONY&	Université d'Antananarivo&	Directeur de thèse\\
			-------------&	------------- &	Examinateur\\
			-------------&	------------- &	Président de jury\\
			-------------& ------------- &	Examinateur\\
			------------- &	------------- &	Rapporteur externe\\
			-------------&	------------- &	Rapporteur interne
		\end{tabular}
	\end{center}
\end{table}
}
	\date{}
}

\begin{document}
 \renewcommand{\tablename}{Tableau}
 \renewcommand{\figurename}{Figure}
 \renewcommand{\listfigurename}{Liste des figures}
\includepdf{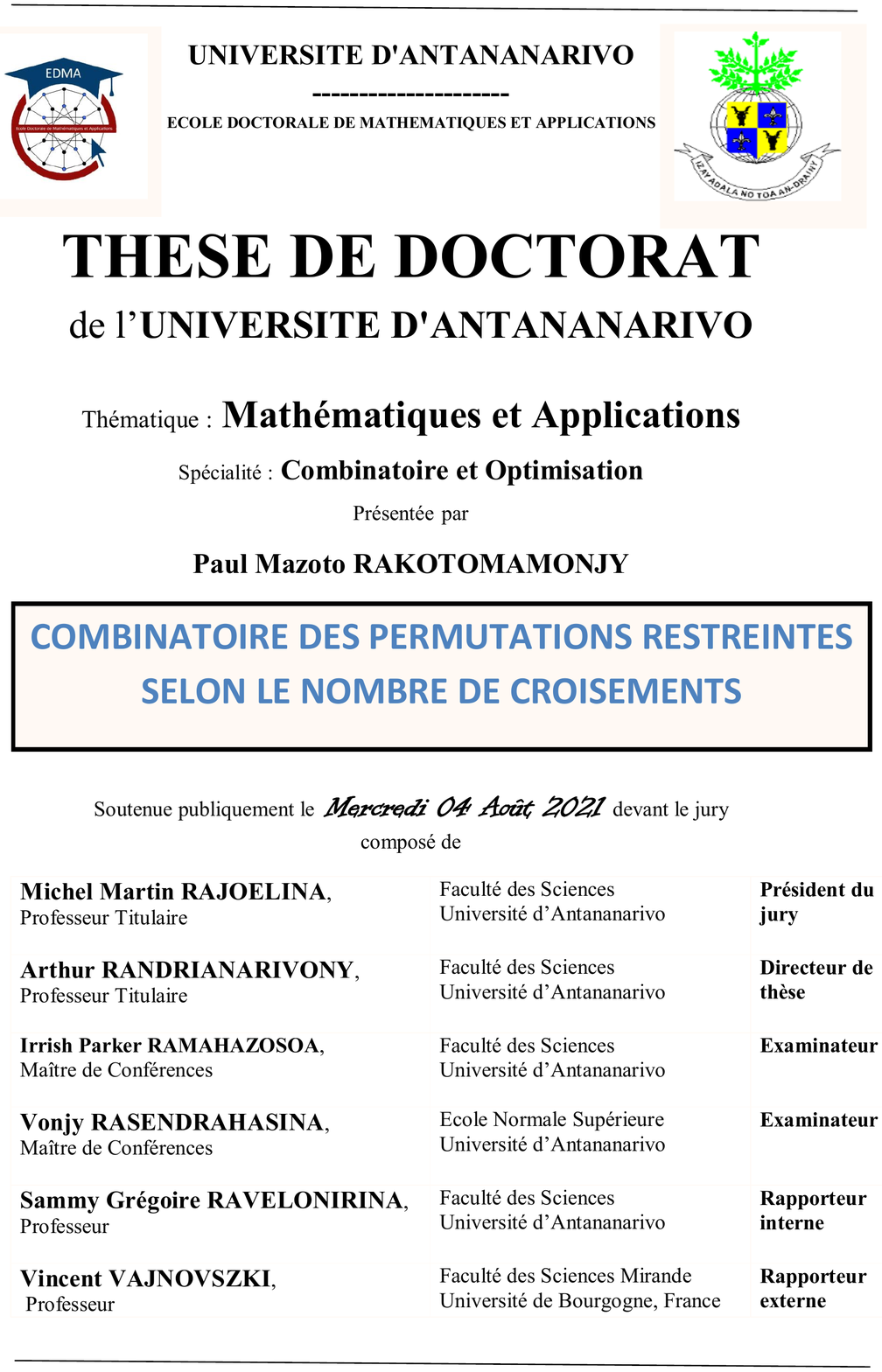}
\frontmatter
\pagestyle{plain}

  \chapter{Remerciements}
    Sans la grâce de notre Seigneur Dieu, cette thèse n'aurait pu être menée à son terme. Je n'ai pas assez de mots pour le remercier. Merci mon Dieu.

     J'exprime ma  profonde gratitude à mon directeur de thèse, Monsieur \textbf{Arthur RANDRIANARIVONY}, qui, depuis le Master, m'a fait découvrir un sujet de recherche passionnant et m'a consacré du temps malgré ses multiples occupations. La confiance qu'il m'a accordée a abouti à des résultats qui ont largement dépassé mes espérances initiales. Je lui en suis  infiniment reconnaissant.
     
     Je tiens également à adresser mes remerciements à tous les membres du jury.
\begin{itemize}
	\item[-] A Monsieur \textbf{Michel Martin RAJOELINA}  qui me fait l'honneur de sa présence en tant que président du jury.
\item[-]  A Messieurs \textbf{Vonjy  RASENDRAHASINA} et \textbf{Irrish Parker RAMAHAZOSOA} qui ont  gentiment accepté d'examiner cette thèse.
\item[-] A Messieurs \textbf{Hanitriniaina Sammy Grégoire RAVELONIRINA} et \textbf{Vincent VAJNOVSZKI}. Je suis très reconnaissant à l'honneur que vous me faites en acceptant d'être les rapporteurs de cette thèse.
\end{itemize}

    Je ne saurais oublier mes remerciements à  ceux  qui ont contribué de près ou de loin à l'avancement de ma thèse:
    \begin{itemize}
    	\item[-]  A Monsieur \textbf{Sandrataniaina ANDRIANTSOA}, un de mes co-auteurs et collègues, un frère avec qui j'ai partagé mes problèmes tout le temps.
    	\item[-]  A Messieurs \textbf{Alain RALAMBO}, \textbf{Olivier ROBINSON} pour les multiples discussions et encouragement.
     	\item[-] A l'\textbf{IT University} car, sans son aide financière qui m'a donné de stabilité,  cette thèse n'aurait pas vu le jour.  Dès fois, pour la recherche, j'ai manqué des cours  et  préparations.     
 		\item[-] A tous les enseignants et les responsables de la mention  Mathématiques et Informatique, Sciences et Technologies, Université d'Antananarivo. Cette thèse est l'un des fruits de leurs efforts dans l'éducation.
  		\item[-] A ma femme \textbf{Voahangy}, mes enfants \textbf{Fy} et \textbf{Mahandry} qui ont souffert de mes longues périodes d'absence lors de la préparation de cette thèse. 
    		\item[-] A toute ma famille pour le soutien et l'encouragement qu'elle a manifestés envers mes études. \textbf{\textit{\og Ny\! hazo\! no\! vanon-ko\! lakana,\! ny\! tany\! naniriany\! no\! tsara\fg}}.
    \end{itemize}

J'espère, à travers ces quelques lignes, n'avoir oublié personne. Que tous ceux qui auraient dû être cités ici reçoivent l'expression de mon plus grand respect et de mon estime.

\newpage
\chapter{Avant-propos}

	Cette thèse est située dans le domaine de la combinatoire énumérative, une branche des mathématiques discrètes où l'étude des statistiques sur les permutations restreintes occupe au cours de  ces dernières années une place importante grâce à ses multiples applications en informatique, en probabilité et en combinatoire elle-même. Elle est le fruit de nos deux articles \cite{Rakot, Rakot2} préparés depuis quelques années:
\begin{itemize}
	\item[$\bullet$]  \cite{Rakot}, \textit{Restricted permutations refined by number of crossings and nestings}, Discrete
Math. 343 (2020) 111950.
	\item[$\bullet$]
\cite{Rakot2}, \textit{Crossings over permutations avoiding some pairs of patterns of length three}, J.~Integer Sequences 23 (2020), Article 20.6.3. (En collaboration avec Andriantsoa et Randrianarivony). 
\end{itemize}
Dans ces papiers, nous avons introduit les études combinatoires des permutations restreintes selon le nombre de croisements. L'objectif principal est de trouver les distributions du nombre de croisements sur les permutations interdisant un ou deux motifs de l'ensemble $S_3$ des permutations de longueur 3.

Les outils que nous avons utilisés viennent naturellement de diverses œuvres. D'abord, Elizalde et Pak \cite{ElizP} ont défini une bijection $\Theta:S_n(321)\rightarrow S_n(132)$ conservant deux statistiques bien connues en combinatoire. Nous avons besoin de cette bijection pour prouver bijectivement l'équidistribution du nombre de croisements sur les ensembles $S_n(321)$ et $S_n(132)$. Ensuite, Randrianarivony \cite{ARandr} a défini un q,p-analogue des nombres de Catalan $C_n(q,p)$ et a interprété combinatoirement $C_n(q,p)$ en termes de croisements sur les permutations sans imbrication. Dû au fait que toute  permutation sans imbrication est une permutation 321-interdite, son résultat est lié directement au nôtre; ce qui nous permet d'exprimer en termes de fraction continue la distribution du nombre de croisements sur les ensembles $S_n(321)$ et $ S_n(132)$ et  $ S_n(213)$. Enfin, en nous inspirant  des techniques utilisées par Dokos et al. \cite{Dokos} et Elizalde \cite{Eliz1} concernant l'étude de diverses statistiques usuelles sur l'ensemble des permutations interdisant un ou plusieurs motifs de $S_3$, nous avons prouvé dans \cite{Rakot} que la composition miroir-complément-inverse des involutions triviales sur les permutations conserve le nombre de croisements. Ce dernier a beaucoup facilité l'énumération des permutations interdisant deux motifs de $S_3$ selon le nombre de croisements \cite{Rakot2}.

\tableofcontents
\listoffigures \addcontentsline{toc}{chapter}{\listfigurename}
\listoftables  \addcontentsline{toc}{chapter}{\listtablename}
\chapter*{Notations}
\addcontentsline{toc}{chapter}{Notations}

\begin{tabular}{ll}
$|\sigma|$ & La longueur de $\sigma$, où $\sigma$ est une permutation.\\ 
$\sigma(a\cdots b)$ & La sous-séquence contiguë de $\sigma$ de la $a$-ème jusqu'à la $b$-ème lettre.\\
	&  $\sigma$ est une permutation, $a$ et $b$ sont des entiers tels que $a\leq b$.\\
$\sigma^{+a}$ & La permutation déduite de  $\sigma$ en ajoutant $a$ à chacun de ses nombres.\\
$\sigma^{i\rtimes a}$ & La permutation déduite de $\sigma$ en ajoutant $a$ à tous les nombres de $\sigma$ \\
			& qui sont plus grands ou égaux à $i$.\\
$\sigma^{(i,a)}$ & La permutation déduite de $\sigma$ en  augmentant d'une unité tous les entiers\\
	& supérieurs ou égaux à $a$ dans $\sigma$, puis en insérant l'entier $a$ à la $i$-ème position\\
$\crs(\sigma)$ & Le nombre de croisements d'une permutation $\sigma$.\\
$\D_n$ & Ensemble de chemins de Dyck de demi-longueur $n$, $n$ est un entier non nul.\\ 
$D^{(L)}$& La moitié gauche d'un chemin de Dyck $D$.\\
$D^{(R)}$& La moitié droite d'un chemin de Dyck $D$.\\
$|E|$ & Le cardinale de $E$, où $E$ est un ensemble.\\ 
$F_n(T;q)$ & Le polynôme distributeur du nombre de croisements sur $S_n(T)$. C'est-à-dire,\\
& $\displaystyle F_n(T;q)=\sum_{\sigma \in S_n(T)}q^{\crs(\sigma)}$.\\
$F(T;q,z)$ & La fonction génératrice ordinaire de $F_n(T;q)$, c'est-à-dire, \\
& $\displaystyle F(T;q,z)=\sum_{n\geq 0}F_n(T;q)z^n.$\\
$F_n^k(q)$ & Le polynôme distributeur du nombre de croisements sur $S_n^k$. \\
$F_n^k(T;q)$ & Le polynôme distributeur du nombre de croisements sur $S_n^k(T)$.\\
$F_{n,k}(q)$ & Le polynôme distributeur du nombre de croisements sur $S_n^k$. \\
$F_{n,k}(T;q)$ & Le polynôme distributeur du nombre de croisements sur $S_{n,k}(T)$.\\
$[n]$ & L'ensemble $\{1,2,\ldots, n\}$ pour tout entier $n\geq 1$.\\
$\nes(\sigma)$ & Le nombre d'imbrications d'une permutation $\sigma$.\\
$S_n$ & Ensemble de permutations de $[n]$ pour tout entier $n\geq 1$.\\
$S_n^k$ & Ensemble de permutations $\sigma$ de $[n]$ tel que $\sigma(k)=1$.\\
$S_{n,k}$ & Ensemble de permutations $\sigma$ de $[n]$ tel que $\sigma(n)=k$.\\ 
$S_n(\tau)$ & Ensemble de permutations de $[n]$ interdisant le motif $\tau$, où $\tau$ est une\\
& permutation.\\
$S_n(T)$ & Ensemble de permutations de $[n]$ interdisant tous les motifs dans $T$,\\
& où $T$ est un ensemble de motifs.
\end{tabular}

\mainmatter
\addcontentsline{toc}{chapter}{Introduction générale}	
	\chapter*{Introduction générale}
	
Les concepts de croisement et d'imbrication sur les  permutations ont été introduits  par M\'edicis et Viennot \cite{MedVienot} et plusieurs auteurs poursuivent encore leurs études \cite{Burril, Cort, Cort2, ARandr1, ARandr}. Il est bien connu que les deux statistiques $\crs$ et $\nes$ sont équidistribuées. C'est-à-dire, pour tout $n\geq 0$, 
\begin{equation*}
\sum_{\sigma \in S_n}q^{\crs(\sigma)}=\sum_{\sigma \in S_n}q^{\nes(\sigma)}.
\end{equation*}
Il a aussi été prouvé que le polynôme $\sum_{\sigma \in S_n}x^{\crs(\sigma)}y^{\nes(\sigma)}$ est symétrique, c'est-à-dire,
\begin{equation*}
\sum_{\sigma \in S_n}x^{\crs(\sigma)}y^{\nes(\sigma)}=\sum_{\sigma \in S_n}y^{\crs(\sigma)}x^{\nes(\sigma)}
\end{equation*}
et que le développement en fraction continue de sa fonction génératrice est
\begin{equation*}
\sum_{n\geq 0}\left( \sum_{\sigma \in S_n}x^{\crs(\sigma)}y^{\nes(\sigma)}\right) z^n=
\frac{1}{1-\displaystyle\frac{[1]_{x,y}~.z}{				
		1-\displaystyle\frac{[1]_{x,y}~.z}{
			1-\displaystyle\frac{[2]_{x,y}~.z}{								
				1-\displaystyle\frac{[2]_{x,y}~.z}{
					1-\displaystyle\frac{[3]_{x,y}~.z}{
						1-\displaystyle\frac{[3]_{x,y}~.z}{				
							\ddots}
					}
			}}		
}}}, \label{arthurfc}
\end{equation*}
où $[n]_{a,b}:=a^{n-1}+a^{n-2}b+\ldots+ab^{n-2}+b^{n-1}$. On peut trouver deux versions différentes des preuves bijectives de ces résultats dans \cite{Cort, ARandr}.

Qu'en est-il des distributions de ces deux statistiques sur les permutations à motifs interdits? Quelles expressions, relations et interprétations combinatoires peut-on trouver? 

Motivé par ces questions ouvertes, cette thèse  est le fruit de nos deux articles \cite{Rakot, Rakot2}, dans lesquels nous sommes particulièrement intéressés par l'énumération des permutations interdisant un ou deux motifs de longueur 3 selon le nombre de croisements. Dans notre premier article \cite{Rakot}, nous avons trouvé les identités suivantes:
\begin{equation}
\sum_{\sigma\in S_n(321)} q^{\crs(\sigma)}=\sum_{\sigma\in S_n(132)} q^{\crs(\sigma)}=\sum_{\sigma\in S_n(213)} q^{\crs(\sigma)}.\label{eqmain1}
\end{equation}
Outre les techniques que nous nous sommes inspirées de certains travaux \cite{Dokos,Eliz1,Eliz2,ARob,ARob2}, nous utilisons une bijection  $\Theta:S_n(321)\rightarrow S_n(132)$ qui était initialement construite par Elizalde et Pak dans \cite{ElizP}. Pour tout $\tau \in \{321,132,213\}$, nous obtenons grâce au résultat de Randrianarivony \cite{ARandr} l'identité inattendue suivante: 
\begin{equation}
\sum_{\sigma\in S(\tau)} q^{\crs(\sigma)}z^{|\sigma|}=\frac{1}{1-\displaystyle\frac{z}{				
		1-\displaystyle\frac{z}{
			1-\displaystyle\frac{qz}{								
				1-\displaystyle\frac{qz}{
					1-\displaystyle\frac{q^2z}{
						1-\displaystyle\frac{q^2z}{				
							\ddots}
					}
			}}		
}}}. \label{eqmain2}
\end{equation}
Pour l'instant, trouver les distributions du nombre de croisements sur les permutations interdisant l'un des motifs $123$, $231$ et $312$ reste un problème ouvert. Nous observons que la fraction continue  \eqref{eqmain2} est apparue dans \cite{BEliz} comme la distribution du nombre d'occurrences d'un motif sur les permutations 231-interdites. La recherche d'éventuelles correspondances entre ces résultats s'avère intéressante, vu que  Corteel \cite{Cort} a établi la liaison entre les occurrences de motif, les croisements et les imbrications sur les permutations.\\
Pour tout ensemble de motifs $T$, notons $\displaystyle F_n(T;q):=\sum_{\sigma \in S_n(T)}q^{\crs(\sigma)} \text{ et }  
F(T;q,z):=\sum_{\sigma \in S(T)}q^{\crs(\sigma)}z^{|\sigma|}=\sum_{n\geq 0}F_n(T;q)z^n$.
Nous écrivons $F_n(\tau_1,\tau_2,\ldots;q)$ et $F(\tau_1,\tau_2,\ldots;q,z)$ si $T=\{\tau_1,\tau_2,\ldots\}$.  Incité par les résultats cités précédemment, nous avons poursuivi nos recherches dans \cite{Rakot2}. Comme résultats, nous avons trouvé les différentes relations suivantes qui lient les distributions de $\crs$ sur les permutations interdisant les motifs 231 et 312:
\begin{align*}
		& F(312;q,z)=\displaystyle \frac{1}{1-zF(231;q,z)},\\
		& F(312,123; q,z)= 1+\left(\frac{z}{1-z} \right)^2 +zF(231,123; q,z),\\	
		& F(312,\tau; q,z)= 1+\left( \frac{z}{1-z}\right)F(231,\tau'; q,z), \text{ pour tout  $(\tau,\tau') \in   \{132,213\}^2$}.	
\end{align*}
Nous avons aussi trouvé les deux résultats d'énumérations suivants:
	\begin{align}
		F(231,321;q,z)&=\frac{1-qz}{1-(1+q)z-(1-q)z^2},\label{main1x}\\	
		F(123,\tau;q,z)&=1+\frac{(1-qz)z}{(1-z)(1-(1+q)z)}, \text{ pour tout  $\tau \in   \{132,213\}$}\label{main1y}.	
	\end{align}
A travers une œuvre récente de  Bukata et al. \cite{BKLPRW}, nous avons observé que les identités \eqref{main1x} et \eqref{main1y} sont respectivement des nouvelles interprétations combinatoires des triangles \href{https://oeis.org/A076791}{A076791} et \href{https://oeis.org/A299927}{A299927} de OEIS \cite{OEIS} (ou On-line Encyclopedia of Integer Sequences). Bukata et al. ont interprété ces triangles en termes d'autres statistiques sur les permutations évitant une paire de motifs de longueur 3 (voir \cite[Prop. 7 et Prop. 11]{BKLPRW}).
Pour prouver ces résultats, nous avons combiné les méthodes bijectives et les fonctions génératrices.

Cette thèse sera organisée en trois chapitres. 
Le premier chapitre est destiné aux définitions et outils préliminaires qui seront utiles pour la suite. En particulier, nous y rappellerons la bijection $\Theta:S_n(321)\rightarrow S_n(132)$ de Elizalde et Pak qui  aura  une importance capitale pour la preuve des résultats du chapitre suivant.  
Dans les deux derniers chapitres, nous examinerons les distributions du nombre de croisements sur les permutations interdisant un ou deux  motifs de $S_3$ et nous discuterons comment nos résultats sont liés à ceux de \cite{BKLPRW,Dokos,MansSh,ARandr,ARob,Sarino}.

	\chapter{Outils préliminaires}

		Cette partie a pour objectif de fournir quelques définitions et outils préliminaires qui seront nécessaires pour la suite. 
\section{Quelques définitions}
\subsection{Permutation à motif interdit}
\bd {\rm Soit $E\subset\mathbb{N}$ tel que $|E|=n>0$. 
Une \textit{permutation} $\sigma$ de $E$ est une bijection de $E$ dans $E$. La \textit{réduction} d'une permutation $\sigma$ de $E$, notée $\red(\sigma)$, est une permutation de $[n]:=\{1,2,\ldots,n\}$ définie par $\red(\sigma):=\tau\circ \sigma \circ \tau^{-1}$, où $\tau$ est l'unique bijection croissante de $E$ vers $[n]$.}
\ed

Habituellement, si $E=\{e_1,e_2,\ldots,e_n\}$ tel que $e_1<e_2<\ldots<e_n$, une permutation $\sigma$ de $E$ peut s'écrire en terme de bi-mot comme suit
\begin{equation*}
\sigma=\left( \begin{array}{cccc}
e_1 & e_2 & \ldots & e_n\\
\sigma(e_1) & \sigma(e_2) & \ldots & \sigma(e_n)
\end{array} \right).
\end{equation*}
Mais généralement, nous préférons l'écriture linéaire qui est en une seule ligne  
\begin{equation*}
\sigma=\sigma(e_1)~ \sigma(e_2)~ \ldots~ \sigma(e_n).
\end{equation*} 
 Nous notons $S_n$ l'ensemble des permutations de $[n]$ et nous désignons par $|\sigma|$ la longueur d'une permutation $\sigma$. \\
 Par exemple, la permutation
	$\sigma=\left( \begin{array}{ccccc}
	1 & 3 &6 & 8& 9\\
	3 & 9& 1&  8 & 6
	\end{array} \right)$ de $\{1,3,6,8,9\}$ s'écrit simplement $\sigma=3~ 9~ 1~ 8~ 6$ et sa réduction est $\red(\sigma)=2~ 5~ 1~ 4~ 3 \in S_5$. 
\bd	
{\rm Soit $\sigma\in S_n$ et $\tau\in S_m$ deux permutations telles que $m\leq n$. On dit que $\sigma$ \textit{contient} le motif $\tau$ s'il existe une sous-suite $\alpha=\sigma(i_1)\cdots \sigma(i_m)$ de $\sigma$ (avec $i_1<i_2<\cdots<i_m$) telle que  $\red(\alpha)=\tau$. 
Dans le cas contraire, on dit que $\sigma$ est \textit{$\tau$-interdite}.}
\ed
Pour tout motif $\tau$ donné, on note $S_n(\tau)$ l'ensemble des permutations $\tau$-interdite de $[n]$. La permutation $\sigma=25143\notin S_5(132)$ car $\alpha_1=253$, $\alpha_2=254$ et $\alpha_3=143$ sont des sous-suites telles que $\red(\alpha_1)=\red(\alpha_2)=\red(\alpha_3)=132$. Par contre, on a $\sigma \in S_5(123)$.

 Généralement, si $T$ est un ensemble de motifs, nous notons
l'ensemble des permutations de $[n]$ interdisant tous les motifs dans $T$ par $S_n(T)=\bigcap_{\tau \in T}S_n(\tau)$
et l'ensemble de toutes permutations interdisant tous les motifs dans $T$ par
$S(T)=\bigcup_{n\geq 0}S_n(T)$. 

\subsection{Statistiques sur les permutations}
\bd
{\rm Une \textit{statistique} $s$ sur un ensemble $E$ est une application de $E$ dans $\mathbb{N}$. Le \textit{polynôme distributeur} ou simplement la \textit{distribution} d'une statistique $s$  sur  $E$ est le polynôme en $q$:   
	\begin{equation*}
	S(q)=\sum_{e \in E} q^{s(e)}.
	\end{equation*} 
	Dans le cas général, la distribution jointe de $r$ statistiques $s_1,s_2,\ldots ,s_r$ sur $E$ est 
	\begin{equation*}
	S(q_1,q_2,\ldots, q_r)=\sum_{e \in E} q_{1}^{s_1(e)}q_{2}^{s_2(e)}\ldots q_{r}^{s_r(e)}. 
	\end{equation*} }
\ed
\bex
{\rm L'application $\card$ définie sur l'ensemble $E_n$  des parties d'un ensemble $X$ à $n$ éléments par $\card(A)=|A|$ =nombre d'éléments de $A$ pour tout $A\in E_n$ est une statistique sur $E_n$  dont la distribution est
\begin{equation*}
D_n(q)=\sum_{A \in E_n} q^{\card(A)}= (1+q)^n.
\end{equation*}

En effet, nous avons
\begin{eqnarray*}
	\sum_{A \in E_n} q^{\card(A)}&=&\sum_{k=0}^{n}\sum_{A \in E_n, \card(A)=k} q^{k}= \sum_{k=0}^{n}\binom{n}{k} q^{k}=(1+q)^n.
\end{eqnarray*}
}
\eex
\bd
{\rm On dit que deux statistiques $s_1$ et $s_2$ sont \textit{équidistribuées} sur un ensemble $E$ si elles ont la même distribution, c'est-à-dire, 
\begin{equation*}
\sum_{e \in E} q^{s_1(e)}= \sum_{e \in E} q^{s_2(e)}.
\end{equation*}}
\ed
Voici quelques statistiques usuelles sur l'ensemble des permutations. Pour toute permutation $\sigma\in S_n$, on définit par

$\begin{array}{ll}
\exc(\sigma)&:=|\{i: \sigma(i)>i\}| \text{ le \textit{nombre d'excédances} de $\sigma$};\\
\fp(\sigma)&:=|\{i: \sigma(i)=i\}| \text{ le \textit{nombre de points fixes} de $\sigma$};\\
\des(\sigma)&:=|\{i: \sigma(i)>\sigma(i+1)\}| \text{ le \textit{nombre de descentes} de $\sigma$};\\
\inv(\sigma)&:=|\{(i,j): i<j \text{ et } \sigma(i)>\sigma(j)\}| \text{ le \textit{nombre d'inversions } de $\sigma$};\\
\maj(\sigma)&:=\displaystyle \sum_{\substack{i\in [n-1], \sigma(i)>\sigma(i+1)}}i \quad \text{ l'\textit{indice majeur} de  $\sigma$ }.
\end{array}$

Sur $S_n$, nous avons, par exemple, des équidistributions entre les statistiques $\exc$ et $\des$, d'une part, et entre  les statistiques $\inv$ et $\maj$, d'autre part. Ainsi, d'après MacMahon \cite{MacM}  et Foata \cite{Foata}, nous avons les identités
\begin{equation}\label{eul-mah}
\sum_{\sigma \in S_n} q^{\des(\sigma)}=\sum_{\sigma \in S_n} q^{\exc(\sigma)} \text{ et } \sum_{\sigma \in S_n} q^{\inv(\sigma)}=\sum_{\sigma \in S_n} q^{\maj(\sigma)}, \text{ pour tout }n\geq 1.
\end{equation}

Nombreux auteurs ont étudié les statistiques $\exc$, $\fp$, $\des$, $\maj$ et $\inv$, non seulement sur les permutations en général, mais aussi sur les permutations à motif interdit (voir par exemple \cite{Dokos, ElizD, ElizP, Eliz1, Eliz2, ARob,ARob2}). Dans ce manuscrit, nous nous intéressons à la statistique \textit{\og nombre de croisements\fg{}}. 
\bdf
{\rm Un \textit{croisement} d'une permutation $\sigma$ de $S_n$ est 
un couple $(i,j)$ tel que $i<j<\sigma(i)<\sigma(j)$ ou $\sigma(i)<\sigma(j)\leq i<j$. 

Une \textit{imbrication}  d'une permutation $\sigma$ de $S_n$ est également un couple $(i,j)$ vérifiant $i<j<\sigma(j)<\sigma(i)$ ou $\sigma(j)<\sigma(i)\leq i<j$. }
\edf
Pour toute permutation $\sigma$, nous dénotons respectivement par $\crs(\sigma)$ et $\nes(\sigma)$  le nombre de croisements et le nombre d'imbrications de $\sigma$.  Habituellement, pour une meilleure compréhension, on peut dessiner les diagrammes d'arcs d'une permutation.

	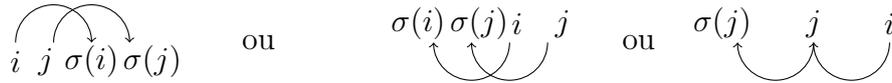
\begin{figure}[h]
	\centering
	
	\begin{tikzpicture}
	\draw[<-] (2.3,0) arc(0:180:0.5);
	\draw[<-] (2.8,0) arc(0:180:0.5);
	\draw (1.3,-0.25) node {$i$};
	\draw (2.3,-0.25) node {$\sigma(i)$};
	\draw (1.7,-0.25) node {$j$};
	\draw (3.1,-0.25) node {$\sigma(j)$};
		\draw (4.5 ,0) node {ou};

	\draw[->] (7.8,0) arc(0:-180:0.5);
	\draw[->] (8.3,0) arc(0:-180:0.5);
	\draw (8.5,0.25) node {$j$};
	\draw (7.4,0.25) node {$\sigma(j)$};
	\draw (7.9,0.25) node {$i$};
	\draw (6.6,0.25) node {$\sigma(i)$};
	\draw (9.5 ,0) node {ou};

	\draw[->] (11.8,0) arc(0:-180:0.5);
	\draw[->] (12.8,0) arc(0:-180:0.5);
	\draw (12.8,0.25) node {$i$};
	\draw (11.8,0.25) node {$j$};
	\draw (10.6,0.25) node {$\sigma(j)$};
	\end{tikzpicture}
	\label{fig:cross}
		\caption{Diagramme  d'arcs pour un croisement.}
\end{figure}
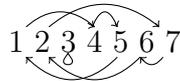
\begin{figure}[h]
	\centering
	
	\begin{tikzpicture}
	\draw[black] (0,1) node {$1\ 2\ 3\ 4\ 5\ 6\ 7$}; 
	\draw (-0.9,1.2) parabola[parabola height=0.2cm,red] (-0,1.2);\draw[->,black] (-0.09,1.25)--(-0,1.2);			
	\draw (-0.6,1.2) parabola[parabola height=0.3cm,red] (0.85,1.2); \draw[->,black] (0.8,1.25)--(0.85,1.2);			
	\draw (-0.45,0.7) -- (-0.35,0.85)[rounded corners=0.1cm] -- (-0.25,0.7) -- cycle;			
	\draw (-0.03,1.2) parabola[parabola height=0.2cm,red] (0.3,1.2); \draw[->,black] (0.27,1.25)--(0.3,1.2);			
	\draw (-0.9,0.8) parabola[parabola height=-0.3cm,red] (0.3,0.8); \draw[->,black] (-0.85,0.75)--(-0.9,0.8);			
	\draw (-0.6,0.8) parabola[parabola height=-0.3cm,red] (0.6,0.8); \draw[->,black] (-0.55,0.75)--(-0.6,0.8);			
	\draw (0.6,0.8) parabola[parabola height=-0.1cm,red] (0.9,0.8); \draw[->,black] (0.65,0.75)--(0.6,0.8);			
	\end{tikzpicture}
	\label{fig:arcdiag}
	\caption{Diagramme d'arcs de la permutation $\pi=4735126 \in S_7$.}
\end{figure}

\bex
{\rm Les croisements de la permutation $\pi=4735126\in S_7$ dessinée dans la Figure \ref{fig:arcdiag} sont $(1,2)$, $(5,6)$ et $(6,7)$. Nous avons alors $\crs(\pi)=3$. La permutation $\pi$ possède également trois  imbrications qui sont $(2,4)$, $(3,5)$ et $(3,6)$, c'est-à-dire, $\nes(\pi)=3$.
}
\eex

Il est bien connu que les deux statistiques $\crs$ et $\nes$ sont équidistribuées, c'est-à-dire,  
\begin{equation*}
\sum_{\sigma \in S_n}q^{\crs(\sigma)}=\sum_{\sigma \in S_n}q^{\nes(\sigma)} \text{ pour tout $n\geq 0$}.
\end{equation*}

\subsection{Fonction génératrice d'une suite}
Parlons maintenant d'un des outils d'énumération les plus utilisés en combinatoire énumérative.
\bd
{\rm La \textit{fonction génératrice ordinaire} (f.g.o) d'une suite $(U_n)_n$ est la série formelle
	\begin{equation*}
	U(z)=\sum_{n\geq 0}U_n z^n.
	\end{equation*}}
\ed	

\bex
{\rm Soit $A(z)$ la f.g.o d'une suite arithmétique $(a_n)_n$ de raison $r$ et de premier terme $a_0$.
Nous avons
\begin{align*}
A(z)&=a_0+\sum_{n\geq 1}(a_{n-1}+r)z^n \\
&= a_0+z\sum_{n\geq 0}a_{n}z^{n}+rz \sum_{n\geq 0}z^{n}\\
&= a_0+zA(z)+ \frac{rz}{1-z}.
\end{align*}
Alors, nous obtenons  
\begin{equation*}
A(z)=\frac{a_0}{1-z}+\frac{rz}{(1-z)^2}.
\end{equation*}}
\eex
\bex
 {\rm La f.g.o. de la suite double $(\binom{n}{k})_{n,k}$ connue sous le nom de \textit{triangle de Pascal} est la série bivariée
 \begin{align*}
\sum_{n,k\geq 0} \binom{n}{k}q^kz^n=\frac{1}{1-(1+q)z}.
 \end{align*}} 
\eex

\bex
{\rm Randrianarivony \cite{ARandr} a exprimé la f.g.o. du polynôme distributeur du nombre de croisements sur $S_n(321)$ en termes de fraction continue et il a prouvé l'identité suivante:
\begin{equation*}
\sum_{n\geq 0}\left( \sum_{\sigma \in S_n(321)} q^{\crs(\sigma)}\right) z^n=\frac{1}{1-\displaystyle\frac{z}{				
		1-\displaystyle\frac{z}{
			1-\displaystyle\frac{qz}{								
				1-\displaystyle\frac{qz}{
					1-\displaystyle\frac{q^2z}{
						1-\displaystyle\frac{q^2z}{				
							\ddots}
					}
			}}		
}}}.
\end{equation*}
 }
\eex
Voici une proposition fondamentale à laquelle nous ferons toujours référence dans les chapitres qui suivent. La preuve est laissée au lecteur.
\bpro \label{profgo}
Soit $(a_n)$ et $(b_n)$ deux suites de f.g.o respectives $A(z)$ et $B(z)$. Soit $c_n=a_n+b_n$ et $d_n=\sum_{k=0}^{n}a_{k}b_{n-k}$ pour tout $n\geq 0$. Les f.g.os de $(c_n)$ et $(d_n)$ sont respectivement
\begin{align*}
C(z)=A(z)+B(z) \text{ et } D(z)=A(z).B(z).
\end{align*}
\epro
\subsection{Preuve bijective d'une identité}

\bd
{\rm Soit $A$, $B$ deux ensembles. Soit $s_1$, $s_2$ deux statistiques respectives sur $A$ et  $B$.\\
\begin{itemize}
\item \textit{Prouver bijectivement} l'identité  $|A|=|B|$  consiste à trouver une bijection entre $A$ et $B$,
\item \textit{Prouver bijectivement} l'identité  $ \sum_{a \in A} q^{s_1(a)}= \sum_{b \in B} q^{s_2(b)}$  consiste à trouver une bijection $f:A\rightarrow B$ qui échange les statistiques $s_1$ sur $A$ et $s_2$ sur $B$, c'est-à-dire, $s_2(f(a))=s_1(a)$ pour tout $a \in A$.
\end{itemize}
}
\ed

\bex
{\rm MacMahon \cite{MacM} a montré algébriquement les identités de \eqref{eul-mah} et Foata \cite{Foata} en a donné les preuves bijectives.}
\eex
\bex
{\rm On peut trouver les preuves bijectives des identités suivantes dans \cite{ElizP,Eliz1,Eliz2}.
\begin{equation*}\label{fpex}
 \sum_{\sigma \in S_n(321)}x^{\fp(\sigma)}y^{\exc(\sigma)}= \sum_{\sigma \in S_n(132)}x^{\fp(\sigma)}y^{\exc(\sigma)}= \sum_{\sigma \in S_n(213)}x^{\fp(\sigma)}y^{\exc(\sigma)}.
 \end{equation*}
}
\eex 
 \bex
{ \rm Dans \cite{Cort,ARandr}, on peut trouver deux versions différentes des preuves bijectives de l'identité  
 \begin{equation*}\label{sym}
 \sum_{\sigma \in S_n}x^{\crs(\sigma)}y^{\nes(\sigma)}=\sum_{\sigma \in S_n}y^{\crs(\sigma)}x^{\nes(\sigma)}.
 \end{equation*}}
 \eex

\bex
{\rm Dans le Chapitre \ref{chap2} de ce manuscrit, on va établir les preuves bijectives des identités suivantes: 
\begin{align*}
\sum_{\sigma \in S_n(321)} q^{\crs(\sigma)}=\sum_{\sigma \in S_n(132)} q^{\crs(\sigma)}=\sum_{\sigma \in S_n(213)} q^{\crs(\sigma)}.
\end{align*}
}
\eex

\subsection{Interprétation combinatoire d'une suite}
Dans cette section, on va parler d'une notion fondamentale en combinatoire énumérative.
\bd
{\rm \textit{Interpréter combinatoirement} une suite consiste à trouver une famille d'objets combinatoires énumérée par cette suite.}
\ed
 
 \bex {
 \rm Supposons que $(T_{n,k})_{n,k}$ soit un triangle (suite double). Un exemple d'interprétation combinatoire de $(T_{n,k})_{n,k}$ est
\begin{align*}
|\{a\in E_n|s(a)=k\}|=T_{n,k} \text{ pour tous entiers $n$ et $k\geq 0$},
\end{align*} 
où $E_n$ est une famille d'objets combinatoires dépendant de $n$ et $s$ une statistique sur $E_n$. On dit que le triangle $(T_{n,k})_{n,k}$ énumère la famille  $E_n$ selon la statistique $s$. En d'autres termes, si $T(x,y)=\sum_{n,k\geq 0}T_{n,k}x^k y^n$ désigne la fonction génératrice de $(T_{n,k})_{n,k}$, alors
\begin{align*}
\sum_{n\geq 0} \sum_{a \in E_n} q^{s(a)} z^n=T(q,z).
\end{align*}
}
\eex
\bex
{\rm Considérons par exemple les triangles \href{https://oeis.org/A076791}{A076791} et \href{https://oeis.org/A299927}{A299927} de OEIS \cite{OEIS} (ou On-line Encyclopedia of Integer Sequences) ayant les fonctions génératrices ordinaires respectives $A076791(q,z)=\frac{1-qz}{1-(1+q)z-(1-q)z^2}$ et $A299927(q,z)=1+\frac{(1-qz)z}{(1-z)(1-(1+q)z)}$. Dans le Chapitre  \ref{chap3} de cette thèse, on va prouver que ces triangles énumèrent les permutations interdisant certaines paires de motifs de $S_3$ selon le nombre de croisements. Plus précisément, on a
	\begin{align*}
		\sum_{\sigma \in S(231,321)}q^{\crs(\sigma)}z^{|\sigma|}&=A076791(q,z), \\	
			\sum_{\sigma \in S(123,\tau)}q^{\crs(\sigma)}z^{|\sigma|}&=A299927(q,z), \text{ pour tout $\tau \in \{132,213\}$}.
	\end{align*}}
\eex

\section{La bijection de Elizalde et Pak}
Il est bien connu que, pour tout $\sigma \in S_3$, nous avons $|S_n(\tau)|=\frac{1}{n+1}\binom{2n}{n}$ (le $n$-ème nombre de Catalan). La preuve bijective de l'identité $|S_n(321)|=|S_n(132)|$ intéresse beaucoup de combinatoristes \cite{Knuth}. C'est ainsi que diverses bijections entre les deux familles d'objets $S_n(321)$ et $S_n(132)$ sont construites \cite{ClaKitaev} et chacune d'elles a sa propre propriété. Dans cette section, nous nous sommes particulièrement intéressés à l'une de ces bijections. 

 Dans \cite{ElizP}, Elizalde et Pak ont construit une  bijection   $\Theta: S_n(321)\rightarrow S_n(132)$ conservant à la fois le nombre de points fixes et le nombre d'excédances, c'est-à-dire,  $(\fp,\exc)(\Theta(\sigma))=(\fp,\exc)(\sigma)$ pour tout $\sigma \in S_n(321)$. La formulation de la bijection $\Theta$ n'était pas directe car elle passe par deux  familles d'objets intermédiaires, notamment les tableaux de Young et les chemins de Dyck. 
 
 Dans cette section,  nous rappellerons la bijection $\Theta$ de Elizalde et Pak avant de proposer une nouvelle formulation qui n'utilise plus aucune famille d'objets auxiliaires. 
 
 \subsection{Quelques définitions, notations et algorithme}
 \subsubsection*{Diagramme et tableaux de Young}
 \bd
{\rm Un \textit{diagramme de Young}   est un arrangement de cases juxtaposées et alignées à gauche tel que les longueurs des lignes (de haut vers le bas) décroissent au sens large. La suite des longueurs des lignes d'un diagramme de Young forme une partition $\lambda=(\lambda_1,\ldots,\lambda_k)$ d'un entier $n$ (c'est-à-dire, $n=\sum \lambda_i$) tel que $\lambda_1\geq \ldots\geq \lambda_k\geq 1$, $n$ étant le nombre de cases du diagramme. La partition $\lambda$  est appelée \textit{forme} du diagramme.}
 \ed 
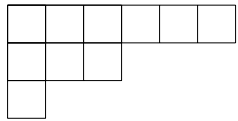
\begin{figure}[h]
	\begin{center}

		\begin{tikzpicture}			
		\draw[step=0.5cm, black, very thin] (0,0.5) grid (3,1); 
		\draw[step=0.5cm, black, very thin] (0,0) grid (1.5,1);
		\draw[step=0.5cm, black, very thin] (0,-0.5) grid (0.5,1);
		\end{tikzpicture}
		\label{fig1x1}
			\caption{Un diagramme de Young de forme $\lambda=(6, 3, 1)$.}
	\end{center}
\end{figure}
\bd 
{\rm Un \textit{tableau de Young} est un diagramme de Young rempli d'entiers strictement positifs selon les conditions suivantes 
\begin{itemize}
	\item[$\bullet$] les nombres inscrits dans les cases croissent de gauche à droite,
	\item[$\bullet$]  les nombres inscrits dans les cases croissent strictement de haut en bas.
\end{itemize}}
 \ed

\begin{figure}[h]
	\begin{center}
	
		\begin{tikzpicture}			
		\draw[step=0.5cm, black, very thin] (0,0.5) grid (3,1); 
		\draw[step=0.5cm, black, very thin] (0,0) grid (1.5,1);
		\draw[step=0.5cm, black, very thin] (0,-0.5) grid (0.5,1);
		
		\draw[black] (0.3,0.8) node {$1$}; \draw[black] (0.8,0.8) node {$3$}; \draw[black] (1.3,0.8) node {$5$};\draw[black] (1.8,0.8) node {$6$};\draw[black] (2.3,0.8) node {$8$}; \draw[black] (2.75,0.8) node {$10$};
		
		\draw[black] (0.3,0.3) node {$2$}; \draw[black] (0.8,0.3) node {$7$}; \draw[black] (1.3,0.3) node {$9$};
		
		\draw[black] (0.3,-0.2) node {$4$};	
		\end{tikzpicture}
		\label{fig1x2}
			\caption{Un tableau de Young de forme $\lambda=(6, 3, 1)$ rempli par $[10]$.}
	\end{center}
\end{figure}
\subsubsection*{Chemin de Dyck}
\bdf
	{\rm Un \textit{chemin de Dyck de demi-longueur $n$} est une suite de sommets $(x_i,y_i)_{0\leq i\leq 2n}$ du premier quadrant qui part de l'origine $(0,0)$ et se termine au point $(2n, 0)$ telle que $(x_{i+1},y_{i+1})=(x_i+1,y_i\pm 1)$ pour tout $0\leq i<2n$. 
		
		Un \textit{pas} d'un chemin de Dyck est un segment $(x,y)-(x',y')$ de $\mathbb{N}\times \mathbb{N}$ tel que $(x',y')=(x+1,y\pm 1)$. Le pas est dit \textit{ascendant} (resp. \textit{descendant}) si $(x',y')=(x+1,y+1)$ (resp. $(x',y')=(x+1,y-1)$).}
\edf
Habituellement, on code chaque pas ascendant par la lettre $u$ et chaque pas descendant par $d$. On désigne par $\mathcal{D}_{n}$ l'ensemble des chemins de Dyck de demi-longueur $n$. 	Il est à noter que, si $D=D_1D_2\ldots D_{2n}\in \D_n$, alors nous avons $|D|_u=|D|_d=n$ et $|D(k)|_u\geq |D(k)|_d$ pour tout sous-mot initial $D(k)=D_1D_2\ldots D_{k}$ de longueur $k$ de $D$, où $|w|_a$ désigne le nombre d'occurrences de la lettre $a$ dans un mot $w$.
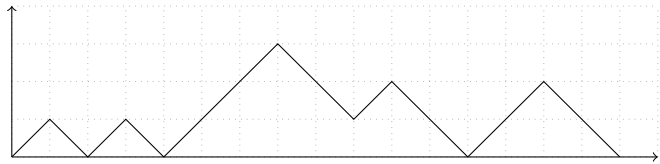
\begin{figure}[h]
	\begin{center}
	
		\begin{tikzpicture}
		\draw[step=0.5cm, gray, very thin,dotted] (0, 0) grid (8.5,2);
		\draw[black] (0, 0)--(0.5, 0.5)--(1, 0)--(1.5, 0.5)--(2, 0)--(2.5, 0.5)--(3, 1)--(3.5, 1.5)--(4, 1)--(4.5, 0.5)--(5, 1)--(5.5, 0.5)--(6, 0)--(6.5, 0.5)--(7, 1)--(7.5, 0.5)--(8, 0);
		\draw[->,black] (0, 0)--(0,2);	\draw[->,black] (0, 0)--(8.5,0);
		\end{tikzpicture}
		\caption{Un chemin de Dyck $D=ududuuuddudduudd$.}\label{fig11}	
	\end{center}
\end{figure} 

 Dans \cite{ElizD,ElizP},  Elizalde et al. ont introduit la statistique nombre de tunnels sur l'ensemble des chemins de Dyck. Ils ont utilisé cette nouvelle statistique pour énumérer les permutations à motif interdit selon les deux statistiques $\fp$ et $\exc$. 
 
 \bdf
 {\rm Soit $D\in \mathcal{D}_{n}$. 
 
 Un \textit{tunnel} de $D$ est un segment horizontal qui n'intersecte  $D$ que seulement en deux points et reste toujours en dessous de $D$. 
 
Un tunnel est dit \textit{gauche} (resp. \textit{centré}, \textit{droite})  si l'abscisse du milieu d'un tunnel de $D$ est inférieure (resp. égale, supérieure) à $n$. }
\edf
Graphiquement, les tunnels  droites et  gauches d'un chemin de Dyck sont respectivement ceux dont le milieu se trouve à droite et à gauche de la droite verticale $x=n$ séparant le chemin $D$ en deux. Nous désignons respectivement par $\rt(D)$, $\ct(D)$ et  $\lt(D)$ le nombre de tunnels droites, le nombre de tunnels centrés et le nombre de tunnels gauches. Le chemin de Dyck $D$ de la Figure \ref{tunnel}   possède les caractéristiques suivantes: $\rt(D)=3$, $\ct(D)=1$ et $\lt(D)=4$. 
	
	\begin{figure}[h]
		\begin{center}
			
			\begin{tikzpicture}
			\draw[step=0.5cm, gray, very thin,dotted] (0, 0) grid (8,1.5);
			\draw [-, black](4, -0.2)--(4, 1.5);
			\draw[black] (0, 0)--(0.5, 0.5)--(1, 0)--(1.5, 0.5)--(2, 0)--(2.5, 0.5)--(3, 1)--(3.5, 1.5)--(4, 1)--(4.5, 0.5)--(5, 1)--(5.5, 0.5)--(6, 0)--(6.5, 0.5)--(7, 1)--(7.5, 0.5)--(8, 0);	
			\draw[gray] (0.5,0) node {$\bullet$};  \draw [-,thin, gray](0,0)--(1,0);
			\draw[gray] (1.5,0) node {$\bullet$};  \draw [-,thin, gray](1,0)--(2,0);
			\draw[gray] (3.5,1) node {$\bullet$};  \draw [-,thin, gray](3,1)--(4,1);
			\draw[gray] (3.5,0.5) node {$\bullet$};  \draw [-,thin, gray](2.5,0.5)--(4.5,0.5);
			\draw[blue] (4,0) node {$\bullet$};  \draw [-,thin, blue](2,0)--(6,0);
			\draw[black] (5,0.5) node {$\bullet$};  \draw [-,thin, black](4.5,0.5)--(5.5,0.5);
			\draw[black] (7,0.5) node {$\bullet$};  \draw [-,thin, black](6.5,0.5)--(7.5,0.5);
			\draw[black] (7,0) node {$\bullet$};  \draw [-,thin, black](6,0)--(8,0);
			\end{tikzpicture}
		\caption{Les tunnels du chemin de Dyck $D=ududuuuddudduudd$.}\label{tunnel}	
		\end{center}
	\end{figure} 
	
	Comme illustré dans la Figure \ref{tunnel}, chaque tunnel est un segment dont l'extrémité gauche est l'origine d'un pas ascendant et l'extrémité droite est la destination d'un pas descendant. Il est mentionné dans \cite{ElizD} que chaque tunnel d'un chemin de Dyck $D$ est en correspondance avec la décomposition du chemin $D = AuBdC$, où $B$ est le sous-chemin de Dyck non vide au-dessus du tunnel et $AC$ est aussi un chemin de Dyck (mais pouvant être vide).	

\subsubsection*{Notations}
Étant donné une permutation $\sigma$ de $[n]$, nous adoptons les notations suivantes:
	\begin{itemize}
		\item[$\bullet$] pour tous entiers naturels $a$ et $b$ tels que $a\leq b$, $\sigma(a\cdots b) :=\sigma(a)\sigma(a+1)\cdots\sigma(b)$ désigne la sous-séquence contiguë de $\sigma$ de la $a$-ème jusqu'à la $b$-ème lettre et pour tout opérateur $* \in \{<, \leq, \neq, \geq, >\}$ et un nombre $x$ donné, nous pouvons écrire $\sigma(a\cdots b)* x$ si et seulement si $\sigma(i)* x$ pour tout $i\in [a;b]$.  
		
		Exemple : si nous considérons la permutation $\pi=6413275$, alors nous avons $\pi(2\cdots 5)=4132\leq 4$.
		\item[$\bullet$] pour tout entier naturel non nul $a$,  $\sigma^{+a}$ désigne la permutation de $\{1+a,\ldots, n+a\}$ déduite de  $\sigma$ en ajoutant $a$ à chacun de ses nombres. Exemple: $312^{+2}=534$.
		\item[$\bullet$] pour tous entiers naturels non nuls $a$ et $b$, $\sigma^{a \rtimes b}$ désigne la permutation  de $\{1,\ldots, a-1,a+b,\ldots,n+b\}$ déduite de $\sigma$ en ajoutant $b$ à tous les nombres de $\sigma$ qui sont plus grands ou égaux à $a$. 
		
		Exemple: $4132^{3\rtimes 2}=\textcolor{gray}{6}1\textcolor{gray}{5}2$.
		\item[$\bullet$] pour tous entiers  $i$ et $x \in [n+1]$, $\sigma^{(i,x)}$  désigne la permutation de $[n+1]$ déduite de $\sigma$ en  augmentant d'une unité tous les entiers supérieurs ou égaux à $x$ dans $\sigma$, puis en insérant l'entier $x$ à la $i$-ème position. Autrement dit, nous avons $\sigma^{(i,x)}:=\sigma^{x\rtimes 1}(1\cdots i-1)\cdot x \cdot \sigma^{x\rtimes 1}(i\cdots |\sigma|)$.

		Exemple: $3142^{(2,\textcolor{gray}{3})}=4\textcolor{gray}{3}152$. 
	\end{itemize}

\subsubsection*{Un algorithme}
Notons $\mathcal{M}(\sigma)$ l'ensemble des couples (valeur d'excédance, non-excédance) obtenu par l'Algorithme \ref{algo} ci-dessous.  Par exemple, nous avons $\mathcal{M}(43152)=\{(4,3),(3,5)\}$.
	 
\begin{algorithm}[ht!]
	\KwData{A permutation $\sigma$ 	with excedences $e_1<\ldots<e_k$ and non-excedences $a_1<\ldots<a_{n-k}$.}
	\KwResult{A set of pairs (excedence value, non-excedence)  of $\sigma$ denoted by $\mathcal{M}(\sigma)$.}
	
	Let  $p \leftarrow 1$, $q \leftarrow 1$ and $\mathcal{M}(\sigma)\leftarrow\{\}$\;
	\Begin{
		\Repeat{$p>k \vee q>n-k$}{
			\eIf{$e_p>a_q$}{$q\leftarrow q+1$\;}{
				\eIf{$\sigma(e_p)<\sigma(a_q)$}{ $p \leftarrow p + 1$\;}{ 
					$\mathcal{M}(\sigma)\leftarrow\mathcal{M}(\sigma)\cup \{(\sigma(e_p),a_q)\}$\;
					$p \leftarrow p + 1$\; 
					$q \leftarrow q + 1$\;
				}
			}
		}	
	}
	\caption{\og Matching algorithm\fg.}
	\label{algo}
\end{algorithm}

Il est à noter que cet algorithme, qui associe certaines non-excédances aux valeurs d'excédances d'une permutation donnée, n'est autre que l'algorithme d'affectation présenté par Elizalde et Pak dans \cite{ElizP} mais nous avons seulement modifié sa sortie pour une raison particulière correspondant à notre problème. 
	
	\subsection{La bijection $\Theta:S_n(321)\rightarrow S_n(132)$}\label{sec12}
La bijection $\Theta$ est définie originalement comme la composition de deux correspondances bien connues dans la littérature. La première est  une bijection $\Psi: S_n(321)\rightarrow\D_n$ due à Knuth \cite{Knuth} qui utilise les tableaux de Young comme objet auxiliaire. La deuxième est  une bijection $\Phi:S_n(132)\rightarrow\D_n$ qui peut être obtenue par une simple modification de Krathenthaler \cite{Krat}.

	\subsubsection{La bijection $\Psi:S_n(321)\rightarrow \mathcal{D}_{n}$}\label{sec112}	
	La bijection $\Psi$ est essentiellement due à Knuth \cite{Knuth} et est  une composée de deux bijections que nous allons décrire dans la suite. 
	
La première est une bijection entre $S_n$ et l'ensemble des couples de tableaux de Young de même forme $\lambda \vdash n$. La construction de l'image est basée sur l'algorithme d'insertion  de Robinson-Schensted-Knuth ou simplement l'algorithme de $\RSK$ (voir \cite{Knuth,STAN}). Soit $\sigma \in S_n$ et $(P,Q)=\RSK(\sigma)$. Assumons que nous avons $(P^{(i-1)},Q^{(i-1)})=\RSK(\sigma(1\cdots i-1))$. Nous obtenons $(P^{(i)},Q^{(i)})$ en insérant $\sigma(i)$ dans $P^{(i-1)}$ et $i$ dans $Q^{(i-1)}$ de la manière suivante:  
\begin{itemize}
	\item[-] si $\sigma(i)$ est plus grand que tous les nombres de la première ligne de $P^{(i-1)}$, alors on crée une nouvelle case  à la fin de la première ligne de $P^{(i-1)}$ et on y place $\sigma(i)$. Sinon, $\sigma(i)$ prend la place de l'élément $x$  plus grand que $\sigma(i)$ qui est le plus à gauche (on dit que $\sigma(i)$ \textit{bouscule} $x$),  puis on insère $x$ dans la ligne suivante de la même façon. 
	\item[-] puisque le tableau $Q^{(i)}$ a la même forme que $P^{(i)}$, il suffit de créer une nouvelle case pour $i$ à la même position que la case nouvellement créée dans $P^{(i)}$ lors de l'insertion de $\sigma(i)$.
\end{itemize}	
	  De cette façon, nous avons $(P,Q)=(P^{(n)},Q^{(n)})$. Une des fameuses propriétés bien connues de la correspondance de $\RSK$ est la notion de dualité. Autrement dit, nous avons $\RSK(\sigma^{-1})=(Q,P)$ si et seulement si $\RSK(\sigma)=(P,Q)$ (voir \cite{Knuth}). Une autre propriété de $\RSK$  est que le nombre de lignes de $P$ est égale à la longueur de la plus longue séquence croissante dans $\sigma$. Par conséquent,  la permutation $\sigma$ est 321-interdite si et seulement si $P$ (et évidement $Q$) possède au plus deux lignes. Comme illustré dans la figure \ref{fig:rsk}, nous détaillons la construction de $(P,Q)$ à partir de la  permutation $\sigma=24135867 \in S_8(321)$.
	
	\begin{figure}[h]
		\begin{center}		
		
			\begin{tikzpicture}
			\draw (0,0.6)-- (0,1);\draw (0.4,0.6)-- (0.4,1); \draw (0,0.6)-- (0.4,0.6);\draw (0,1)-- (0.4,1);
			\draw[step=0.4cm, black, very thin] (0, -0.4) grid (0.4,0);
			\draw[black] (0.2,0.8) node {$2$};
			\draw[black] (0.2,-0.2) node {$1$};
			
			\draw (0.8, 0.6) -- (1.6,0.6); \draw (0.8, 1) -- (1.6,1); \draw (0.8, 0.6) -- (0.8,1); \draw (1.6, 0.6) -- (1.6,1); 
			\draw (1.6, 0.6) -- (1.6,1); \draw (1.2, 0.6) -- (1.2,1);
			\draw (0.8, -0.4) -- (1.6,-0.4);\draw (0.8, 0) -- (1.6,0);\draw (0.8, -0.4) -- (0.8,0); \draw (1.6, -0.4) -- (1.6,0); 
			\draw (1.6, -0.4) -- (1.6,0); \draw (1.2, -0.4) -- (1.2,0);
			\draw[black] (1,0.8) node {$2$};\draw[black] (1.4,0.8) node {$4$};
			\draw[black] (1,-0.2) node {$1$}; \draw[black] (1.4,-0.2) node {$2$};
			
			\draw (2, 0.6) -- (2.8,0.6);\draw (2, 1) -- (2.8,1);\draw (2, 0.2) -- (2,1); \draw (2.4, 0.2) -- (2.4,1); 
			\draw (2.8, 0.6) -- (2.8,1); \draw (2, 0.2) -- (2.4,0.2);
			\draw (2, -0.4) -- (2.8,-0.4);\draw (2, 0) -- (2.8,0);\draw (2, -0.8) -- (2,0); \draw (2.4, -0.8) -- (2.4,0); 
			\draw (2.8, -0.4) -- (2.8,0); \draw (2, -0.8) -- (2.4,-0.8);
			\draw[black] (2.2,0.8) node {$1$};\draw[black] (2.6,0.8) node {$4$};
			\draw[black] (2.2,0.4) node {$2$};
			\draw[black] (2.2,-0.2) node {$1$}; \draw[black] (2.6,-0.2) node {$2$};
			\draw[black] (2.2,-0.6) node {$3$};			
			
			\draw (3.2, 0.6) -- (4,0.6);\draw (3.2, 1) -- (4,1);\draw (3.2, 0.2) -- (3.2,1); \draw (3.6, 0.2) -- (3.6,1); 
			\draw (4, 0.2) -- (4,1); \draw (3.2, 0.2) -- (4,0.2);
			\draw (3.2, -0.4) -- (4,-0.4);\draw (3.2, 0) -- (4,0);\draw (3.2, -0.8) -- (3.2,0); \draw (3.6, -0.8) -- (3.6,0); 
			\draw (4, -0.8) -- (4,0); \draw (3.2, -0.8) -- (4,-0.8);
			\draw[black] (3.4,0.8) node {$1$};  \draw[black] (3.8,0.8) node {$3$};
			\draw[black] (3.4,0.4) node {$2$};  \draw[black] (3.8,0.4) node {$4$};
			\draw[black] (3.4,-0.2) node {$1$}; \draw[black] (3.8,-0.2) node {$2$};
			\draw[black] (3.4,-0.6) node {$3$}; \draw[black] (3.8,-0.6) node {$4$};
			
			\draw (4.6, 0.6) -- (5.8,0.6);\draw (4.6, 1) -- (5.8,1);\draw (4.6, 0.2) -- (5.4,0.2);
			\draw (4.6, 0.2) -- (4.6,1); \draw (5, 0.2) -- (5,1);  \draw (5.4, 0.2) -- (5.4,1);  \draw (5.8, 0.6) -- (5.8,1);
			\draw (4.6, -0.4) -- (5.8,-0.4);\draw (4.6, 0) -- (5.8,0);  \draw (4.6, -0.8) -- (5.4,-0.8);
			\draw (4.6, -0.8) -- (4.6,0); \draw (5, -0.8) -- (5,0); 	  \draw (5.4, -0.8) -- (5.4,0); \draw (5.8, -0.4) -- (5.8,0);
			\draw[black] (4.8,0.8) node {$1$};  \draw[black] (5.2,0.8) node {$3$}; \draw[black] (5.6,0.8) node {$5$};
			\draw[black] (4.8,0.4) node {$2$};  \draw[black] (5.2,0.4) node {$4$};
			\draw[black] (4.8,-0.2) node {$1$}; \draw[black] (5.2,-0.2) node {$2$}; \draw[black] (5.6,-0.2) node {$5$};
			\draw[black] (4.8,-0.6) node {$3$}; \draw[black] (5.2,-0.6) node {$4$};			
			
			\draw (6.2, 0.6) -- (7.8,0.6);\draw (6.2, 1) -- (7.8,1);\draw (6.2, 0.2) -- (7,0.2); 
			\draw (6.2, 0.2) -- (6.2,1); \draw (6.6, 0.2) -- (6.6,1);  \draw (7, 0.2) -- (7,1); \draw (7.4, 0.6) -- (7.4,1); \draw (7.8, 0.6) -- (7.8,1);
			\draw (6.2, -0.4) -- (7.8,-0.4);\draw (6.2, 0) -- (7.8,0);  \draw (6.2, -0.8) -- (7,-0.8);
			\draw (6.2, -0.8) -- (6.2,0); \draw (6.6, -0.8) -- (6.6,0); 	  \draw (7, -0.8) -- (7,0);\draw (7.4, -0.4) -- (7.4,0); \draw (7.8, -0.4) -- (7.8,0);
			\draw[black] (6.4,0.8) node {$1$};  \draw[black] (6.8,0.8) node {$3$}; \draw[black] (7.2,0.8) node {$5$}; \draw[black] (7.6,0.8) node {$8$};
			\draw[black] (6.4,0.4) node {$2$};  \draw[black] (6.8,0.4) node {$4$};
			\draw[black] (6.4,-0.2) node {$1$}; \draw[black] (6.8,-0.2) node {$2$}; \draw[black] (7.2,-0.2) node {$5$}; \draw[black] (7.6,-0.2) node {$6$};
			\draw[black] (6.4,-0.6) node {$3$}; \draw[black] (6.8,-0.6) node {$4$};

			\draw (8.4, 0.6) -- (10,0.6); \draw (8.4, 1) -- (10,1);\draw (8.4, 0.2) -- (9.6,0.2); 
			\draw (8.4, 0.2) -- (8.4,1); \draw (8.8, 0.2) -- (8.8,1);  \draw (9.2, 0.2) -- (9.2,1); \draw (9.6, 0.2) -- (9.6,1); \draw (10, 0.6) -- (10,1);	  	  
			\draw (8.4, -0.4) -- (10,-0.4);\draw (8.4, 0) -- (10,0);  \draw (8.4, -0.8) -- (9.6,-0.8);
			\draw (8.4, -0.8) -- (8.4,0); \draw (8.8, -0.8) -- (8.8,0); 	\draw (9.2, -0.8) -- (9.2,0);  \draw (9.6, -0.8) -- (9.6,0);\draw (9.6, -0.4) -- (9.6,0); \draw (10, -0.4) -- (10,0);
			
			\draw[black] (8.6,0.8) node {$1$};  \draw[black] (9,0.8) node {$3$}; \draw[black] (9.4,0.8) node {$5$}; \draw[black] (9.8,0.8) node {$6$};
			\draw[black] (8.6,0.4) node {$2$};  \draw[black] (9,0.4) node {$4$}; \draw[black] (9.4,0.4) node {$8$};
			\draw[black] (8.6,-0.2) node {$1$}; \draw[black] (9,-0.2) node {$2$}; \draw[black] (9.4,-0.2) node {$5$}; \draw[black] (9.8,-0.2) node {$6$};
			\draw[black] (8.6,-0.6) node {$3$}; \draw[black] (9,-0.6) node {$4$};\draw[black] (9.4,-0.6) node {$7$};

			\draw (10.6, 0.6) -- (12.6,0.6); \draw (10.6, 1) -- (12.6,1);\draw (10.6, 0.2) -- (11.8,0.2); 
			\draw (10.6, 0.2) -- (10.6,1); \draw (11, 0.2) -- (11,1);  \draw (11.4, 0.2) -- (11.4,1); \draw (11.8, 0.2) -- (11.8,1); \draw (12.2, 0.6) -- (12.2,1);	\draw (12.6, 0.6) -- (12.6,1);  	  
			\draw (10.6, -0.4) -- (12.6,-0.4);\draw (10.6, 0) -- (12.6,0);  \draw (10.6, -0.8) -- (11.8,-0.8);
			\draw (10.6, -0.8) -- (10.6,0); \draw (11, -0.8) -- (11,0); 	\draw (11.4, -0.8) -- (11.4,0);  \draw (11.8, -0.8) -- (11.8,0);\draw (12.2, -0.4) -- (12.2,0); \draw (12.6, -0.4) -- (12.6,0);
			
			\draw[black] (10.8,0.8) node {$1$};  \draw[black] (11.2,0.8) node {$3$}; \draw[black] (11.6,0.8) node {$5$}; \draw[black] (12,0.8) node {$6$}; \draw[black] (12.4,0.8) node {$7$};
			\draw[black] (10.8,0.4) node {$2$};  \draw[black] (11.2,0.4) node {$4$}; \draw[black] (11.6,0.4) node {$8$};
			\draw[black] (10.8,-0.2) node {$1$}; \draw[black] (11.2,-0.2) node {$2$}; \draw[black] (11.6,-0.2) node {$5$}; \draw[black] (12,-0.2) node {$6$}; \draw[black] (12.4,-0.2) node {$8$};
			\draw[black] (10.8,-0.6) node {$3$}; \draw[black] (11.2,-0.6) node {$4$};\draw[black] (11.6,-0.6) node {$7$};
			
			\draw[black] (13.2,0.6) node {$=P$};
			\draw[black] (13.2,-0.4) node {$=Q$};  	 	
			\end{tikzpicture}
				\caption{Construction de $(P,Q)=\RSK(24135867)$}\label{fig:rsk}
		\end{center}
	\end{figure}

\hspace{-0.59cm}Pour toute permutation 321-interdite $\sigma$, nous remarquons que l'algorithme \ref{algo}  nous permet d'obtenir facilement $(P,Q)=\RSK(\sigma)$. De plus, si  $\mathcal{M}(\sigma)=\{(E_1,a_1),(E_2,a_2),\\ \ldots,(E_l,a_l)\}$, alors les deuxièmes lignes de $P$ et $Q$ sont respectivement $[E_1,E_2,\ldots, E_l]$  et $[a_1,a_2,\ldots, a_l]$ et on peut immédiatement en déduire les premières lignes.
Par exemple, comme  $\mathcal{M}(\sigma)=\{(2,\textcolor{gray}{3}), (4,\textcolor{gray}{4}), (8,\textcolor{gray}{7})\}$ pour $\sigma=24135867 \in S_8(321)$, alors, si  $(P,Q)=\RSK(\sigma)$, les deuxièmes lignes de $P$ et $Q$ sont respectivement  $[2,4,8]$ et $[3,4,7]$. Ainsi, les premières lignes $P$ et $Q$ sont respectivement $[1,3,5,6,7]$ et $[1,2,5,6,8]$.

La deuxième correspondance est une simple transformation de la paire de tableaux de Young standards $(P,Q)$, résultat de $\RSK$,  en chemin de Dyck $D=\Psi(\sigma)$. Il est d'abord à noter que $P$ (resp. $Q$) possède au plus deux lignes car $\sigma$ est 321-interdite. On construit le chemin $D$ de la façon suivante:
\begin{itemize}
	\item[-] la moitié gauche du chemin $D$ est déduite de $P$ en joignant, pour $i$ de $1$ à $n$, un pas montant si le nombre $i$ se trouve dans la première ligne de $P$ et un pas descendant si le nombre $i$ se trouve dans la deuxième ligne de $P$. 
	\item[-] La moitié droite du chemin $D$ est alors déduite de $Q$ en joignant, pour $j$ de $n$ à $1$, un pas ascendant si le nombre $j$ se trouve dans la deuxième ligne de $Q$ et un pas descendant si le nombre $j$ se trouve dans la première ligne de $Q$. 
\end{itemize}
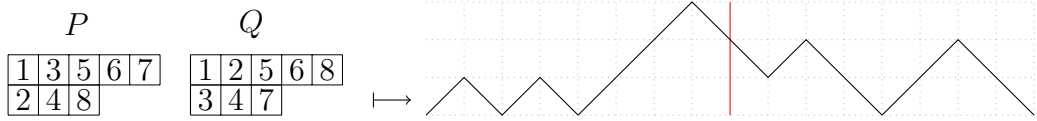
\begin{figure}[h]
		\begin{center}
			
			\begin{tikzpicture}
			\draw[black] (0.9,1.2) node {$P$}; \draw[black] (3.2,1.2) node {$Q$};
			\draw[step=0.4cm, black, very thin] (0,0.4) grid (2,0.8); \draw[step=0.4cm, black, very thin] (2.4,0.4) grid (4.4,0.8);
			\draw[step=0.4cm, black, very thin] (0, 0) grid (1.2,0.4);     \draw[step=0.4cm, black, very thin] (2.4, 0)  grid (3.6,0.4);
			\draw (2.4,0)-- (2.4,0.8);	
			\draw[black] (0.2,0.6) node {$1$};	\draw[black] (0.6,0.6) node {$3$};  \draw[black] (1,0.6) node {$5$}; \draw[black] (1.4,0.6) node {$6$}; \draw[black] (1.8,0.6) node {$7$};
			\draw[black] (0.2,0.2) node {$2$}; 	\draw[black] (0.6,0.2) node {$4$}; \draw[black] (1,0.2) node {$8$};
			
			\draw[black] (2.6,0.6) node {$1$};	\draw[black] (3,0.6) node {$2$}; \draw[black] (3.4,0.6) node {$5$}; \draw[black] (3.8,0.6) node {$6$}; \draw[black] (4.2,0.6) node {$8$};
			\draw[black] (2.6,0.2) node {$3$}; 	\draw[black] (3,0.2) node {$4$}; \draw[black] (3.4,0.2) node {$7$};
			
			\draw[|->,thin, black] (4.8, 0.2)--(5.3,0.2);
			\draw[-, red] (9.5, 0)--(9.5,1.5);
			\draw[step=0.5cm, gray, very thin,dotted] (5.5, 0) grid (13.5,1.5);
			\draw[black] (5.5, 0)--(6, 0.5)--(6.5, 0)--(7, 0.5)--(7.5, 0);
			\draw[black](7.5, 0)--(8, 0.5);
			\draw[black](8, 0.5)--(8.5, 1)--(9, 1.5)--(9.5, 1)--(10, 0.5)--(10.5, 1)--(11, 0.5);
			\draw[black](11, 0.5)--(11.5, 0);
			\draw[black](11.5, 0)--(12, 0.5)--(12.5, 1)--(13, 0.5)--(13.5, 0);
			\end{tikzpicture}
			\caption{Le chemin de Dyck correspondant à $(P,Q)={\rm RSK}(24135867)$.}
		\end{center}
	\end{figure}
Si nous notons respectivement $D^{(L)}$ et $D^{(R)}$ la moitié gauche et la moitié droite du chemin $D$, alors nous pouvons écrire $D=D^{(L)} D^{(R)}$. De plus, nous avons la proposition évidente suivante. 
	\bpro\label{prop22}
	Pour tout chemin de Dyck $D$, nous avons $|D^{(L)}|_d=|D^{(R)}|_u$.
	\epro
	\begin{proof}
	Pour tout chemin de Dyck $D$, il existe une permutation $\sigma \in S(321)$ telle que $D=\Psi(\sigma)$. Si $(P,Q)=\RSK(\sigma)$, alors, d'après la définition de la deuxième correspondance décrite ci-dessus, $|D^{(L)}|_d$ et $|D^{(R)}|_u$ sont respectivement égaux aux tailles des secondes lignes de $P$ et $Q$. Ainsi, nous avons $|D^{(L)}|_d=|D^{(R)}|_u$.
	\end{proof}
	\subsubsection{La bijection $\Phi^{-1}:\mathcal{D}_{n}\rightarrow S_n(132)$} \label{sec113}
La bijection $\Phi$ n'est autre qu'une simple touche de celui de Krattenthaler \cite{Krat}. Ici, nous allons présenter la bijection $\Phi^{-1}$ d'une manière légèrement différente de celle de Elizalde et Pak.

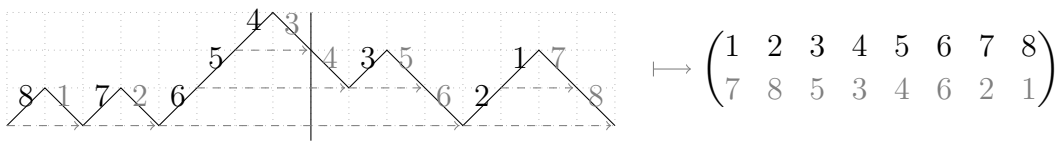
\begin{figure}[h]
	\begin{center}
		
		\begin{tikzpicture}
		\draw[step=0.5cm, gray, very thin,dotted] (0, 0) grid (8,1.5);
		\draw [-, black](4, -0.2)--(4, 1.5);
		\draw[black] (0, 0)--(0.5, 0.5)--(1, 0)--(1.5, 0.5)--(2, 0)--(2.5, 0.5)--(3, 1)--(3.5, 1.5)--(4, 1)--(4.5, 0.5)--(5, 1)--(5.5, 0.5)--(6, 0)--(6.5, 0.5)--(7, 1)--(7.5, 0.5)--(8, 0);
		\draw[black] (0.25,0.4) node {$8$};
		\draw[black] (1.25,0.4) node {$7$};
		\draw[black] (2.25,0.4) node {$6$}; 
		\draw[black] (2.75,0.9) node {$5$};
		\draw[black] (3.25,1.4) node {$4$}; 
		\draw[black] (4.75,0.9) node {$3$};
		\draw[black] (6.25,0.4) node {$2$}; 
		\draw[black] (6.75,0.9) node {$1$};   
		
		\draw[gray] (0.75,0.4) node {$1$};
		\draw[gray] (1.75,0.4) node {$2$};
		\draw[gray] (3.75,1.35) node{$3$};
		\draw[gray] (4.25,0.9) node {$4$};
		\draw[gray] (5.25,0.9) node {$5$}; 
		\draw[gray] (5.75,0.4) node {$6$};
		\draw[gray] (7.25,0.9) node {$7$}; 
		\draw[gray] (7.75,0.4) node {$8$}; 
		
		\draw [->,dashed, gray](0,0)--(0.95,0);
		\draw [->,dashed, gray](1,0)--(1.95,0);
		\draw [->,dashed, gray](3,1)--(3.95,1);
		\draw [->,dashed, gray](2.5,0.5)--(4.45,0.5);
		\draw [->,dashed, gray](2,0)--(5.95,0);
		\draw [->,dashed, gray](4.5,0.5)--(5.45,0.5);
		\draw [->,dashed, gray](6.5,0.5)--(7.45,0.5);
		\draw [->,dashed, gray](6,0)--(7.95,0);
		\draw[|->,thin, gray] (8.5, 0.75)--(9,0.75);
		\draw[black] (11.5,0.75) node {$\begin{pmatrix}
			\textcolor{black}{1}& \textcolor{black}{2}& \textcolor{black}{3}& \textcolor{black}{4}& \textcolor{black}{5}& \textcolor{black}{6}& \textcolor{black}{7}& \textcolor{black}{8} \\
			\textcolor{gray}{7}& \textcolor{gray}{8}& \textcolor{gray}{5}& \textcolor{gray}{3}& \textcolor{gray}{4}& \textcolor{gray}{6}& \textcolor{gray}{2}& \textcolor{gray}{1}
			\end{pmatrix}$};
		\end{tikzpicture}
		\caption{La permutation 132-interdite correspondante à un chemin de Dyck.}\label{fig:form1}
	\end{center}
\end{figure}

En partant d'un chemin de Dyck $D\in \mathcal{D}_{n}$, nous construisons la permutation correspondante par la procédure suivante.
De gauche à droite, attribuons des numéros aux pas ascendants de $D$ de $n$ à $1$ et aux pas descendants de $1$ à $n$.  Ensuite, nous avons $\Phi^{-1}(D)(n+1-i)=j$ si et seulement si l'origine du pas ascendant numéroté $n+1-i$ (i.e. le $i$-ième pas ascendant) et la destination du pas descendant numéroté $j$ (i.e.  le $j$-ième pas descendant) de $D$ sont les extrémités d'un tunnel de $D$. De cette façon, il n'est pas difficile de montrer que l'application $\Phi^{-1}$ est une bijection bien définie. Voir Figure \ref{fig:form1} pour une illustration graphique.

Notre première observation de la bijection $\Phi$ décrite par la procédure ci-dessus conduit à la proposition suivante.
\bpro\label{prop23}
	Supposons que $D\in \D_n$ et $j=|D^{(R)}|_u+1$. La permutation $\sigma=\Phi^{-1}(D)$ satisfait les propriétés suivantes:
	\begin{itemize}
		\setlength\itemsep{-0.3em}
		\item[{\rm (i)}] Si $j\geq 2$, alors nous avons $\sigma^{-1}(1\cdots j-1)\geq j\leq \sigma(1\cdots j-1)$,
		\item[{\rm (ii)}] Pour tout entier $i\geq j$, si $\sigma(i)>i$, alors nous avons $\sigma^{-1}(i)<i$.
	\end{itemize}
\epro
\begin{proof}
Considérons $D\in \D_n$ et supposons que $j=|D^{(R)}|_u+1\geq 2$. Nous résumons dans le Tableau \ref{table:tabx0} les numéros assignés aux pas ascendants et descendants de $D$ pour obtenir $\sigma=\Phi^{-1}(D)$. 

\begin{table}[h]
	\begin{center}
	
		\begin{tabular}{|c|c|c|c|}
			\hline
			Sous-chemin & $D^{(L)}$ & $D^{(R)}$\\
			\hline
			Pas ascendants & $n,\ldots,j+1,j$& $j-1,\ldots, 2,1$\\
			\hline
			Pas descendants &$1,2,\ldots,j-1$ &$j,j+1,\ldots,n$\\
			\hline
		\end{tabular}
		\label{table:tabx0}
		\caption{Les numéros assignés aux pas du chemin $D$.}	
	\end{center}
\end{table}

En regardant la colonne de $D^{(R)}$ du tableau, nous obtenons $\sigma(1\cdots j-1)\geq j$. De même, en regardant la  colonne de $D^{(L)}$, nous obtenons également $\sigma^{-1}(1\cdots j-1)\geq j$. Ceci prouve la propriété (i).

Considérons maintenant un entier $i$ tel que $j\leq i<\sigma(i)$. Dans le chemin de Dyck $D$, le tunnel liant le pas ascendant  numéroté $i$ avec le pas descendant numéroté $\sigma(i)$ décompose le chemin de Dyck $D = \cdots uBd \cdots$, où $u$ est le $(n+1-i)$-ème pas ascendant de $D$, $d$ est le $\sigma(i)$-ème pas descendant de $D$, $B$ est un sous-chemin de Dyck de $D$ qui contient au moins les pas descendants numérotés $j,j+1,\ldots, \sigma(i)-1$ et les pas ascendants de $B$ sont évidemment numérotés par des nombres inférieurs à $i$. Cela implique que nous devons avoir $\sigma^{-1}(j\cdots \sigma(i)-1)<i$. Puisque $i\in \{j,j+1,\cdots, \sigma(i)-1\}$, nous obtenons donc $\sigma^{-1}(i)<i$. Ceci termine également la preuve de la  propriété (ii).			
\end{proof}
Concluons simplement cette section par la remarque évidente suivante qui implique que $\Phi^{-1}$ échange les tunnels gauches et centrés d'un chemin de Dyck et les non-excédances de la permutation 132-interdite correspondante.
\brem\label{rem23}
{\rm Le tunnel d'un chemin de Dyck $D\in \D_n$ associant le pas ascendant numéroté $n+1-i$ au  pas descendant numéroté $j$ est un tunnel gauche ou centré si et seulement si $n+1-i\geq j=\Phi^{-1}(D)(n+1-i)$}.
\erem
\subsubsection{Illustration graphique}
\begin{figure}[h]
	 	\begin{center}
	 	
	 		\begin{tikzpicture}
	 		
	 		\draw[black] (2,4) node {$24135867 \in S_8(321)$}; 		
	 		\draw[|->,dashed, black] (2, 3.5)--(2,2);
	 		
	 		\draw[black] (0.9,1.2) node {$P$}; \draw[black] (3.2,1.2) node {$Q$};
	 		\draw[step=0.4cm, black, very thin] (0,0.4) grid (2,0.8); \draw[step=0.4cm, black, very thin] (2.4,0.4) grid (4.4,0.8);
	 		\draw[step=0.4cm, black, very thin] (0, 0) grid (1.2,0.4);     \draw[step=0.4cm, black, very thin] (2.4, 0)  grid (3.6,0.4);
	 		\draw (2.4,0)-- (2.4,0.8);	
	 		\draw[black] (0.2,0.6) node {$1$};	\draw[black] (0.6,0.6) node {$3$};  \draw[black] (1,0.6) node {$5$}; \draw[black] (1.4,0.6) node {$6$}; \draw[black] (1.8,0.6) node {$7$};
	 		\draw[black] (0.2,0.2) node {$2$}; 	\draw[black] (0.6,0.2) node {$4$}; \draw[black] (1,0.2) node {$8$};
	 		
	 		\draw[black] (2.6,0.6) node {$1$};	\draw[black] (3,0.6) node {$2$}; \draw[black] (3.4,0.6) node {$5$}; \draw[black] (3.8,0.6) node {$6$}; \draw[black] (4.2,0.6) node {$8$};
	 		\draw[black] (2.6,0.2) node {$3$}; 	\draw[black] (3,0.2) node {$4$}; \draw[black] (3.4,0.2) node {$7$};
	 		
	 		\draw[|->,dashed, black] (4.8, 0.2)--(5.3,0.2);
	 			\draw[|->,dashed, black] (4, 3.5)--(5.5,1.5);
	 				\draw[black] (5,2.5) node {$\Psi$};
	 			
	 		\draw[-,thin, black] (9.5, 0)--(9.5,1.5);
	 		\draw[step=0.5cm, gray, very thin,dotted] (5.5, 0) grid (13.5,1.5);
	 		\draw[black] (5.5, 0)--(6, 0.5)--(6.5, 0)--(7, 0.5)--(7.5, 0);
	 		\draw[black](7.5, 0)--(8, 0.5);
	 		\draw[black](8, 0.5)--(8.5, 1)--(9, 1.5)--(9.5, 1)--(10, 0.5)--(10.5, 1)--(11, 0.5);
	 		\draw[black](11, 0.5)--(11.5, 0);
	 		\draw[black](11.5, 0)--(12, 0.5)--(12.5, 1)--(13, 0.5)--(13.5, 0);
	 			
	 		\draw[black] (2.7,-0.3) node {Tableaux de Young}; \draw[black] (9.5,-0.3) node {Chemin de Dyck}; 
	 		\draw[|->,dashed, black] (9.5, 2)--(9.5,3.5);
	 		\draw[black] (10,2.2) node {$\Phi^{-1}$};
	 		
	 		\draw[black] (9.5,4) node {$78534621\in S_8(132)$};	 		
	 			\draw[|->,thin, black] (4.5,4)--(6.5,4);
	 				\draw[black] (5.6,4.2) node {$\Theta$};
	 			
	 		\end{tikzpicture}
		\caption{Illustration graphique d'une correspondance par la bijection $\Theta$.}
	 		\label{fig13}
	 		
	 	\end{center}
	 \end{figure}
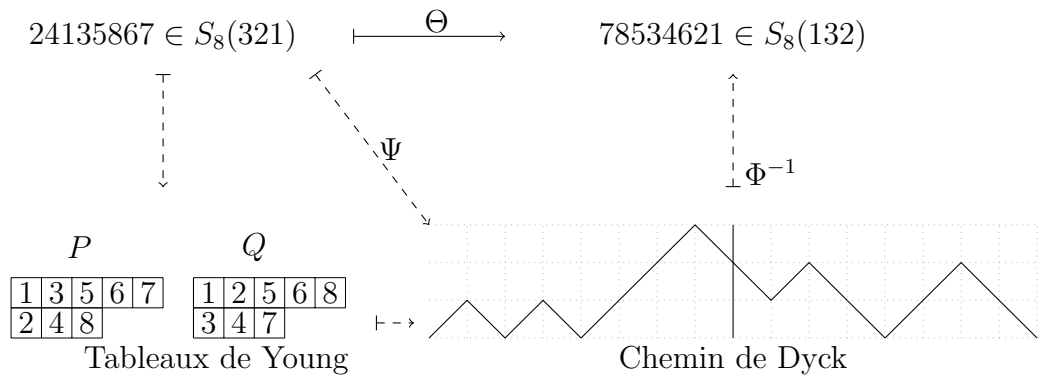

Nous présentons dans la Figure \ref{fig13} un exemple d'une correspondance par la bijection $\Theta$. Les détails sur les étapes de construction de $\Theta(\sigma)$ à partir d'une permutation $321$-interdite $\sigma=24135867$ ont été décrits dans les deux sous-sections précédentes.

	\chapter{Permutations évitant un motif de $S_3$}\label{chap2}

	\section{Introduction}\label{sec21}
	Plusieurs travaux récents étudient des statistiques sur les permutations interdisant un ou plusieurs motifs de $S_3$. 	Dans \cite{ARob}, Robertson et al.~ se sont concentrés sur les distributions du nombre de points fixes. Elizalde \cite{Eliz1} a généralisé le résultat de Robertson et al.~ et a étudié les distributions jointes du nombre de points fixes et du nombre d'excédances. Récemment, Dokos et al.~\cite{Dokos} ont étudié les distributions des statistiques $\inv$ et $\maj$. En nous inspirant des travaux de \cite{Dokos, ElizD, ElizP, Eliz1,Eliz2,ARob} et en particulier celui de \cite{ARandr}, notre présente contribution concerne la distribution du nombre de croisements sur les permutations interdisant un motif de $S_3$ et nous avons trouvé deux résultats. 

Le premier résultat concerne l'équidistribution du nombre de croisements sur les permutations 321, 132 et 213-interdites, c'est-à-dire, pour tout $n\geq 1$
	\begin{align}\label{res1}
	\sum_{\sigma\in S_n(321)} q^{\crs(\sigma)}=\sum_{\sigma\in S_n(132)} q^{\crs(\sigma)}=\sum_{\sigma\in S_n(213)} q^{\crs(\sigma)}.
	\end{align}
Pour prouver ces identités, nous allons montrer que la bijection $\Theta$ de Elizalde et Pak ainsi que la composée miroir-complément-inverse $\rci$ des involutions triviales préservent le nombre de croisements (la définition de l'involution $\rci$ sur l'ensemble des permutations sera donnée dans la Section \ref{sec23}). Nous ferons appel au raffinement d'un résultat de Randrianarivony prouvé dans \cite{ARandr} pour obtenir le développement en fraction continue de la fonction génératrice de ces distributions.
 
Le deuxième résultat concerne une relation entre les fonctions génératrices des polynômes énumérateurs de $\crs$ sur les permutations 231 et 312-interdites. Plus précisément, nous avons
\begin{align}\label{res2}
F(312;q,z)=\frac{1}{1-zF(231;q,z)},
\end{align}
où $\displaystyle F(\tau;q,z):=\sum_{\sigma\in S(\tau)} q^{\crs(\sigma)}z^{|\sigma|}=\sum_{n\geq 0}\left( \sum_{\sigma\in S_n(\tau)} q^{\crs(\sigma)}\right) z^{n}$ pour tout  motif $\tau$.
Pour ce résultat, nous allons prouver que l'application  $f:S_{n-1} \longrightarrow S_{n}^n:=\{a\in S_n: a(n)=1\}$ associant $\sigma$ à $\sigma^{-(n ,1)}$ est une bijection qui conserve également le nombre de croisements. Puisque $f(S_{n-1}(231))=S_{n}^n(312)$, nous obtiendrons une relation de récurrence menant à notre résultat (cf Section \ref{sec23}).
	 
Le reste de ce chapitre sera organisé comme suit. D'abord, dans la Section \ref{sec22}, nous prouverons que la bijection $\Theta$ étudiée dans le chapitre précédent préserve le nombre de croisements. Ensuite, dans les Sections \ref{sec23} et \ref{sec24}, nous établirons les preuves des \eqref{res1} et \eqref{res2}.  Enfin, avant de conclure, nous discuterons les liaisons entre nos résultats et ceux de  \cite{Bloom2,MansSh,ARandr,ARob}.

\section{Retour sur la bijection $\Theta:S_n(321)\rightarrow S_n(132)$}\label{sec22}

\subsection{Nouvelle formulation de la bijection $\Theta$}
En exploitant la connaissance de la bijection $\Theta$ présentée dans la section précédente, nous avons trouvé une nouvelle formulation de $\Theta$ qui n'utilise plus d'autres objets intermédiaires (cf Théorème \ref{thm12}). Pour prouver notre résultat, nous avons besoin de la proposition fondamentale suivante.
\bpro \label{prop15}
Soit $\sigma \in S_n(321)$ tel que $k=\sigma(n)$. Supposons que $j=|\{(a,b)\in \mathcal{M}(\sigma):a\leq k\}|+1$. 
Nous avons les propriétés suivantes:
\begin{itemize}
	\item[{\rm(i)}] $\Theta(\sigma)=\Theta(\pi)^{(n-k+j,j)}$ où $\pi=\red(\sigma(1\cdots n-1))$,
	\item[{\rm(ii)}] $n-k+j$ est la plus petite non-excédance de $\Theta(\sigma)$.
\end{itemize} 
\epro
	\begin{proof}
	Soit $\sigma \in S_n(321)$, $k=\sigma(n)$ et $\pi=\red[\sigma(1\cdots n-1)]$. Puisque $\sigma=\pi^{(n,k)}$, le couple $(P,Q)=\RSK(\sigma)$ est obtenu à partir de $(P',Q')=\RSK(\sigma(1\cdots n-1))$ en y insérant $(n,k)$. En suivant la logique de l'insertion de $(n,k)$ dans $(P',Q')$ pour obtenir $(P,Q)$, on observe que $\Psi(\sigma)$ peut être également  obtenu à partir de $\Psi(\pi)$. Pour une meilleure compréhension, nous raisonnons graphiquement et nous distinguons trois cas selon les valeurs de $k$. Dans tous les cas, soit $i-1$ (resp. $j-1$) la position du plus grand nombre inférieur à $k$ dans la première ligne (resp. la deuxième ligne) de $P'$. En d'autres termes, le nombre à la  $i$-ème (resp $j$-ème) colonne de la première ligne (resp. deuxième ligne) de $P'$ est supérieur à $k$ s'il en existe.  De plus, $j-1=|\{(a,b)\in \mathcal{M}(\sigma):a\leq k\}|$ est le nombre de valeurs d'excédance de $\sigma$ inférieur à $k$ affecté par l'algorithme d'affectation. Comme la permutation $\sigma$ est bi-croissante (voir \cite{ARef}), nous avons $i=k-(j-1)=k+1-j$, la somme des nombre de  non-excédances et le nombre de valeurs d'excédance inférieur à $k$ non bousculé de $\RSK$.

Supposons d'abord que $k=n$. Pour obtenir $(P,Q)$, il suffit d'ajouter $n$ à la fin des premières lignes de $P'$ et $Q'$ (voir Tableau \ref{fig:proof1}). Par conséquent, lorsque nous traduisons $(P,Q)$ en chemin de Dyck, nous obtenons $\Psi(\sigma)=\Psi(\pi)^{(L)}. ud . \Psi(\pi)^{(R)}$. Ainsi, lorsque nous numérotons les pas de $\Psi(\sigma)$ pour obtenir $\Theta(\sigma)$ selon la procédure décrite dans la Section \ref{sec23}, ceux du pas ascendant $u$ et du pas descendant $d$ nouvellement ajoutés sont respectivement $n+1-i$ et $j$. 

\begin{table}[h]
	\begin{center}

		\begin{tabular}{|cc|c|}	
			\hline	
			\begin{tikzpicture}
			\draw[black] (-0.25,0.25) node[scale=0.6pt] {P'=};
			\draw[black] (1,1) node[scale=0.7pt] {Avant insertion de $(n,n)$};
			\draw (0, 0.5) -- (2,0.5);\draw (0, 0.25) -- (2,0.25); \draw (0, 0) -- (1.25,0); 
			\draw (0, 0) -- (0,0.5);   \draw (2,0.25) -- (2,0.5); \draw (1.25, 0) -- (1.25,0.25); 
			\draw[black] (1,0.35) node[scale=0.6pt] {$\ldots$};	\draw[black] (1.4,0.75) node[scale=0.6pt] {$i-1$-ème}; \draw[->,thin, black] (1.75,0.7)--(1.95,0.5);
			\draw[black] (0.4,0.1) node[scale=0.6pt] {$\ldots$};
			\draw[black] (0.8,-0.3) node[scale=0.6pt] {$(j-1)$-ème}; \draw[->,thin, black] (0.8,-0.25)--(1.2,0); 	
			\end{tikzpicture}
			&\hspace{1cm} &
			\begin{tikzpicture}			
			\draw[black] (-0.25,0.25) node[scale=0.6pt] {P=};
			\draw[black] (2,1) node[scale=0.7pt] {Après insertion de $(n,n)$};
			\draw (0, 0.5) -- (2,0.5);\draw (0, 0.25) -- (2,0.25); \draw (0, 0) -- (1.25,0); 
			\draw (0, 0) -- (0,0.5); \draw (1.25, 0) -- (1.25,0.25); \draw (1.65,0.25) -- (1.65,0.5); \draw (2,0.25) -- (2,0.5); 
			\draw[black] (1,0.35) node[scale=0.6pt] {$\ldots$};	\draw[black] (0.4,0.1) node[scale=0.6pt] {$\ldots$};
			\draw[black] (0.8,-0.3) node[scale=0.6pt] {$(j-1)$-ème}; \draw[->,thin, black] (0.8,-0.25)--(1.2,0);
			\draw[black] (1.7,0.75) node[scale=0.6pt] {$i$-ème}; \draw[->,thin, black] (1.75,0.7)--(1.85,0.5);	\draw[red] (1.85,0.35) node[scale=0.6pt] {n};

			\draw[black] (2.5,0.25) node[scale=0.6pt] {Q=};
			\draw (2.75, 0.5) -- (4.75,0.5);\draw (2.75, 0.25) -- (4.75,0.25); \draw (2.75, 0) -- (4,0);
			\draw (2.75, 0) -- (2.75,0.5);  \draw (4, 0) -- (4,0.25); \draw (4.4,0.25) -- (4.4,0.5); \draw (4.75,0.25) -- (4.75,0.5);
			\draw[black] (3.75,0.35) node[scale=0.6pt] {$\ldots$};	\draw[black] (3.15,0.1) node[scale=0.6pt] {$\ldots$};
			\draw[red] (4.6,0.35) node[scale=0.6pt] {n};
			\end{tikzpicture}
			\\ 
			\begin{tikzpicture}
			\draw[black] (0.25,0.75) node[scale=0.6pt] {$\Psi(\pi)$=};
			\draw[dotted,black] (1, 0.75)--(1.5, 0.75);
			\draw[black] (1.5, 0)--(1.5, 1.5);  
			\draw[dotted,black] (1.5, 0.75)--(2, 0.75);
			\end{tikzpicture}
			& &
			
			\begin{tikzpicture}
			\draw[black] (0.25,0.75) node[scale=0.6pt] {$\Psi(\sigma)$=};
			\draw[dotted,black] (1, 0.75)--(1.5, 0.75);
			\draw[black] (1.3, 0.95) node[scale=0.6pt] {n+1-i};
			\draw[red] (1.5, 0.75)--(1.75, 1);
			\draw[red] (1.75, 1)--(2, 0.75);  \draw[black] (2, 0.9) node[scale=0.6pt] {j}; 
			\draw[black] (1.75, 0)--(1.75, 1.5);  
			\draw[dotted,black] (2, 0.75)--(2.5, 0.75); 
			\end{tikzpicture}
			\\		\hline
		\end{tabular}
	\caption{Insertion de $(n,n)$.}\label{fig:proof1}
	\end{center}
\end{table}

Supposons maintenant que $k<n$ et $k$ supérieur à tous les éléments de la première ligne de $P'$. Dans ce cas,  l'entier $k$ ne  bouscule aucun élément de la première ligne de $P'$, mais il existe une partie  non vide $A$ dans la deuxième ligne de $P'$ telle que ses éléments sont tous supérieurs à $k$. Comme on peut le voir dans le tableau \ref{fig:proof2}, $\Psi(\pi)^{(L)}$ se termine par une séquence de pas descendants produite par $A$.  Le chemin de Dyck correspondant $\Psi(\sigma)$ à $(P,Q)$ peut être obtenu à partir de $\Psi(\pi)$. En effet, nous avons $\Psi(\sigma)^{(R)}=d.\Psi(\pi)^{(R)}$ et nous obtenons $\Psi(\sigma)^{(L)}$ à partir de $\Psi(\pi)^{(L)}$ en insérant un nouveau pas ascendant (produit par $k$) juste avant le $j$-ème pas descendant.

\begin{table}[h]
	\begin{center}

		\begin{tabular}{|cc|c|}	
			\hline	
			\begin{tikzpicture}
			\draw[black] (-0.25,0.25) node[scale=0.6pt] {P'=};
			\draw[black] (1,1) node[scale=0.7pt] {Avant insertion de $(n,k)$};
			\draw (0, 0.5) -- (2,0.5);\draw (0, 0.25) -- (2,0.25); \draw (0, 0) -- (1.25,0); 
			\draw (0, 0) -- (0,0.5); \draw (0.75, 0) -- (0.75,0.25); \draw (1.25, 0) -- (1.25,0.25); \draw (2,0.25) -- (2,0.5); 
			\draw[black] (1,0.35) node[scale=0.6pt] {$\ldots$};	\draw[black] (1.4,0.75) node[scale=0.6pt] {$i-1$-ème}; \draw[->,thin, black] (1.75,0.7)--(1.95,0.5);				
			\draw[black] (0.4,0.1) node[scale=0.6pt] {$\ldots$};
			\draw[black] (0.4,-0.3) node[scale=0.5pt] {$j-1$-ème}; \draw[->,thin, black] (0.5,-0.25)--(0.65,0); \draw[black] (1,0.1) node[scale=0.6pt] {$A$};	
			\end{tikzpicture}
			&\hspace{1cm} &
			\begin{tikzpicture}			
			\draw[black] (-0.25,0.25) node[scale=0.6pt] {P=};
			\draw[black] (1,1) node[scale=0.7pt] {Après insertion de $(n,k)$};
			\draw (0, 0.5) -- (2,0.5);\draw (0, 0.25) -- (2,0.25); \draw (0, 0) -- (1.25,0); 
			\draw (0, 0) -- (0,0.5); \draw (0.75, 0) -- (0.75,0.25); \draw (1.25, 0) -- (1.25,0.25); \draw (1.65,0.25) -- (1.65,0.5); \draw (2,0.25) -- (2,0.5);
			\draw[black] (1,0.35) node[scale=0.6pt] {$\ldots$};	\draw[black] (0.4,0.1) node[scale=0.6pt] {$\ldots$};
			\draw[black] (0.4,-0.3) node[scale=0.5pt] {$j-1$-ème}; \draw[->,thin, black] (0.5,-0.25)--(0.65,0); \draw[black] (1,0.1) node[scale=0.6pt] {$A$};
			\draw[black] (1.7,0.75) node[scale=0.6pt] {$i$-ème}; \draw[->,thin, black] (1.75,0.7)--(1.85,0.5);	\draw[red] (1.85,0.35) node[scale=0.6pt] {k};

			\draw[black] (2.5,0.25) node[scale=0.6pt] {Q=};
			\draw (2.75, 0.5) -- (4.75,0.5);\draw (2.75, 0.25) -- (4.75,0.25); \draw (2.75, 0) -- (4,0);
			\draw (2.75, 0) -- (2.75,0.5);  \draw (4, 0) -- (4,0.25); \draw (4.4,0.25) -- (4.4,0.5); \draw (4.75,0.25) -- (4.75,0.5);
			\draw[black] (3.75,0.35) node[scale=0.6pt] {$\ldots$};	\draw[black] (3.15,0.1) node[scale=0.6pt] {$\ldots$};
			\draw[red] (4.6,0.35) node[scale=0.6pt] {n};
			\end{tikzpicture}
			\\ 
			\begin{tikzpicture}
			\draw[black] (0.25,0.75) node[scale=0.6pt] {$\Psi(\pi)$=};
			\draw[dotted,black] (1, 1.5)--(1.5, 1.5);
			\draw[black] (1.82,1.4) node[scale=0.6pt] {j}; \draw[black] (1.5, 1.5)--(1.75, 1.25);
			\draw[dotted,black] (1.75, 1.25)--(2, 1); 
			\draw[black] (2, 1)--(2.25, 0.75);
			\draw[black] (2.25, 0)--(2.25, 2);  
			\draw[dotted,black] (2.25, 0.75)--(2.75, 0.75);
			\end{tikzpicture}
			& &
			
			\begin{tikzpicture}
			\draw[black] (0.25,0.75) node[scale=0.6pt] {$\Psi(\sigma)$=};
			\draw[dotted,black] (1, 1.5)--(1.5, 1.5);
			\draw[black] (1.25, 1.75) node[scale=0.6pt] {n+1-i};\draw[red] (1.5, 1.5)--(1.75, 1.75);
			\draw[black] (1.75, 1.75)--(2, 1.5);  \draw[black] (2, 1.75) node[scale=0.6pt] {j}; 
			\draw[dotted,black] (2, 1.5)--(2.25, 1.25);
			\draw[black] (2.25, 1.25)--(2.5, 1);
			\draw[black] (2.5, 0)--(2.5, 2);  
			\draw[red] (2.5,1)--(2.75, 0.75);
			\draw[dotted,black] (2.75, 0.75)--(3.25, 0.75); 
			\end{tikzpicture}
			\\	\hline
		\end{tabular}
			\caption{Insertion de $(n,k)$ avec $k<n$ et sans bousculement.}\label{fig:proof2}
	\end{center}
\end{table}

Le dernier cas que nous discutons ici est celui où l'entier $k$ bouscule un élément de la première ligne de $P'$ lors de l'insertion de $(n;k)$, c'est-à-dire ~$k<n$ et  il y a au moins un élément de la première ligne de $P'$ qui est supérieur à $k$. Désignons par $x$ l'extrémité gauche de cet élément. Il se trouve donc à la $i$-ième colonne de $P'$. Dans la deuxième ligne de $P'$, nous indiquons encore par $A$ la partie non vide des éléments supérieurs à $k$ qui commence à  la $j$-ème colonne. Notons que tout élément de $A$ est également inférieur à $x$, c'est-à-dire ~$A\subseteq\{k+1,k+2,\ldots, x-1\}$. C'est pourquoi, dans le Tableau \ref{fig:proof3}, le chemin $\Psi(\pi)^{(L)}$ se termine par une séquence de pas descendants produite par $A$ suivie d'une séquence de pas ascendants produite par $x$ et les nombres à sa droite.  Lorsqu'on insère $k$ dans $P'$, $k$ remplace $x$ (c'est-à-dire, $k$ bouscule $x$) et $x$ crée une nouvelle cellule à la fin de la deuxième ligne. Comme dans les cas précédents, la déduction de $\Psi(\sigma)$ à partir de $\Psi(\pi)$ suit simplement la logique de l'insertion de $(n,k)$ dans $(P',Q')$ pour obtenir $(P,Q)$. En fait, nous avons $\Psi(\sigma)^{(R)}=u.\Psi(\pi)^{(R)}$ et $\Psi(\sigma)^{(R)}$ est obtenu à partir de $\Psi(\pi)^{(R)}$ en remplaçant le $i$-ème pas ascendant (produit par $x$) par un pas descendant parce que $x$ devient un élément de la deuxième ligne, puis en insérant un nouveau pas ascendant produit par $k$ juste avant le $j$-ème pas descendant. 

		\begin{table}[h]
	\begin{center}

		\begin{tabular}{|cc|c|}	
			\hline	
			\begin{tikzpicture}
			\draw[black] (-0.25,0.25) node[scale=0.6pt] {P'=};
			\draw[black] (1,1) node[scale=0.7pt] {Avant insertion de $(n,k)$};
			\draw (0, 0.5) -- (2.5,0.5);\draw (0, 0.25) -- (2.5,0.25); \draw (0, 0) -- (1.25,0); 
			\draw (0, 0) -- (0,0.5); \draw (0.75, 0) -- (0.75,0.25); \draw (1.25, 0) -- (1.25,0.25);\draw (1.75, 0.25) -- (1.75,0.5); \draw (2, 0.25) -- (2,0.5); \draw (2.5,0.25) -- (2.5,0.5); 
			\draw[black] (1,0.35) node[scale=0.6pt] {$\ldots$};	\draw[black] (0.4,0.1) node[scale=0.6pt] {$\ldots$}; \draw[black] (2.25,0.35) node[scale=0.6pt] {$\ldots$};
			\draw[black] (0.4,-0.3) node[scale=0.5pt] {$j-1$-ème}; \draw[->,thin, black] (0.5,-0.25)--(0.65,0); \draw[black] (1,0.1) node[scale=0.6pt] {$A$}; 
			\draw[black] (1.7,0.75) node[scale=0.6pt] {$i$-ème}; \draw[->,thin, black] (1.75,0.7)--(1.85,0.5); \draw[black] (1.85,0.35) node[scale=0.6pt] {$x$};	
			\end{tikzpicture}
			&\hspace{1cm} &
			\begin{tikzpicture}			
			\draw[black] (-0.25,0.25) node[scale=0.6pt] {P=};
			\draw[black] (1.5,1) node[scale=0.7pt] {Après insertion de $(n,k)$};
			\draw (0, 0.5) -- (2.5,0.5);\draw (0, 0.25) -- (2.5,0.25); \draw (0, 0) -- (1.5,0); 
			\draw (0, 0) -- (0,0.5); \draw (0.75, 0) -- (0.75,0.25); \draw (1.25, 0) -- (1.25,0.25);\draw (1.75, 0.25) -- (1.75,0.5); \draw (2, 0.25) -- (2,0.5); \draw (2.5,0.25) -- (2.5,0.5);
			\draw (1.5,0) -- (1.5,0.25); 
			\draw[black] (1,0.35) node[scale=0.6pt] {$\ldots$};	\draw[black] (0.4,0.1) node[scale=0.6pt] {$\ldots$}; \draw[black] (2.25,0.35) node[scale=0.6pt] {$\ldots$};
			\draw[black] (0.4,-0.3) node[scale=0.5pt] {$j-1$-ème}; \draw[->,thin, black] (0.5,-0.25)--(0.65,0); \draw[black] (1,0.1) node[scale=0.6pt] {$A$};
			\draw[black] (1.7,0.75) node[scale=0.6pt] {$i$-ème}; \draw[->,thin, black] (1.75,0.7)--(1.85,0.5);	\draw[red] (1.85,0.35) node[scale=0.6pt] {k};
			\draw[gray] (1.35,0.1) node[scale=0.6pt] {$x$};
			
			\draw[black] (3,0.25) node[scale=0.6pt] {Q=};
			\draw (3.25, 0.5) -- (5.25,0.5);\draw (3.25, 0.25) -- (5.25,0.25); \draw (3.25, 0) -- (4.5,0);
			\draw (3.25, 0) -- (3.25,0.5);  \draw (4.5, 0) -- (4.5,0.25); \draw (5.25,0.25) -- (5.25,0.5); \draw (4.5,0.) -- (4.5,0.25); \draw (4.25,0.) -- (4.25,0.25);
			\draw[black] (4.25,0.35) node[scale=0.6pt] {$\ldots$};	\draw[black] (3.75,0.1) node[scale=0.6pt] {$\ldots$};
			\draw[red] (4.4,0.1) node[ scale=0.6pt] {$n$};
			\end{tikzpicture}
			\\
			\begin{tikzpicture}
			\draw[black] (0.25,0.75) node[scale=0.6pt] {$\Psi(\pi)$=};
			\draw[dotted,black] (1, 1.5)--(1.5, 1.5);
			\draw[black] (1.8,1.45) node[scale=0.6pt] {j}; \draw[black] (1.5, 1.5)--(1.75, 1.25);
			\draw[dotted,black] (1.75, 1.25)--(2, 1);
			\draw[black] (2, 1)--(2.25, 0.75);
			\draw[black] (2.25, 0.75)--(2.5, 1);
			\draw[dotted] (2.5, 1)--(2.75, 1.25); \draw[black] (2.5, 0.9) node[scale=0.5pt] {n+1-i}; 
			\draw[black] (2.75, 0)--(2.75, 2); 
			\draw[dotted,black] (2.75,1.25)--(3.25,1.25);
			\end{tikzpicture}
			& &
			\begin{tikzpicture}
			\draw[black] (0.25,0.75) node[scale=0.6pt] {$\Psi(\sigma)$=};
			\draw[dotted,black] (1, 1.5)--(1.5, 1.5);
			\draw[black] (1.25, 1.75) node[scale=0.6pt] {n+1-i};\draw[red] (1.5, 1.5)--(1.75, 1.75);
			\draw[black] (1.75, 1.75)--(2, 1.5);  \draw[black] (2, 1.75) node[scale=0.6pt] {j}; 
			\draw[dotted,black] (2, 1.5)--(2.25, 1.25);
			\draw[black] (2.25, 1.25)--(2.5, 1);
			\draw[gray] (2.5,1)--(2.75, 0.75);
			\draw[dotted] (2.75, 0.75)--(3, 1); 
			\draw[red] (3,1)--(3.25,1.25);
			\draw[black] (3, 0)--(3, 2);
			\draw[dotted,black] (3.25,1.25)--(3.75,1.25); 
			\end{tikzpicture}
			\\		\hline
		\end{tabular}
			\caption{Insertion de $(n,k)$ avec $k<n$ et avec bousculement.}\label{fig:proof3}	
			
	\end{center}
\end{table}

		\hspace{-0.65cm}Dans tous les cas évoqués ci-dessus, lorsque nous regardons le chemin $\Psi(\sigma)$, nous pouvons en déduire trois choses.
\begin{itemize}
	
	 \item[$\bullet$] Premièrement, nous avons $\Theta(\sigma)(n+1-i)=j$.
	\item[$\bullet$] Deuxièmement, lorsque nous supprimons le pas ascendant numéroté $n+1-i$ et le pas descendant numéroté $j$ de $\Psi(\sigma)$, le chemin Dyck obtenu est $\Psi(\pi)$. Cela implique également que, lorsque nous retirons $\Theta(\sigma)(n+1-i)$ (qui est égal à $j$) dans $\Theta(\sigma)$, la réduction de la permutation obtenue n'est autre que $\Theta(\pi)$. En d'autres termes, nous avons $\Theta(\sigma)=\Theta(\pi)^{(n+1-i,j)}$. De plus, puisque $i=k+1-j$, alors nous avons $n+1-i=n-k+j$ et nous obtenons finalement $\Theta(\sigma)=\Theta(\pi)^{(n-k+j,j)}$. D'où la propriété (i).
	\item[$\bullet$] Troisièmement, puisque $n-k+j\geq j$ alors $n-k+j$ est une non-excédance de $\Theta(\sigma)$. De plus, sachant que $n-k+j$ est le minimum des numéros assigné aux pas ascendants de $\Psi(\sigma)^{(L)}$ auquel le tunnel associé est gauche ou centré, c'est aussi la plus petite non-excédance de $\Theta(\sigma)$ (voir Remarque \ref{rem23}). D'où la propriété (ii).
\end{itemize}
Ceci termine ainsi la preuve de la Proposition \ref{prop15}.
\end{proof}

\bthm \label{thm12}
	Pour tout $\sigma=\pi^{(n,k)} \in S_n(321)$ (c'est-à-dire $\sigma(n)=k$), on a
	$$\Theta(\sigma) =\begin{cases}
	\sigma & \text{ si\ \  $|\sigma|=1$}\\
	\Theta(\pi)^{(n-k+j,j)} & \text{ si\ \ $|\sigma|>1$}
	\end{cases},$$ 
	où $j=|\{(a,b)\in \mathcal{M}(\sigma):a\leq k\}|+1$  et $\pi=red(\sigma(1\cdots n-1))$.

	Pour tout $\alpha\in S(132)$, la relation suivante permet de calculer $\Theta^{-1}(\alpha)$
	$$\Theta^{-1}(\alpha)  =\begin{cases}
	\alpha & \text{ si $|\alpha|=1$}\\
	\Theta^{-1}(\beta)^{(|\alpha|,|\alpha|+\alpha(l)-l)}& \text{ si $|\alpha|>1$}
	\end{cases},$$
	où $\displaystyle l=\min_i\{i:\alpha(i)\leq i\}$ est la plus petite non-excédance de $\alpha$ et $\beta=\red(\alpha(1\cdots l-1)\alpha(l+1\cdots |\alpha|))$.
\ethm

\bpr
La propriété (i) de la Proposition \ref{prop15} assure que pour toute permutation $\sigma \in S_n(321)$, nous avons
$$\Theta(\sigma) =\begin{cases}
\sigma & \text{ si\ \  $|\sigma|=1$};\\
\Theta(\pi)^{(n-\sigma(n)+j,j)} & \text{ si\ \ $|\sigma|>1$ et $\pi=red(\sigma(1\ldots n-1))$}.
\end{cases}$$ 
où   $j=|\{(a,b)\in \mathcal{M}(\sigma):a\leq \sigma(n)\}|+1$. 

Supposons que $\alpha \in S_n(132)$ et $\displaystyle l=\min_i\{i:\alpha(i)\leq i\}$. Nous avons $\alpha=\beta^{(l,\alpha(l))}$ avec $\beta=\red(\alpha(1\ldots l-1)\alpha(l+1\ldots |\alpha|))$.
Les deux propriétés de la Proposition \ref{prop15} assurent qu'il existe un entier $k$ vérifiant $l=n-k+\alpha(l)$ tel que $\Theta^{-1}(\alpha)=\Theta^{-1}(\beta)^{(n,k)}$. Nous avons donc $k=n-l+\alpha(l)$ et $\Theta^{-1}(\alpha)  =
	\Theta^{-1}(\beta)^{(n,n-l+\alpha(l))}$. Plus précisément, pour toute permutation $\alpha\in S_n(132)$, nous avons
$$\Theta^{-1}(\alpha)  =\begin{cases}
\alpha & \text{ si $|\alpha|=1$};\\
\Theta^{-1}(\beta)^{(n,n-l+\alpha(l))}& \text{ si $|\alpha|>1$ et $\beta=\red(\alpha(1\ldots l-1)\alpha(l+1\ldots |\alpha|))$}.
\end{cases}$$
où $\displaystyle l=\min_i\{i:\alpha(i)\leq i\}$. Ce qui achève la démonstration du Théorème \ref{thm12}.
\epr

Nous allons terminer cette section par un exemple d'utilisation de cette nouvelle formulation. Considérons la permutation $\sigma=24135867$ de $S_8(321)$ de la Figure \ref{fig13}. D'abord, nous avons $\mathcal{M}(\sigma)=\{(2,3),(4,4),(8,7)\}$. Puisque nous sommes dans le cas $n=8$ et $\sigma(n)=7$, alors nous avons $j=|\{(a,b)\in \mathcal{M}(\sigma):a\leq \sigma(n)\}|+1=3$ et $n-\sigma(n)-j=4$. D'après la formule du Théorème \ref{thm12}, nous avons 
$\Theta(\sigma)=\Theta(\pi)^{(4,3)}$, avec $\pi=6743521$. Pour plus de clarté, nous avons détaillé dans le Tableau \ref{tab21} le calcul de $\Theta(\sigma)$ en utilisant la récursion du Théorème \ref{thm12}.  

	\begin{table}[h]
		\begin{center}
			
			\begin{tabular}{|c|c|c|c|}
				\hline
				$l$ & $\sigma_l=red(\sigma(1)\ldots \sigma(l))$ & $(l-\sigma(l)+j,j)$ & $\Theta(\sigma_{l-1})^{(l-\sigma(l)+j,j)}$\\
				\hline
				1& \textcolor{gray}{1} & - &\textcolor{gray}{1} \\
				2& 1\textcolor{gray}{2} & $(1,1)$&$\textcolor{gray}{1}2$\\
				3& 23\textcolor{gray}{1} & $(3,1)$ &$23\textcolor{gray}{1}$\\
				4& 241\textcolor{gray}{3} & $(3,2)$ &$34\textcolor{gray}{2}1$\\
				5& 2413\textcolor{gray}{5} &$(3,3)$ &$45\textcolor{gray}{3}21$\\
				6& 24135\textcolor{gray}{6} & $(3,3)$ &$56\textcolor{gray}{3}421$\\
				7& 241357\textcolor{gray}{6} & $(4,3)$ &$674\textcolor{gray}{3}521$\\
				8& $\sigma$=2413586\textcolor{gray}{7} & $(4,3)$&$\Theta(\sigma)=785\textcolor{gray}{3}4621$\\
				\hline
			\end{tabular}
			\label{tab21}
		\caption{Calcul de  $\Theta(24135867)$  en utilisant la formule du Théorème \ref{thm12}.}	
		\end{center}		
	\end{table}

Le calcul de l'inverse est aussi simple. Si $\alpha=78534621$, nous avons  $\Theta^{-1}(\alpha)=\Theta^{-1}(6743521)^{(8,7)}$ car $|\alpha|=8$, $k=4$ et $\alpha(k)=3$.
			

\subsection{Autres propriétés de la bijection $\Theta$}
Rappelons que la bijection $\Theta$ préserve le nombre de point fixes et d'excédances, c'est-à-dire, $(\fp,\exc)\left( \Theta(\sigma)\right) =(\fp,\exc)(\sigma)$ \text{ pour tout $\sigma \in S(321)$}( voir \cite{ElizP}). Pour prouver cela, Elizalde et Pak ont établi que les deux bijections  $\Psi$ et $\Phi$ échangent les statistiques (\fp,\exc) sur $S_n(321)$ et $S_n(132)$  et (\ct,\rt) sur les Chemins de Dyck. Plus précisément, nous avons
$\Theta: (\fp,\exc) \stackrel{\Psi}{\longrightarrow} (\ct,\rt)  \stackrel{\Phi^{-1}}{\longrightarrow} (\fp,\exc)$.  En utilisant la nouvelle formulation décrite dans le chapitre précédent, nous allons montrer qu'elle préserve aussi le nombre de croisements.  Pour cela, nous allons définir une opération $\otimes$ sur  $S(132)$ qui sera en correspondance avec l'opération habituelle $\oplus$ à travers la bijection $\Theta$. 
\bdf
{\rm La \textit{somme directe} de deux permutations $\sigma_1$ et $\sigma_2$, notée $\sigma_1\oplus \sigma_2$, est une permutation définie par $\sigma_1\oplus \sigma_2 =\sigma_1\cdot\sigma_2^{+|\sigma_1|}$, où $\cdot$ désigne la concaténation. On dit qu'une permutation $\sigma$ est $\oplus$-\textit{décomposable} si $\sigma$ peut s'écrire comme une somme directe de deux permutations non vides. Sinon, la permutation $\sigma$ est $\oplus$-\textit{irréductible}.} 
\edf
Il est clair que l'ensemble des permutations 321-interdites  est  stable par l'opération $\oplus$. De plus, l'opération $\oplus$ est associative  et nous avons $\crs(\sigma_1\oplus \sigma_2)=\crs(\sigma_1)+\crs(\sigma_2)$ pour tout $\sigma_1$ et $\sigma_2 \in S(321)$.  Avant de définir l'opération $\otimes$ sur les permutations 132-interdites, nous allons prouver la proposition suivante.
	\bpro\label{propo}
	Soit $\sigma \in S_n(132)$. Notons $T(\sigma):=	\{i\in [n]|\ \sigma^{-1}(i)>i<\sigma(i)\}$. Nous avons les propriétés suivantes:
		\begin{itemize}
		\item[{\rm (a)}]  $T(\sigma)=\emptyset$ si et seulement si $\sigma(1)=1$.
		\item[{\rm (b)}] Si $j\in T(\sigma)$, alors $i \in T(\sigma)$ pour tout $i\leq j$.
		\item[{\rm (c)}] Nous avons $|T(\sigma)|\leq \frac{n}{2}$. 
		\item[{\rm (d)}] Si $D=\Phi(\sigma)$, alors $|T(\sigma)|=|D^{(R)}|_u$.
		\item[{\rm (e)}] Si $k=1+|T(\sigma)|$, alors nous avons $\sigma^{(k,k)}\in S_{n+1}(132)$.
	\end{itemize}	
\epro
\begin{proof}
Soit $\sigma \in S_n(132)$. La propriété (a) est évidente puisque nous avons   $\sigma(1)=1$ si et seulement si $\sigma=12\cdots n$.  En utilisant le fait que $\sigma$ est 132-interdite, on peut facilement prouver par l'absurde la propriété $(b)$. Supposons que $t=\max T(\sigma)$. Selon $(b)$, nous avons $t=|T(\sigma)|$. En d'autres termes, nous avons $\sigma^{-1}(1\cdots t)>t<\sigma(1\cdots t)$. Cela n'est possible que si nous avons $n-t\geq t$ (c'est-à-dire  $t\leq \frac{n}{2}$), d'où la propriété $(c)$. La propriété $(d)$ vient de la Proposition \ref{prop23}.

Supposons maintenant que $k=1+|T(\sigma)|$. Premièrement, nous avons par définition $\sigma^{(k,k)}=\sigma^{k\rtimes 1}(1\cdots k-1)\cdot k\cdot\sigma^{k\rtimes 1}(k\cdots n)$. Selon $(b)$, nous avons $\sigma^{k\rtimes 1}(1\cdots k-1)=(\sigma(1\cdots k-1))^{+1}$ puisque $\sigma(1\cdots k-1)\geq k$. Deuxièmement, si $\sigma^{(k,k)}\notin S_{n+1}(132)$,  alors la séquence $k\cdot\sigma^{k\rtimes 1}(k\cdots n)$ contient au moins une occurrence de 132 car  $\sigma^{k\rtimes 1}(1\cdots k-1) \cdot k$ est 132-interdite. Supposons que nous avons $k<i_1<i_2$ tel que  $\red[k\sigma^{k\rtimes 1}(i_1)\sigma^{k\rtimes 1}(i_2)]=132$. Sachant  que $t=k-1$ est maximum, nous devrons examiner deux cas.
\begin{itemize}
	\item[$\bullet$] Si $\sigma^{-1}(k)> k$, alors $\red[\sigma(k)\sigma(i_1-1)\sigma(i_2-1)]=132$  puisque $\sigma(k)<k$.
	\item[$\bullet$] Si $\sigma^{-1}(k)\leq k$, alors  $\red[k\sigma(i_1-1)\sigma(i_2-1)]=132$.
\end{itemize}
Ceci contredit le fait que $\sigma$ est 132-interdite. Finalement, nous avons $\sigma^{(k,k)}\in S_{n+1}(132)$.
\end{proof}
Soit $\sigma \in S_n$ et  $p\geq 1$ un entier. Soit également $1\leq a_i,b_i\leq n+i$ pour tout $i\in [p]$. Nous écrivons
\begin{equation*}
\sigma^{\{(a_1,b_1),\ldots,(a_p,b_p)\}}:=(\cdots(\sigma^{(a_1,b_1)})\cdots)^{(a_p,b_p)} \in S_{n+p}. 
\end{equation*}
Si $\pi$ est une permutation de longueur $p$, nous écrivons aussi 
\begin{equation*}
\sigma^{(a,\pi)}:=\sigma^{\{(a,\pi(1)),(a+1,\pi(2)),\ldots,(a+p-1,\pi(p))\}}.
\end{equation*}
Par exemple, nous avons  $3142^{(3,\textcolor{gray}{213})}=3142^{\{(3,\textcolor{gray}{2}), (4,\textcolor{gray}{1}), (5,\textcolor{gray}{3})\}}=41\textcolor{gray}{2}53^{\{ (4,\textcolor{gray}{1}), (5,\textcolor{gray}{3})\}}=52\textcolor{gray}{31}64^{ (5,\textcolor{gray}{3})}=62\textcolor{gray}{413}75$. 

\bd
{\rm On définit le \textit{produit direct} de deux permutations $\alpha$ et $\beta$ par
$\alpha\otimes\beta:=\beta^{(k,\alpha^{+(k-1)})}=\beta^{k \rtimes |\alpha|}(1\cdots k-1)\cdot\alpha^{+(k-1)}\cdot\beta^{k \rtimes |\alpha|}(k\cdots |\beta|)$, où $k=1+|T(\beta)|$.}
\ed
 Exemple: $\textcolor{gray}{312}\otimes 543612=543612^{(3, \textcolor{gray}{534})}=87\textcolor{gray}{534}6912$ car $k=3$,  $312^{+2}=534$ et  $543612^{3\rtimes 6}=876912$.
On définit également de la même manière que celles avec l'opérateur $\oplus$ les notions de  $\otimes$-\textit{décomposable} et $\otimes$-\textit{irréductible} mais sur $S(132)$. 

\bpro\label{propxxx}
	Pour toutes permutations $\sigma_1$ et $\sigma_2\in S(132)$, nous avons \begin{equation*}
	\crs(\sigma_1\otimes\sigma_2)=\crs(\sigma_1)+\crs(\sigma_2).
	\end{equation*}
\epro
\begin{proof}
	Supposons que $\sigma=\sigma_1\otimes\sigma_2$ avec $\sigma_1, \sigma_2\in S(132)$. Par définition, nous avons $\sigma(k\cdots k+|\sigma_1|-1)=\sigma_1^{+(k-1)}$ et $\sigma(1\cdots k-1)\cdot\sigma(k+|\sigma_1|\cdots |\sigma_1|+|\sigma_2|)=\sigma_2^{k\rtimes |\sigma_1|}$, où $k=1+|T(\sigma_2)|$. Puisque $\sigma_1^{+(k-1)}$ est une permutation de $\{k,\ldots,|\sigma_1|+k-1\}$, en se référant aux arcs diagrammes du produit $\sigma$, aucun arc dans $\sigma$ lie une entrée dans $\sigma_1^{+(k-1)}$ et une entrée dans $\sigma_2^{k\rtimes |\sigma_1|}$.  Par conséquent, nous avons $\crs(\sigma)=\crs(\sigma_1^{+(k-1)})+\crs(\sigma_2^{k\rtimes |\sigma_1|})$. Ainsi, le résultat souhaité découle du fait que $\crs(\pi^{+a})=\crs(\pi)$ et $\crs(\pi^{a\rtimes b})=\crs(\pi)$ pour toute permutation $\pi$ et  pour tous entiers $a$ et $b\geq 1$.
\end{proof}

\bpro\label{prop25}
Pour toutes permutations  $\sigma_1, \sigma_2\in S(321)$, nous avons 
\begin{equation*}
\Theta(\sigma_1 \oplus \sigma_2)=\Theta(\sigma_2)\otimes\Theta(\sigma_1).
\end{equation*}
\epro	
	\begin{proof}
	Soit $\sigma=\sigma_1\oplus\sigma_2$ avec $\sigma_1 \in S_n(321)$ et $\sigma_2\in S_m(321)$. Notons d'abord $(P_1,Q_1)=\RSK(\sigma_1)$. Pour obtenir $(P,Q)=\RSK(\sigma)$, il faut insérer $\sigma_2'=\sigma_2^{+n}$ dans $(P_1,Q_1)$. Comme toutes les lettres de $\sigma_2'$ sont supérieures à toutes celles de $(P_1,Q_1)$, on obtient $(P,Q)=(P_1\cdot P_2,Q_1\cdot Q_2)$, où $(P_2,Q_2)=RSK(\sigma_2')$ et les points indiquent des concaténations. Ainsi, le chemin de Dyck produit par $(P,Q)$ est $D_1^{(L)}(D_2^{(L)}D_2^{(R)})D_1^{(R)}=\Psi(\sigma)$, où 	
	$D_1^{(L)}$, $D_1^{(R)}$, $D_2^{(L)}$ et $D_2^{(R)}$ sont respectivement les sous-chemins produits par $P_1$, $Q_1$, $P_2$ et $Q_2$.  En outre, nous avons $D_1^{(L)}D_1^{(R)}=\Psi(\sigma_1)$ et $D_2^{(L)}D_2^{(R)}=\Psi(\sigma_2)$. 
	
	Soit $k=|D_1^{(R)}|_u+1$. Nous résumons dans le tableau suivant tous les numéros assignés aux pas du chemin $\Psi(\sigma)$ avant d'obtenir $\Theta(\sigma)$.
		
		\begin{table}[h]
			\begin{center}
			
				\begin{tabular}{|c|c|c|c|}
					\hline
					Sous-chemin	& $D_1^{(L)}$  &  $ D_2^{(L)}D_2^{(R)}$ & $D_1^{(R)}$\\
					\hline
					Pas ascendants & $n+m,\ldots,k+m$& $k+m-1,\ldots,k+1,k$& $k-1,\ldots, 2,1$\\
					\hline
					Pas descendants &$1,2,\ldots,k-1$ &$k,k+1,\ldots,k+m-1$ &$k+m,\ldots,n+m$\\
					\hline
				\end{tabular}
				\label{table:tabx1}
					\caption{Les numéros assignés aux pas du chemin $\Psi(\sigma)$.}
			\end{center}				
		\end{table}

		\hspace{-0,68cm}Lorsque nous calculons $\pi=\Phi^{-1}(D_1^{(L)}(D_2^{(L)}D_2^{(R)})D_1^{(R)})=\Theta(\sigma)$ selon la procédure décrite dans la Section \ref{sec113}, nous obtenons deux sous-séquences à  examiner.
\begin{itemize}
	
	\item[$\bullet$] La séquence $\pi(k\cdots k+m-1)$ est produite par $D_2^{(L)}D_2^{(R)}$ et est une permutation de $\{k,k+1,\ldots, k+m-1\}$. En d'autres termes, nous avons $\pi(k\cdots k+m-1)=\pi_2^{+(k-1)}$ où $\pi_2=red(\pi(k\cdots k+m-1))=\Phi^{-1}(D_2^{(L)}D_2^{(R)})=\Theta(\sigma_2)\in S_{m}(132)$.
	\item[$\bullet$] La séquence $\pi(1\cdots k-1)\cdot\pi(k+m\cdots m+n)$ est produite par $D_1^{(L)}D_1^{(R)}$ et est une permutation de $[n+m]-\{k,k+1,\ldots k+m-1\}$. En examinant les première et troisième colonnes du Tableau \ref{table:tabx1}, on obtient $\pi^{-1}(1\cdots k-1)\geq k+m \leq \pi(1\cdots k-1)$. De plus, nous avons $\pi(1\cdots k-1)\cdot\pi(k+m\cdots m+n)=\pi_1^{k\rtimes m}$ avec $\pi_1=red(\pi(1\cdots k-1)\cdot\pi(k+m\cdots m+n))=\Phi^{-1}(D_1^{(L)}D_1^{(R)})=\Theta(\sigma_1)\in S_n(132)$. 
\end{itemize}
À partir de ces deux points, on peut écrire $\pi=\pi_1^{k\rtimes m}(1\cdots k-1)\cdot\pi_2^{+(k-1)}\pi_1^{k\rtimes m}(k\cdots n)$. 
D'après la Proposition \ref{propo}, nous avons $k=|T(\pi_1)|+1$. Donc, nous obtenons $\pi=\pi_2\otimes \pi_1=\Theta(\sigma_2)\otimes\Theta(\sigma_1)$. Ceci complète la preuve de la proposition.
	\end{proof}
Comme conséquence directe  de la Proposition \ref{prop25}, l'ensemble $S(132)$ est stable par le produit direct $\otimes$. L'opérateur $\otimes$ est également associatif sur $S(132)$ puisque la somme directe $\oplus$ l'est sur $S(321)$. De plus, une permutation $\sigma \in S(321)$ est $\oplus$-irréductible si et seulement si $\Theta(\sigma)\in S(132)$ est $\otimes$-irréductible. Un autre intérêt majeur de la Proposition \ref{prop25}, c'est qu'elle permet de prouver par récurrence que la bijection $\Theta$ conserve non seulement le nombre de points fixes et le  nombre d'excédances mais aussi le nombre de croisements (cf Théorème \ref{thm31}).

Pour toute permutation $\pi$ et pour tous entiers $a$ et $b$ satisfaisant $b\leq a\leq |\pi|+1$, nous adoptons les notations suivantes:
	\begin{itemize}
		\item[$\bullet$]  $A_1(\pi, a,b):=|\{b\leq i<a/\pi(i)<b\}|$, 
		\item[$\bullet$] $A_2(\pi,a,b):=|\{b\leq i<a/a<\pi^{-1}(i)\}|$, 
		\item[$\bullet$] $A_3(\pi,a,b):=|\{b\leq i<a/\pi^{-1}(i)<i<\pi(i)\}|$,
		\item[$\bullet$] $ A_4(\pi,a,b) :=|\{b\leq i<a/\pi(i)<i<\pi^{-1}(i)\}|$.
	\end{itemize} Nous allons prouver le lemme suivant qui a un rôle important pour une démonstration ultérieure.

	\blem \label{lem26} Soit $\pi$ une permutation et $a, b$ deux entiers tels que $1\leq b\leq a\leq 1+|\pi|$. Nous avons $\crs(\pi^{(a,b)})=\crs(\pi)+ \crs(\pi,a,b)$, où $\crs(\pi,a,b)=A_1(\pi, a,b)+A_2(\pi, a,b)+A_3(\pi, a,b)-A_4(\pi, a,b).$
	\elem
\begin{proof}
Considérons une permutation $\pi\in S_{n-1}$ et deux entiers $a$ et $b$ satisfaisant $b\leq a\leq n$. Soit $\sigma=\pi^{(a,b)}$. Les deux sous-séquences $s_1=\sigma(1\cdots b-1)$ et $s_2=\pi(1\cdots b-1)$ sont en ordre isomorphique, c'est-à-dire, $s_1(i)>s_1(j)$ si et seulement si $s_2(i)>s_2(j)$ pour tous entiers $i$ et $j$. Il en va de même pour $\sigma(a+1\cdots n)$ et $\pi(a\cdots n-1)$. Soit $i\in \{b,\ldots, a-1\}$,
	\begin{itemize}		
		\item[(a)] il  est facile de vérifier que tous les croisements de $\pi$, excepté ceux de la forme $(i,\pi^{-1}(i))$ (c'est-à-dire, $\pi(i)< i<\pi^{-1}(i)$), restent des croisements de $\sigma$, 
		\item[(b)] le nouvel arc inférieur $(a,b)$ croise avec tout arc $(i,\pi(i))$ tel que $\pi(i)<b\leq i$ 
		et avec tout arc $(i,\pi^{-1}(i))$ tel que $i<a<\pi^{-1}(i)$. Donc $(i,a)$ ou $(a,\pi^{-1}(i))$ qui n'est pas un croisement de $\pi$ devient un croisement pour $\sigma$, 	
		\item[(c)] si $\pi^{-1}(i)< i<\pi(i)$, alors $(\pi^{-1}(i),i)$ devient un croisement de $\sigma$ puisque $\pi^{-1}(i)<i< \sigma(\pi^{-1}(i))=i+1<\sigma(i)=\pi(i)+1$.		
	\end{itemize}
	De (a) nous obtenons les $\crs(\pi)-A_4(\pi, a,b)$, de (b) les $A_1(\pi, a,b)+A_2(\pi, a,b)$ et de (c) les $A_3(\pi, a,b)$. Ensemble, cela donne la relation souhaitée du lemme.
\end{proof}
	\brem\label{rem:rem2}
	{\rm Si $b<x<a$, nous avons $\crs(\pi^{(x,x)},a,b)=\crs(\pi,a-1,b)$.	Plus généralement, si $b<x_1<x_2<\cdots<x_p<a$, nous avons $\crs(\pi^{\{(x_1,x_1),(x_2,x_2),\ldots, (x_p,x_p)\}},a,b)=\crs(\pi,a-p,b)$}.
	\erem
	
	\blem \label{lem:lem33}
	Soit $\sigma=\pi^{(n,k)}\in S_n(321)$. Pour tout entier $i$ tel que $k\leq i <n$, on a ($\sigma(i)<k$ et $\sigma^{-1}(i)<i$) ou bien $\sigma^{-1}(i)< i<\sigma(i)$.
	\elem
	\begin{proof}
	Soit $\sigma=\pi^{(n,k)}\in S_n(321)$. Soit  $i$ un entier tel que $k\leq i <n$. Il est facile de montrer par l'absurde que \begin{itemize}
		\item[-] si $i$ est une excédance de $\sigma$, alors $\sigma^{-1}(i)<i$,
		\item[-] si $i$ est une non-excédance de $\sigma$, alors $\sigma(i)<k$ et $\sigma^{-1}(i)<i$.
	\end{itemize}
	\end{proof}
	
	\blem \label{lem:lem34}
	Soit $\sigma=\pi^{(n,k)}\in S_n(321)$ une permutation $\oplus$-irréductible. Nous avons $\crs(\Theta(\pi),n-k+j,j)=\crs(\pi,n,k)$, où $j-1=|\{(a,b)\in \mathcal{M}(\sigma):a\leq k\}|$. En outre, nous avons les propriétés suivantes:	
	\begin{itemize}		
		\item[{\rm (i)}] Si $\pi$ est $\oplus$-irréductible, alors $\crs(\pi,n,k)=n-k$. 
		\item[{\rm (ii)}] Si $\pi$ n'est pas $\oplus$-irréductible, alors il existe $l>1$ tel que $\alpha=\pi(1\cdots n-l-1)$ est $\oplus$-irréductible et $\crs(\pi,n,k)=\crs(\alpha,n-l,k)=n-l-k$. 
	\end{itemize} 	
	\elem
	\begin{proof}
Soit $\sigma=\pi^{(n,k)}\in S_n(321)$ une permutation $\oplus$-irréductible. Soit également $j-1=|\{(a,b)\in \mathcal{M}(\sigma):a\leq k\}|$. Nous allons examiner deux cas.

(i) $\pi$ est $\oplus$-irréductible: En utilisant le Lemme \ref{lem:lem33}, nous obtenons $\crs(\pi, n,k)=A_1(\pi,n,k)+A_3(\pi,n,k)=n-k$ car $A_2(\pi,n,k)=A_4(\pi,n,k)=0$. Puisque $\Theta(\pi)$ est également $\otimes$-irréductible et que $n-k+j$ est le minimum des non-excédances de $\Theta(\sigma)=\Theta(\pi)^{(n-k+j,j)}$, nous devons avoir ($\Theta(\pi)^{-1}(i)\geq n-k+j$ et $\Theta(\pi)(i)>i$) ou bien $\Theta(\pi)^{-1}(i)<i< \Theta(\pi)(i)$ pour $j\leq i<n-k+j$. Cela implique également que $\crs(\Theta(\pi), n-k+j,j)=A_2(\Theta(\pi),n-k+j,j)+A_3(\Theta(\pi),n-k+j,j)=n-k+j-j=n-k$ parce que  $A_1(\Theta(\pi),n-k+j,j)=A_4(\Theta(\pi),n-k+j,j)=0$.
Nous obtenons donc $\crs(\Theta(\pi),n-k+j,j)=\crs(\pi,n,k)=n-k$. 

(ii) $\pi$ est $\oplus$-décomposable :  D'une part, il existe un entier $m>1$ tel que $\pi=\pi_1\oplus \pi_2\oplus\cdots\oplus\pi_{m}$. Il n'est pas difficile de montrer que $|\pi_i|=1$ pour tout $i\neq 1$ et puisque $\sigma$ est $\oplus$-irréductible, il faut que $|\pi_1|\geq k$. Par conséquent, si nous désignons par $\alpha=\pi_1$ et $n-1-l=|\alpha|$ la longueur de $\alpha$, nous avons $\pi=\alpha \oplus 12\cdots l=\alpha^{\{(n-l,n-l),\ldots,(n-1,n-1)\}}$. Donc, selon la Remarque \ref{rem:rem2}, on obtient $\crs(\pi, n,k)=\crs(\alpha, n-l,k)$. Puisque $\alpha$ est $\oplus$-irréductible, d'après la propriété (i), nous obtenons $\crs(\pi, n,k)=\crs(\alpha, n-l,k)=n-l-k$. 
En revanche, nous avons $\Theta(\pi)=12\cdots l\otimes\Theta(\alpha)=\Theta(\alpha)^{(i,12\cdots l^{+(i-1)})}$ où $i=1+|T(\Theta(\alpha))|$.  De plus, comme $n-k+j$ est le minimum de non-excédance de $\Theta(\sigma)=\Theta(\pi)^{(n-k+j,j)}$, nous devons avoir $j<i<i+l-1<n-k+j$. Sinon, $\Theta(\sigma)$ peut ne pas être $\otimes$-irréductible ou ne pas être 132-interdite. Par conséquent, en utilisant de nouveau la Remarque \ref{rem:rem2} et la propriété (i) du lemme, nous obtenons $\crs(\Theta(\pi), n-k+j,j)=\crs(\Theta(\alpha), n-l-k+j,j)=n-l-k$. 

Dans tous les cas, on a $\crs(\Theta(\pi),n-k+j,j)=\crs(\pi,n,k)$. D'où le Lemme \ref{lem:lem34}.
\end{proof}

\bthm \label{thm31}
		Pour toute permutation  $\sigma \in S_n(321)$, nous avons $\crs(\Theta(\sigma))=\crs(\sigma)$.
\ethm
	\begin{proof}
			En combinant les  deux Lemmes \ref{lem26} et \ref{lem:lem34}, on peut procéder par récurrence sur $n$. Dans ce cas, le théorème  est évident pour $n =1,2,3$. Supposons que le théorème est vrai pour $k<n$. Considérons une permutation $\sigma \in S_n(321)$. 
		
		Supposons d'abord que $\sigma$ est $\oplus$-décomposable. On peut décomposer $\sigma$ comme une somme directe $\sigma=\oplus_{i=1}^l\sigma_i$ de $l$ permutations $\oplus$-irréductibles.  D'après la Proposition \ref{prop25}, on obtient $\Theta(\sigma)=\otimes_{i=1}^l\Theta(\sigma_{l+1-i})$, où les $\Theta(\sigma_i)$ sont tous $\otimes$-irréductibles. En appliquant l'hypothèse de récurrence, nous obtenons 
		\begin{equation*}
		\crs(\Theta(\sigma))=\sum_{k=1}^{l}\crs(\Theta(\sigma_k))=\sum_{k=1}^{l}\crs(\sigma_k)=\crs(\sigma). \end{equation*}		
		Supposons maintenant que $\sigma$ est $\oplus$-irréductible.
		Soit $\pi=red[\sigma(1\cdots n-1)] \in S_{n-1}(321)$ et $\Theta(\sigma)=\Theta(\pi)^{(n-\sigma(n)+j,j)}$, où $j-1=|\{(e,a)\in \mathcal{M}(\sigma)|e\leq \sigma(n)\}|$. Lorsque nous appliquons l'hypothèse de récurrence avec le Lemma \ref{lem:lem34}, on obtient 
		\begin{equation*}
		\crs(\Theta(\sigma))=\crs(\Theta(\pi))+\crs(\Theta(\pi),n-\sigma(n)+j,j)=\crs(\pi)+\crs(\pi,n,\sigma(n))=\crs(\sigma).
		\end{equation*} 
		Ceci termine ainsi la preuve du théorème.	
	\end{proof}

\begin{figure}[h]
		\begin{center}
			\begin{tikzpicture}
			\draw[black] (0,1) node {$1\ 2\ 3\ 4\ 5\ 6\ 7$}; 
			\draw (-0.9,1.2) parabola[parabola height=0.2cm,red] (-0,1.2); \draw[->,black] (-0.08,1.25)--(-0,1.2);
			\draw (-0.4,1.2) parabola[parabola height=0.2cm,red] (0.6,1.2); \draw[->,black] (0.55,1.25)--(0.6,1.2);
			\draw (-0.4,0.8) parabola[parabola height=-0.2cm,red] (0.6,0.8); \draw[->,black] (-0.35,0.75)--(-0.4,0.8);
			\draw (0.3,1.2) parabola[parabola height=0.2cm,red] (0.9,1.2); \draw[->,black] (0.85,1.25)--(0.9,1.2);
			\draw (0.3,0.8) parabola[parabola height=-0.2cm,red] (0.9,0.8); \draw[->,black] (0.35,0.75)--(0.3,0.8);
			
			\draw (-0,0.8) parabola[parabola height=-0.2cm,red] (-0.7,0.8); \draw[->,black] (-0.65,0.75)--(-0.7,0.8);		
			\draw (-1.,0.8) parabola[parabola height=-0.1cm,red] (-0.7,0.8);  \draw[->,black] (-0.95,0.75)--(-1.,0.8);		
			\draw[|->,thin, black] (1.5, 1)--(2,1);		
			\draw[black] (3.5,1) node {$1\ 2\ 3\ 4\ 5\ 6\ 7$}; 
			\draw (2.6,1.2) parabola[parabola height=0.3cm,red] (4.4,1.2); \draw[->,black] (4.32,1.25)--(4.4,1.2);
			\draw (2.9,1.2) parabola[parabola height=0.2cm,red] (4.1,1.2); \draw[->,black] (4.02,1.25)--(4.1,1.2);
			\draw (3.2,1.2) parabola[parabola height=0.1cm,red] (3.8,1.2); \draw[->,black] (3.7,1.25)--(3.8,1.2);
			
			\draw (4.4,0.8) parabola[parabola height=-0.2cm,red] (3.5,0.8);	\draw[->,black] (3.55,0.75)--(3.5,0.8);	
			\draw (4.1,0.8) parabola[parabola height=-0.2cm,red] (3.2,0.8); \draw[->,black] (3.25,0.75)--(3.2,0.8);
			\draw (3.5,0.8) parabola[parabola height=-0.1cm,red] (2.9,0.8); \draw[->,black] (3,0.75)--(2.9,0.8);
			\draw (3.8,0.8) parabola[parabola height=-0.2cm,red] (2.6,0.8); \draw[->,black] (2.7,0.75)--(2.6,0.8);
			
			\end{tikzpicture}
		\caption{Diagrammes d'arcs de $\sigma=4162735$ et $7652134=\Theta(\sigma)$.}
			\label{fig:arcdiag2}
			
		\end{center}
	\end{figure}
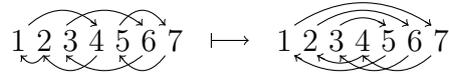

Nous  concluons cette section par un petit exemple d'illustration.  Dans la Figure \ref{fig:arcdiag2}, nous dessinons  les diagrammes d'arcs  pour $\sigma=4162735$ et $7652134=\Theta(\sigma)$. Observer que nous avons  $\crs(\Theta(\sigma))=\crs(\sigma)=5$.

\section{Equidistribution de $\crs$ sur $S_n(321)$, $S_n(132)$ et $S_n(213)$} \label{sec23}
Dans cette section, nous allons établir les preuves bijectives de \eqref{res1} à l'aide de la bijection $\Theta$ et quelques involutions triviales sur les permutations définies comme suit.
\bdf
{\rm Pour tout $\sigma \in S_n$,
\begin{itemize}
	\setlength\itemsep{-0.3em}
	\item[(a)] le \textit{miroir} $r$ de $\sigma$ est  $r(\sigma)=\sigma(n)\sigma(n-1)\cdots \sigma(1)$. En d'autres termes, nous avons $r(\sigma)(j)=\sigma(n+1-j)$ pour tout $j$. 
	\item[(b)] le \textit{complément} $c$ de $\sigma$ est $c(\sigma)=(n+1-\sigma(1))(n+1-\sigma(2))\ldots (n+1-\sigma(n))$. C'est-à-dire, $c(\sigma)(j)=n+1-\sigma(j)$ pour tout $j$.  
	\item[(c)] l'\textit{inverse} $i$ de $\sigma$ est $i(\sigma)$, une permutation telle que $i(\sigma)(j)=k$ si et seulement si  $\sigma(k)=j$. Nous dénotons habituellement $\sigma^{-1}=i(\sigma)$.
\end{itemize}
}
\edf
Pour simplifier l'écriture, nous écrivons $fg:=f\circ g$ pour toutes involutions $f$ et $g$ de $\{r, c, i\}$. Soit $\sigma\in S_n$.
D'après les définitions ci-dessus, le \textit{miroir-complément} de $\sigma$ est $rc(\sigma)$ tel que $rc(\sigma)(n+1-i)=n+1-\sigma(i)$ pour tout $i$, et le \textit{miroir-complément-inverse} de $\sigma$ est $rci(\sigma)$ tel que $rci(\sigma)(n+1-\sigma(i))=n+1-i$ pour tout $i$.
\bex 
{\rm
Si $\pi=41532$, nous avons $r(\pi)=23514$, $c(\pi)=25134$, $\pi^{-1}=25413$, $rc(\sigma)=43152$  et $rci(\pi)=35214$. }
\eex
Notons que pour tout $f\in \{r, c, i\}$, et pour tout ensemble de motifs $T$, nous avons 
\begin{equation*}
f(T) =\{\textit{f}(\tau): \tau \in T\}.
\end{equation*}
Nous avons également l'équivalence suivante qui est bien connue dans la littérature (voir \cite{SiS,West})
\begin{equation*}
\sigma \in S_n(T) \Longleftrightarrow f(\sigma) \in S_n(f(T)).
\end{equation*}
 Comme dans \cite{Eliz1} et \cite{Dokos}, nous avons souvent besoin de cette équivalence pour prouver nos résultats. 

	\bthm\label{thm23} Pour tout entier naturel $n\geq 1$, nous avons
\begin{equation*}
	\sum_{\sigma \in S_n(321)}q^{\crs(\sigma)}=\sum_{\sigma \in S_n(132)}q^{\crs(\sigma)}=\sum_{\sigma \in S_n(213)}q^{\crs(\sigma)}.
\end{equation*}
	\ethm
\begin{proof}
	D'après le  Théorème \ref{thm31}, puisque la bijection $\Theta$ préserve le nombre de croisements, nous avons la première égalité  
 \begin{equation}\label{id1}
 \sum_{\sigma \in S_n(321)}q^{\crs(\sigma)}=\sum_{\sigma \in S_n(132)}q^{\crs(\sigma)}.
 \end{equation} 
Pour la seconde identité, nous allons d'abord prouver que la bijection $rci$ préserve aussi le nombre de croisements. En effet, si $\pi=rci(\sigma)$ pour $\sigma \in S_n$, nous avons par définition $\pi(n+1-\sigma(i))=n+1-i$ pour tout $i \in [n]$.
 Nous avons par conséquent les équivalences suivantes:
 \begin{eqnarray*}
 	i<j<\sigma(i)<\sigma(j)\! \! & \! \! \Longleftrightarrow & \! \! n+1-\sigma(j)<n+1-\sigma(i)<n+1-j\! \!<\! \! n+1-i\\\
 	\! \! \! &\! \! \! \Longleftrightarrow & \! \! \! n+1-\sigma(j)\! \!<\! \! n+1-\sigma(i)\! \!<\! \! \pi(n+1-\sigma(j))\! \! <\! \! \pi(n+1-\sigma(i)).
 \end{eqnarray*}
De même, nous avons
 \begin{eqnarray*}
	\sigma(i)<\sigma(j)\leq i<j\! \! & \! \! \Longleftrightarrow & \! \! n+1-j\! \!<\! \! n+1-i\leq n+1-\sigma(j)<n+1-\sigma(i)\\
	\! \! \! &\! \! \! \Longleftrightarrow & \! \! \! \pi(n+1-\sigma(j))\! \!<\! \! \pi(n+1-\sigma(i))\! \!\leq\! \! n+1-\sigma(j)\! \! <\! \! n+1-\sigma(i).
\end{eqnarray*}
 Cela signifie que $(i,j)$ est un croisement de $\sigma$ si et seulement si $(n+1-\sigma(j),n+1-\sigma(i))$ est un croisement de $\pi$.  Par conséquent, puisque 213=$rci$(132), nous avons 
\begin{equation}\label{id2}
 \sum_{\sigma \in S_n(132)}q^{\crs(\sigma)}=\sum_{\sigma \in S_n(213)}q^{\crs(\sigma)}.
\end{equation}
Le Théorème \ref{thm23} s'ensuit en combinant les deux identités \eqref{id1} et \eqref{id2}. 
\end{proof}
Dans \cite{ARandr}, Randrianarivony a défini un $q,p-$analogue des nombres de Catalan $C_n(q,p)$ à travers la relation de récurrence
\begin{equation} \label{eq:cat}
C_n(q,p)=C_{n-1}(q,p)+q\sum_{k=0}^{n-2}p^kC_k(q,p)C_{n-1-k}(q,p)
\end{equation}
avec $C_0(q,p)=C_1(q,p)=1$ et il a proposé plusieurs interprétations combinatoires de  $C_n(q,p)$ en termes des statistiques sur les permutations sans croisement et les permutations sans imbrication. Le développement en fraction continue de la f.g.o. de $(C_n(q,p))_n$ est 
	\begin{equation}
\sum_{n\geq 0}C_n(q,p)z^n=\frac{1}{1-\displaystyle\frac{z}{				
		1-\displaystyle\frac{qz}{
			1-\displaystyle\frac{pz}{								
				1-\displaystyle\frac{qpz}{
					1-\displaystyle\frac{p^2z}{
						1-\displaystyle\frac{qp^2z}{				
							\ddots}
					}
			}}		
}}}. \label{eq:fc1}
\end{equation}

Puisque les permutations sans imbrication ne sont autres que les permutations 321-interdites (voir  \cite{Rakot}), un des résultats prouvé dans \cite{ARandr} peut être présenté comme suit.   
\bthm {\rm \cite{ARandr}}\label{thm:randr} Pour tout entier $n\geq 0$, nous avons  
\begin{align*}
\sum_{\sigma \in S_n(321)}q^{\exc(\sigma)}p^{\crs(\sigma)}=C_n(q,p).
\end{align*}
\ethm

Comme conséquence des Théorèmes \ref{thm23} et \ref{thm:randr}, nous avons le corollaire suivant.
	\bco\label{cor21} Soit $\tau \in \{321,132,213\}$. Pour tout entier  $n\geq 0$, nous avons
\begin{equation*}
\sum_{\sigma \in S_n(\tau)}q^{\crs(\sigma)}=C_n(1,q).
\end{equation*}
De plus, nous avons
\begin{equation}\label{fcRandr}
\sum_{\sigma\in S(\tau)} q^{\crs(\sigma)}z^{|\sigma|}=\frac{1}{1-\displaystyle\frac{z}{				
		1-\displaystyle\frac{z}{
			1-\displaystyle\frac{qz}{								
				1-\displaystyle\frac{qz}{
					1-\displaystyle\frac{q^2z}{
						1-\displaystyle\frac{q^2z}{				
							\ddots}
					}
			}}		
}}}.
\end{equation}
\eco

Nous avons observé que la fraction continue  \eqref{eqmain2} est apparue dans \cite{BEliz} comme étant la distribution du nombre d'occurrences d'un motif sur les permutations 231-interdites. La recherche des éventuelles correspondances entre ces résultats s'avère être intéressante, vu que  Corteel \cite{Cort} a établi la liaison entre les occurrences de motif, les croisements et les imbrications sur les permutations.

\section{Relation entre les distributions de $\crs$ sur $S_n(312)$ et $S_n(231)$}\label{sec24}

Nous allons établir la preuve de \eqref{res2} concernant une relation entre les distributions du nombre de croisements sur les permutations interdisant les motifs $231$ et $312$. Pour cela, nous adoptons les notations suivantes:
 \begin{itemize}
	\item[-] $Ut(\sigma):=\{i:\sigma^{-1}(i)<i<\sigma(i)\}$ et $\ut(\sigma)=|Ut(\sigma)|$;
	\item[-] $Lt(\sigma):=\{i/\sigma(i)<i<\sigma^{-1}(i)\}$ et $\lt(\sigma)=|Lt(\sigma)|$;
	\item[-] $Cr^*(\sigma)\!:=\!\{(i,j)/i<j<\sigma(i)<\sigma(j) \text{ ou } \sigma(j)<\sigma(i)\!<\! i<\!j \}$ et $\crs^*(\sigma)=|Cr^*(\sigma)|$.
\end{itemize}
En nous servant de ces notations, nous pouvons prouver le lemme suivant.
\blem \label{lem214}Pour toute permutation $\sigma$, nous avons $\crs(\sigma^{-1})=\crs(\sigma)+\ut(\sigma)-\lt(\sigma)$.
\elem
\begin{proof}
On peut d'abord remarquer que pour toute permutation $\sigma$, on a les propriétés suivantes
\begin{itemize}
	\item[(a)] $i\in Lt(\sigma)$ si et seulement si  $(i,\sigma^{-1}(i))$ est un croisement de $\sigma$;
	\item[(b)] $(i,j)\in Cr^*(\sigma)$ si est seulement si $(\sigma(i),\sigma(j))\in Cr^*(\sigma^{-1})$;
	\item[(c)] $i\in Lt(\sigma)$ si et seulement si $i \in Ut(\sigma^{-1})$.
\end{itemize}
Alors, nous obtenons respectivement  $\crs(\sigma)=\crs^*(\sigma)+\lt(\sigma)$, $\crs^*(\sigma^{-1})=\crs^*(\sigma)$ et $\lt(\sigma^{-1})=\ut(\sigma)$  de (a), (b) et (c). D'où,
\begin{eqnarray*}
\crs(\sigma^{-1})&=&\crs^*(\sigma^{-1})+\lt(\sigma^{-1})\\
&=&\crs^*(\sigma)+\ut(\sigma)\\
&=&\crs(\sigma)-\lt(\sigma)+\ut(\sigma).
\end{eqnarray*}
\end{proof}
\blem \label{lem215} Pour toute permutation $\sigma\in S_{n-1}$, nous avons $\crs(\sigma^{(n,1)})=\crs(\sigma)+\ut(\sigma)-\lt(\sigma)$.
\elem
\begin{proof}
C'est juste un cas particulier du Lemme \ref{lem26} avec $A_1(\sigma,n,1)=A_2(\sigma,n,1)=0$, $A_3(\sigma,n,1)=\ut(\sigma)$ et $A_4(\sigma,n,1)=\lt(\sigma)$.
\end{proof}
Dans toute la suite, nous notons $S_n^k:=\{\sigma \in S_n: \sigma(k)=1\}$.	Pour tout $\sigma\in S_{n-1}$, nous écrivons $\sigma^{-(a ,b)}:=(\sigma^{-1})^{(a ,b)}$ pour tous  $a, b \in [n]$. 
\bd
{\rm Pour tout entier $n\geq 1$, on définit une application $f$  par
	\begin{eqnarray*}
		f:S_{n-1} &\longrightarrow& S_{n}^n\\
		\sigma\hspace{ 0.2cm }&\longmapsto  & \hspace{ 0.1cm }\sigma^{-(n ,1)}. 
	\end{eqnarray*}}
\ed
Notons que l'application $f$ est bien définie et est bijective.  Elle vérifie également la propriété suivante.
\bthm \label{thm231}
La bijection $f$ conserve le nombre de croisements.
\ethm
\begin{proof}
Il suffit de combiner les deux lemmes précédents pour obtenir $\crs(\sigma^{-(n,1)})=\crs(\sigma)$ pour tout $\sigma\in S_{n-1}$. Ainsi, nous obtenons $\crs(f(\sigma))=\crs(\sigma)$.
\end{proof}
Comme conséquence directe du théorème précédent, nous avons 
\begin{equation*}
\sum_{\sigma \in S_n^n}q^{\crs(\sigma)}=\sum_{\sigma \in S_{n-1}}q^{\crs(\sigma)}.
\end{equation*}
De plus, nous avons le corollaire suivant.
\bco Pour tout $n\geq 1$, nous avons
\begin{equation}\label{rel111}
\sum_{\sigma \in S_n^n(231)}q^{\crs(\sigma)}=\sum_{\sigma \in S_{n-1}(312)}q^{\crs(\sigma)}.
\end{equation}
\eco
\bpr
Ce corollaire résulte du Théorème \ref{thm231} et l'identité $f(S_{n-1}(231))=S_{n}^n(312)$.
\epr

\bpro\label{prop33}
	Pour tout entier $n\geq 1$, nous avons
	\begin{equation}\label{rec1}
	F_n(312;q)=\sum_{j=0}^{n-1}F_{j}(231;q)F_{n-1-j}(312;q).
	\end{equation} 
\epro
\begin{proof}
	Sachons d'abord que $S_n^j(312)=\{\sigma_1\oplus \sigma_2|\sigma_1\in S_j^j(312),\sigma_2\in S_{n-j}(312)\}$ pour tout $j\geq 1$. Ainsi, l'application suivante est bien définie et est bijective
		\begin{eqnarray*}
		S_{j-1}(231)\times S_{n-j}(312) &\longrightarrow& \hspace{ 0.5cm }\ S_{n}^j(312)\\
		(\alpha,\beta)\hspace{ 1cm }&\longmapsto  & \hspace{ 0.4cm }\alpha^{-(j ,1)}\oplus\beta. 
	\end{eqnarray*}
En termes de distribution,   nous obtenons
	\begin{eqnarray*}
		\sum_{\sigma \in S_n^j(312)}q^{\crs(\sigma)}&=&  \sum_{\alpha \in S_{j-1}(231)}q^{\crs(\alpha)}\times \sum_{\beta \in S_{n-j}(312)}q^{\crs(\beta)} \\
		&=&  F_{j-1}(231;q)\times F_{n-j}(312;q). 
	\end{eqnarray*}
Par conséquent,  nous obtenons
 	\begin{eqnarray*}
 	F_n(312;q)&=&\sum_{j=1}^{n}\sum_{\sigma \in S_n^j(312)}q^{\crs(\sigma)}=\sum_{j=1}^{n}F_{j-1}(231;q)\times F_{n-j}(312;q). 
 \end{eqnarray*}
D'où la proposition. 
\end{proof}

\bthm \label{thm31x}
Nous avons la relation $\displaystyle F(312;q,z)=\frac{1}{1-zF(231;q,z)}$.
\ethm
\bpr
En utilisant la récurrence  (\ref{rec1}), nous obtenons
\begin{align*}
F(312;q,z)&=1+z\sum_{n\geq 0}\left( \sum_{j=0}^{n}F_{j}(231;q)F_{n-j}(312;q)\right) z^n\\
&=1+zF(231;q,z)F(312;q,z) ~~~~\text{(d'après Prop. \ref{profgo})}.
\end{align*}
En déduisant $F(312;q,z)$ de cette équation fonctionnelle, nous obtenons l'identité voulue du Théorème \ref{thm31x}.
\epr	
	\section{Liaisons avec d'autres résultats et discussions}\label{sec25}
Nos résultats sont évidemment liés avec d'autres résultats connus. Robertson \cite{ARob, ARob2}, par exemple,   a aussi construit une bijection $\Gamma: S_n(321)\rightarrow S_n(132)$ qui conserve également le nombre de points fixes. Bloom et Saracino \cite{Bloom, Bloom2} ont prouvé plus tard que la bijection $\Gamma$ conserve aussi le nombre d'excédances. Puis, dans \cite{Sarino}, Saracino a récemment montré que les deux bijections $\Theta$ et $\Gamma$ sont liées par une relation  simple décrite dans le théorème suivant. 
\bthm {\rm \cite{Sarino}} On a $\Gamma(\sigma)=\Theta\circ rci(\sigma)$ pour tout $\sigma \in S(321)$.
\ethm
Puisque les bijections $\Theta$ et $rci$ conservent le nombre de croisements, alors la bijection $\Gamma$ l'est également. Ainsi, nous avons le théorème suivant.
\bthm  Les bijections $\Gamma$ et $\Theta$ conservent le triplet de statistiques $(\fp, \exc, \crs)$.
\ethm
 De plus, comme la bijection $rci$ préserve aussi le couple (\fp, \exc, \crs), nous obtenons par conséquent les identités suivantes
\begin{equation*}
\sum_{\sigma \in S_n(213)}x^{\fp(\sigma)}q^{\exc(\sigma)}p^{\crs(\sigma)}=\sum_{\sigma \in S_n(132)}x^{\fp(\sigma)}q^{\exc(\sigma)}p^{\crs(\sigma)}=\sum_{\sigma \in S_n(321)}x^{\fp(\sigma)}q^{\exc(\sigma)}p^{\crs(\sigma)}.
\end{equation*}
Avec le résultat de Randrianarivony (cf Théorème \ref{thm:randr}), le cas $x=1$ nous mène à d'autres nouvelles interprétations du q,p-analogue des nombres de Catalan $C_n(q,p)$.   
\bthm
Pour tout $\tau \in \{213, 132, 321\}$, nous avons
	\begin{equation*}
	\sum_{\sigma \in S_n(\tau)}q^{\exc(\sigma)}p^{\crs(\sigma)}=C_n(q,p) \text{ pour tout entier $n\geq 1$}
	\end{equation*}
\ethm
Il est à rappeler que trouver l'expression de $F(\tau; q,z)$ reste encore ouvert pour tout motif $\tau \in \{123, 231, 312\}$.	La relation entre $F(312; q,z)$ et $F(231; q,z)$ que nous avons proposée ici (cf relation (\ref{res2}) ) réduit le reste des motifs à traiter. 

Nous concluons cette section par une remarque sur la récursion de la distribution du nombre d'inversions $\inv$ sur $S_n(321)$.
Si nous notons $I_n(q)=\sum_{\sigma \in S_n(321)}q^{\inv(\sigma)}$, nous avons la relation de récursion suivante qui a été conjecturée par Dokos et al~\cite{Dokos}
\begin{equation}\label{eq:inv}
I_n(q)=I_{n-1}(q) +\sum_{k=0}^{n-2}q^{k+1}I_k(q)I_{n-1-k}(q).
\end{equation}
En utilisant d'autres objets comme les chemins de Motzkin et les polyominos, cette relation à été démontrée par Cheng et al.~dans \cite{CEKS}.  Basés sur un algorithme particulier, Mansour et Shattuck \cite{MansSh} ont fourni une autre preuve. Nous observons qu'on peut facilement obtenir la relation \eqref{eq:inv} de la relation \eqref{eq:cat} en prenant  $p=q$, c'est-à-dire, $I_n(q)=C_n(q,q)$. En effet, il a été prouvé dans  \cite{MedVienot, ARandr} que $\inv(\sigma)=2\nes(\sigma)+\crs(\sigma)+\exc(\sigma)$ pour toute permutation $\sigma$. Ainsi, si $\sigma \in S_n(321)$,  alors nous avons $\inv(\sigma)=\crs(\sigma)+\exc(\sigma)$ et nous obtenons $I_n(q)=C_n(q,q)$. Le  développement en fraction continue de la f.g.o de $(I_n(q))$ présenté dans \cite[Thm. 1]{MansSh} est aussi obtenu de \eqref{eq:fc1} en prenant $p=q$.

\section{Conclusion}
Dans ce chapitre, nous avons étudié les distributions du nombre de croisements sur les permutations interdisant un motif de longueur 3 et nous avons justifié que les résultats obtenus sont liés avec d'autres résultats connus dans la littérature. Nous avons prouvé bijectivement en utilisant la bijection d'Elizalde et Pak que les distributions du nombre de croisements sur les ensembles $S_n(321)$, $S_n(132)$ et $S_n(213)$ sont égales. Nous avons également prouvé une relation entre les distributions du nombre de croisements sur les ensembles $S_n(312)$ et $S_n(231)$. Pour l'instant, nous n'avons aucune information sur la distribution du nombre de croisements sur l'ensemble $S_n(123)$.

	\chapter{Permutations évitant deux motifs de $S_3$}\label{chap3}

\section{Introduction}
Ce chapitre est une extension de notre deuxième article \cite{Rakot2}, un fruit d'une collaboration avec nos collègues Sandrataniaina et Randrianarivony. Nous allons présenter ici une étude complète sur l'énumération des permutations interdisant deux motifs de $S_3$ selon le nombre de croisements.  Notre technique est basée sur la manipulation des structures de $S_n(T)$,  pour tout $T \subset S_3$ tel que $|T|=2$.  

Grâce aux propriétés suivantes, certains cas sont triviaux et nous laissons la vérification au lecteur:
\begin{itemize}
	\item Si $\sigma \in S_n(312,321)$ ou $S_n(312,231)$, alors $\sigma$ est sans croisement.
	\item Si, $T=\{123,321\}$, nous avons 
	\begin{align*}
	F(T;q,z)=1+z+2z^2+(3+q)z^3+(1+2q+q^2)z^4. 
	\end{align*}
	Autrement dit, puisque le coefficient de $z^n$ dans $F(T;q,z)$ est $F_n(T; q)$ pour tout entier $n\geq 0$, nous avons
	\begin{align*}
	F_n(T; q)=
	\begin{cases}
	1 & \text{ si $n=0$ ou $1$};\\
	2 & \text{ si $n=2$};\\
	3+q& \text{ si $n=3$};\\
	1+2q+q^2 & \text{ si $n=4$};\\
	0 & \text{ si $n>4$}.
	\end{cases}.
	\end{align*}
\end{itemize}

Pour toute paire $T$ de la famille  $\mathcal{F}=\{\{321,231\}, \{321,132\}, \{321,213\}, \{123,132\}, \\ \{123,213\},\{312,123\}, \{231,123\}\}$, nous avons trouvé l'expression explicite de $F_n(T;q)$ ou $F(T;q,z)$. Pour les autres paires, nous restons sur une relation de récurrence pour $F_n(T;q)$.

Nous organisons le reste de ce chapitre en trois sections. Dans la Section \ref{sec41},  nous allons d'abord prouver une proposition fondamentale qui est nécessaire pour les preuves de nos résultats. Ensuite, dans la Section \ref{sec42}, en utilisant la proposition fondamentale prouvée dans la première section, nous établirons les calculs de $F_n(T;q)$ ou $F(T;q,z)$ pour toute paire de motifs $T$ de $S_3$. Enfin, nous conclurons le chapitre par un tableau de classement des motifs selon les distributions de $\crs$ sur $S_n(T)$ pour toute paire $T$ de permutations de $S_3$.

\section{Une proposition fondamentale}\label{sec41}
Soit $n$ un entier positif et $k\in [n]$. Nous notons $S_{n}^{k}:=\{\sigma \in S_n|\sigma(k)=1\}$ et $S_{n,k}:=\{\sigma \in S_n|\sigma(n)=k\}$. Pour tout ensemble de motifs $T$, nous notons également $F_n^k(T; q)$ et $F_{n,k}(T;q)$ les polynômes distributeurs du nombre de croisements sur les ensembles respectifs $S_n^k(\T)$ et $S_{n,k}(\T)$ et en particulier,
\begin{equation*}
F_n^k(q):=F_n^k(\emptyset; q) \text{ et } F_{n,k}(q):=F_{n,k}(\emptyset; q)
\end{equation*}
Pour tout entier $k$, nous adoptons aussi les notations suivantes: 
\begin{itemize}
	\item[-] $Ut_k(\sigma)\!:=\!\{i<k/\sigma^{-1}(i)\!<i\!<\!\sigma(i)\}$,  $\ut_k^-(\sigma)\!=\!|Ut_k(\sigma)|$ 
	et $\ut_k^+(\sigma)\!=\!\ut(\sigma)-\ut_k^-(\sigma)$;
	\item[-] $Lt_k(\sigma):=\{i<k/\sigma(i)<i<\sigma^{-1}(i)\}$, $\lt_k^-(\sigma)=|Lt_k(\sigma)|$ et $\lt_k^+(\sigma)=\lt(\sigma)-\lt_k^-(\sigma)$;
	\item[-] $\alpha_k(\sigma)=|\{i\geq k/\sigma(i)<k\}|$.
\end{itemize}

\blem\label{lem31x} Pour toute permutation $\sigma$, nous avons
\begin{equation*}
\crs(\sigma^{(k,1)})=\crs(\sigma)+\ut_k^-(\sigma)-\lt_k^-(\sigma)+\alpha_k(\sigma).
\end{equation*}
\elem
\begin{proof}
	C'est un cas particulier du Lemme \ref{lem26} avec $A_1(\sigma, 1,k)=0$, $A_2(\sigma, 1,k)=\alpha_k(\sigma)$, $A_3(\sigma, 1,k)=\ut_k(\sigma)$ et $A_4(\sigma, 1,k)=\lt_k(\sigma)$.
\end{proof}

\blem\label{lem32x}
Soit $\sigma$ une permutation et $\pi=rc(\sigma)$. Nous avons 
\begin{equation*}
\crs(\pi)=\crs(\sigma)+\ut(\sigma)-\lt(\sigma). 
\end{equation*}
\elem
\begin{proof}
	Soit $\sigma \in S_n$ et $\pi=rc(\sigma)$. Rappelons que par définition $rc(\sigma)(n+1-i)=n+1-\sigma(i)$. Ainsi, 
	\begin{equation*}
	i< \sigma(i)\Leftrightarrow {\rm rc}(\sigma)(n+1-i)<n+1-i.
	\end{equation*}
	Dans ce cas, on a aussi
	\begin{itemize}
		\item[-] $(i,j)\in Cr^*(\sigma)$ si et seulement si $(n+1-j,n+1-i)\in Cr^*(\sigma)$;
		\item[-] $i \in Lt(\sigma)$ si et seulement si $n+1-i \in Ut(\sigma)$.
	\end{itemize}
	Par conséquent,  $\crs^*(\pi)=\crs^*(\sigma)$, $\ut(\pi)=\lt(\sigma)$ et $\lt(\pi)=\ut(\sigma)$.\\ D'où $\crs(\pi)=\crs^*(\pi)+\lt(\pi)=\crs^*(\sigma)+\lt(\pi)=\crs(\sigma)-\lt(\sigma)+\ut(\sigma)$.
\end{proof}
\newpage 
\bd {\rm Soit $k \in [n]$. On définit les deux bijections suivantes				 		    
	\begin{eqnarray*}
		f_{k}:S_{n-1} &\longrightarrow& S_{n}^k\\
		\sigma\hspace{ 0.2cm }&\longmapsto  & \hspace{ 0.1cm }\sigma^{-(k,1)}. 
	\end{eqnarray*}	
	\begin{eqnarray*}
\text{ et }	g_{k}:\hspace{ 0.5cm } S_{n}^{k} \hspace{ 0.5cm } &\longrightarrow& \hspace{ 0.7cm }S_{n}^{n+1-k}\\
	\hspace{ 0.5cm }\sigma^{(k,1)}&\longmapsto  & \hspace{ 0.3cm } rc(\sigma)^{(n+1-k,1)}. 
	\end{eqnarray*}	}
\ed

\bthm\label{thm33x}
La bijection $g_{k}$ préserve le nombre de croisements pour tout $k\in [n]$.
\ethm
\begin{proof}
	Soit $\sigma^{(k,1)}\in S_n^{k}$ et $\pi^{(n+1-k,1)}=g_{k}(\sigma^{(k,1)})$ pour tout $\sigma \in S_{n-1}$. Il n'est pas difficile de prouver que
	\begin{equation} \label{eq26}
	\ut_{n+1-k}^-(\pi)=\lt_{k}^+(\sigma)  \text{ et }   \lt_{n+1-k}^-(\pi)=\ut_{k}^+(\sigma).
	\end{equation} 
	De plus, puisque  $|\{i<k/\sigma(i)\geq k\}|=|\{i\geq k/\sigma(i)<k\}|$, alors
	\begin{equation}\label{eq27}
	\alpha_{n+1-k}(\pi)=\alpha_{k}(\sigma).
	\end{equation}
	En effet, on a
	\begin{eqnarray*}
		\alpha_{n+1-k}(\pi)&=&|\{n-i\geq n+1-k/\pi(n-i)<n+1-k\}|,\\
		&=&|\{i\leq k-1/n-\sigma(i)<n+1-k\}|,\\
		&=&|\{i<k/\sigma(i)>k-1\}|,\\
		&=&|\{i<k/\sigma(i)\geq k\}|,\\
		&=&\alpha_{k}(\sigma).
	\end{eqnarray*}
	Par conséquent, en combinant les relations \eqref{eq26} et \eqref{eq27} avec les Lemmes \ref{lem31x} et \ref{lem32x}, on obtient
	\begin{eqnarray*}
		\crs(\pi^{(n+1-k,1)})&=&\crs(\pi)+\ut_{n+1-k}^-(\pi)-\lt_{n+1-k}^-(\pi)+\alpha_{n+1-k}(\pi),\\		
		&=&\crs(\sigma)+\ut(\sigma)-\lt(\sigma)+\lt_{k}^+(\sigma)-\ut_{k}^+(\sigma)+\alpha_{k}(\sigma), \\
		&=&\crs(\sigma)+\left( \ut(\sigma)-\ut_{k}^+(\sigma)\right) -\left( \lt(\sigma)-\lt_{k}^+(\sigma)\right) +\alpha_{k}(\sigma),\\
		&=&\crs(\sigma)+\ut_{k}^-(\sigma)-\lt_{k}^-(\sigma)+\alpha_{k}(\sigma),\\
		&=&\crs(\sigma^{(k,1)}).
	\end{eqnarray*}
\end{proof} 

\bthm\label{thm3x2}
La bijection $f_{n}$ préserve le nombre de croisements et la bijection $f_{n-1}$ satisfait
\begin{equation}\label{eqthm3x2}
\crs(f_{n-1}(\sigma))=	\crs(\sigma)+1 -\delta_{n-1,\sigma(n-1)},  
\end{equation} 
où $\delta$ désigne le symbole de Kronecker.	
\ethm

\begin{proof}
	Pour tout $\sigma \in S_{n-1}$, nous avons particulièrement 
	\begin{itemize}
		\item[-] $\ut_{n}^-(\sigma)=\ut(\sigma)$, $\lt_{n}^-(\sigma)=\lt(\sigma)$, $\alpha_{n}(\sigma)=0$; 
		\item[-] $\ut_{n-1}^-(\sigma)=\ut(\sigma)$, $\lt_{n-1}^-(\sigma)=\lt(\sigma)$ et $\alpha_{n-1}(\sigma)=1-\delta_{n-1,\sigma(n-1)}$, 
	\end{itemize}
	Ainsi, pour tout $\sigma \in S_{n-1}$, on combine les Lemmes \ref{lem214} et \ref{lem31x} pour obtenir 
	\begin{equation}\label{eq33x}
	\crs(\sigma^{-(n,1)})=\crs(\sigma) \text{ et  } \crs(\sigma^{-(n-1,1)})=\crs(\sigma)+1-\delta_{n-1,\sigma(n-1)}. 
	\end{equation}
	D'où le théorème.
\end{proof}

\bco \label{cor3x1} Pour tout $k \in [n]$, on a $F_n^{n+1-k}(q)=F_n^{k}(q)$.
En particulier, on a
\begin{eqnarray}
F_{n}^n(q)&=&F_{n}^1(q)=F_{n-1}(q) \text{ pour tout } n\geq 1 \label{eq3x1}\\
\text{ et } F_{n}^{n-1}(q)&=&F_{n}^2(q)=qF_{n-1}(q)+(1-q)F_{n-2}(q) \text{ pour tout } n\geq 2.\label{eq3x2}
\end{eqnarray}
\eco
\begin{proof}
Le Théorème \ref{thm33x} nous donne
	\begin{equation*}
	\sum_{\sigma\in  S_{n}^{n+1-k}} q^{\crs(\sigma)}=\sum_{\sigma\in  S_{n}^{k}} q^{\crs(\sigma)} \text{ pour tout $k\in [n]$}.
	\end{equation*}	
	C'est-à-dire, $F_n^{n+1-k}(q)=F_n^{k}(q)$ pour tout $k\in [n]$.
	En particulier, nous avons
	\begin{equation*} F_{n}^n(q)=F_n^1(q)=\sum_{\pi\in  S_{n}^{1}} q^{\crs(\pi)}=\sum_{1\oplus \sigma\in  S_{n}} q^{\crs(1\oplus\sigma)}=\sum_{\sigma\in  S_{n-1}} q^{\crs(\sigma)}=F_{n-1}(q) 
	\end{equation*}
	et, en utilisant l'identité \eqref{eqthm3x2} du Théorème \ref{thm3x2}, 
	\begin{eqnarray*}
		F_{n}^2(q)=F_{n}^{n-1}(q)&=&q\times \sum_{\sigma \in S_{n-1},\sigma(n-1)\neq n-1} q^{\crs(\sigma)}+\sum_{\sigma \in S_{n-1},\sigma(n-1)= n-1} q^{\crs(\sigma)}\\
		&=&q\left( \sum_{\pi\in  S_{n-1}} q^{\crs(\pi)}-\sum_{\sigma\oplus 1\in  S_{n-1}} q^{\crs(\sigma\oplus 1)}\right) +\sum_{\sigma\oplus 1\in  S_{n-1}} q^{\crs(\sigma\oplus 1)}\\
		&=&q\left( \sum_{\pi\in  S_{n-1}} q^{\crs(\pi)}-\sum_{\sigma\in  S_{n-2}} q^{\crs(\sigma)}\right) +\sum_{\sigma\in  S_{n-2}} q^{\crs(\sigma)}\\
		&=& q(F_{n-1}(q)-F_{n-2}(q))+F_{n-2}(q)\\
		&=&qF_{n-1}(q)+(1-q)F_{n-2}(q).
	\end{eqnarray*}	
\end{proof}
Soit $m$ et $n$ deux entiers naturels non nuls. Soit également $T\subset S_m$ et $k\in [n]$. Nous notons $T^{-1}=\{\tau^{-1}|\tau\in T \}$ et $T(i):=\{\tau(i)|\tau \in T\}$ pour tout $i\in [m]$. Dans la proposition suivante, nous verrons  à quoi ressemblent les versions restreintes des identités \eqref{eq3x1} et \eqref{eq3x2} du corollaire précédent.
\bpro\label{prop3x1}
Pour tout entier $n\geq 1$, nous avons les propriétés suivantes
\begin{itemize}
	\item[{\rm (i)}] Si $\min T^{-1}(1)>1$, alors   $F_{n}^{1}(T;q)=F_{n-1}(T;q)$;
	\item[{\rm (ii)}] Si   $\min T^{-1}(1)>2$, alors   $F_{n}^{2}(T;q)=qF_{n-1}(T;q)+(1-q)F_{n-2}(T;q)$; 
	\item[{\rm (iii)}]  Si  $\max T^{-1}(1)<m-1$,  alors 
	  	
	\hspace{3cm}$F_{n}^{n-1}(T;q)=qF_{n-1}(\T^{-1};q)+(1-q)F_{n-1,n-1}(T^{-1};q)$;  
	\item[{\rm (iv)}] Si $\max T^{-1}(1)<m$, alors    $F_{n}^{n}(T;q)=F_{n-1}(T^{-1};q)$. 
\end{itemize}
\epro
\begin{proof}
	Soit $m$ et $n$ deux entiers tels que $n\geq m>1$ et $T\subset S_m$. La preuve de la proposition est basée sur les deux faits évidents suivants:
	\begin{itemize}
		\item[(a)] Si $k<\min T^{-1}(1)$, alors $\sigma^{(k,1)} \in S_n^k(T) \Leftrightarrow\sigma \in S_{n-1}(T)$.
		\item[(b)] Si $n-m+\max T^{-1}(1)<k\leq n$, alors $\sigma^{-(k,1)} \in S_n^k(T) \Leftrightarrow\sigma \in S_{n-1}(T^{-1})$.
	\end{itemize}
	Les deux premières propriétés de la Proposition \ref{prop3x1} utilisent le fait (a). Si $\min \T^{-1}(1)\neq 1$, alors $1\oplus \sigma \in S_n^1(T)$ si et seulement si $\sigma \in S_{n-1}(T)$. Par suite, la propriété (i) se déduit des identités suivantes
	\begin{equation*}
	F_n^1(T;q)=\sum_{1\oplus \sigma \in S_{n}^1(T)} q^{\crs(1\oplus \sigma)}=\sum_{\sigma \in S_{n-1}(T)} q^{\crs(\sigma)}=F_{n-1}(T;q).
	\end{equation*}	
	De même, si $\min T^{-1}(1)>2$, alors $\sigma^{(2,1)} \in S_n^2(T) \Leftrightarrow\sigma \in S_{n-1}(T)$. Puisque $\crs(\sigma^{(2,1)})=\crs(\sigma)+1-\delta_{1,\sigma(1)}$ pour toute permutation $\sigma$, nous obtenons alors
	\begin{eqnarray*}
		F_{n}^{2}(T;q) &=& q\times \sum_{\sigma \in S_{n-1}(T),\sigma(1)\neq 1} q^{\crs(\sigma)}+\sum_{\sigma \in S_{n-1}(T),\sigma(1)= 1}q^{\crs(\sigma)}\\
		&=&q\left( F_{n-1}(T;q)-F_{n-1}^1(T;q)\right) +F_{n-1}^1(T;q)\\
		&=&qF_{n-1}(T;q)+(1-q)F_{n-2}(T;q) \text{ (car $F_{n-1}^1(T;q)=F_{n-2}(T;q)$)}.
	\end{eqnarray*}	
	Pour les deux propriétés restantes (ii) et (iv), en plus du fait (b), nous allons exploiter les bijections $f_{n}$ et $f_{n-1}$. Si $\max T^{-1}(1)<m-1$ (c'est-à-dire, $n-m+\max T^{-1}(1)<n-1$), alors $\sigma^{-(n-1,1)} \in S_n^{n-1}(T) \Leftrightarrow\sigma \in S_{n-1}(T^{-1})$.  Ceci implique  $f_{n-1}(S_{n-1}(T^{-1}))=S_n^{n-1}(\T)$. En utilisant le Théorème \ref{thm3x2}, nous obtenons
	\begin{eqnarray*}
		F_{n}^{n-1}(T;q) &=& q\times \sum_{\sigma \in S_{n-1}(\T^{-1}),\sigma(n-1)\neq n-1} q^{\crs(\sigma)}+\sum_{\sigma \in S_{n-1}(\T^{-1}),\sigma(n-1)= n-1}q^{\crs(\sigma)}\\
		&=&q\left( F_{n-1}(\T^{-1};q)-F_{n-1,n-1}(\T^{-1};q)\right) +F_{n-1,n-1}(\T^{-1};q)\\
		&=&qF_{n-1}(\T^{-1};q)+(1-q)F_{n-1,n-1}(\T^{-1};q).
	\end{eqnarray*}
De façon analogue, on prouve la  propriété (iv) en utilisant la bijection $f_n$. Ceci complète ainsi la preuve de la proposition.
\end{proof}
Nous avons mentionné au début de ce chapitre que la Proposition \ref{prop3x1} est fondamentale pour la suite car elle nous permettra de prouver facilement la majorité de nos résultats.

\section{Résultats d'énumérations}\label{sec42}

L'objectif est de trouver la forme explicite de  $F_n(T;q)$ ou $F(T;q,z)$, pour toute paire $T$ de $S_3$. Outre la manipulation des structures de nos objets combinatoires et la Proposition \ref{prop3x1}, un des outils que nous allons utiliser est la décomposition de la statistique $\crs$ en termes du nombre d'inversions, du nombre d'excédances et du nombre d'imbrications, une propriété prouvée dans \cite{MedVienot,ARandr} qui affirme que, pour toute permutation $\sigma$ 
\begin{equation}\label{eqarth}
\crs(\sigma)=\inv(\sigma)-\exc(\sigma)-2\nes(\sigma).
\end{equation}

\subsection{Permutations $(321,231)$-interdites}
On peut prouver par l'absurde le lemme suivant qui détermine la structure de $S_n(321,231)$.
\blem \label{lemx41}
Si $\sigma\in S_n(321,231)$, alors $\sigma^{-1}(1)\leq 2$.
\elem 
\bpro \label{prop4x2} Pour tout $n\geq 2$, nous avons la récursion suivante
\begin{eqnarray}
F_n(321,231;q)&=&(1+q)F_{n-1}(321,231;q)+(1-q)F_{n-2}(321,231;q). \label{eq31}
\end{eqnarray}
\epro
\begin{proof}
D'après le lemme précédent, $S_n(321,231)=S_n^1(321,231)\cup S_n^2(321,231)$ pour tout $n\geq 2$. Alors,  $F_n(321,231;q)=F_n^1(321,231;q)+F_n^2(321,231;q)$. Grâce aux propriétés (i) et (ii) de la Proposition \ref{prop3x1}, on obtient 
\begin{align*}
F_n^1(321,231;q)&=F_{n-1}(321,231;q);\\
F_n^2(321,231;q)&=qF_{n-1}(321,231;q)+(1-q)F_{n-2}(321,231;q).
\end{align*}
D'où, la récurrence (\ref{eq31}) .
\end{proof}
\bthm \label{thm413}
Nous avons
	\begin{equation*}
F(321,231; q,z)=\frac{1-qz}{1-(1+q)z-(1-q)z^2}.
\end{equation*}
\ethm
\begin{proof}
	Posons $T=\{321,231\}$. En effet, la f.g.o. associée à la récurrence \eqref{eq31} est 
	\begin{eqnarray*}
	F(T; q,z)&=& \sum_{n\geq 0}F_n(T;q)z^n		\\
				&=&1+z+\sum_{n\geq 2}\left( (1+q)F_{n-1}(T;q)+(1-q)F_{n-2}(T;q)\right) z^n\\
				&=&1+z+(1+q)z\left( \sum_{n\geq 0} F_{n}(T;q)z^n-1\right) +(1-q)z^2 \sum_{n\geq 0} F_{n}(T;q)z^n\\
				&=&1-qz+(1+q)zF(T; q,z) +(1-q)z^2F(T; q,z).
	\end{eqnarray*}
	Le théorème s'ensuit.
\end{proof}
Nous avons observé à travers une œuvre récente de  Bukata et al. \cite{BKLPRW} que le second membre de l'identité du Théorème \ref{thm413}
n'est autre que la fonction génératrice du triangle \href{https://oeis.org/A076791}{A076791} de OEIS \cite{OEIS} (voir \cite[Prop. 7]{BKLPRW}). Par conséquent, pour tous entiers $n$ et $k\geq 0$,
\begin{align*}
|\{\sigma \in S_n(321,231)|\crs(\sigma)=k\}|=A076791(n,k).
\end{align*}
\subsection{Permutations $(123,132)$ et $(123,213)$-interdites}
D'abord, on a $F_n(123,213;q)=F_n(123,132;q)$ pour tout $n\geq 0$ car $\{123,213\}=rci(\{123,132\})$. Il nous suffit alors d'étudier la structure de $S_n(123,132)$ et de calculer $F_n(123,132;q)$. De plus, sachant que $\{123,132\}=r(\{123,231\})$, alors chaque permutation (123,132)-interdite est un miroir d'une permutation (321,231)-interdite. Par conséquent, en se référant au Lemme \ref{lemx41}, nous avons le lemme suivant qui détermine la structure de $S_n(123,132)$.
\blem \label{lemx42}
Si une permutation $\sigma\in S_n(123,132)$, alors $\sigma^{-1}(1)\geq n-1$.
\elem

\bpro \label{thm421}Soit $\tau\in \{132,213\}$.  Pour tout $n\geq 2$, nous avons 
\begin{eqnarray}
F_{n}(123,\tau;q)&=&(1+q)F_{n-1}(123,\tau;q)+1-q. \label{eq32}
\end{eqnarray}
\epro
\begin{proof}
Le Lemme \ref{lemx42} entraine $S_n(123,132)=S_n^{n-1}(123,132)\cup S_n^{n}(123,132)$ pour tout  $n\geq 2$. D'où,
	$F_n(123,132;q)=F_n^{n-1}(123,132;q)+F_n^n(123,132;q)$. Selon les propriétés (i) et (iv) de la Proposition \ref{prop3x1}, puisque $F_{n,n}(123,132;q)=1$, on obtient
	\begin{eqnarray}
	F_n^n(123,132;q)&=&F_{n-1}(123,132;q), \label{eqp31}\\ 
	\text{ et }	F_n^{n-1}(123,132;q)&=&qF_{n-1}(123,132;q)+1-q. \label{eqp32}
	\end{eqnarray}
	En sommant \eqref{eqp31} et \eqref{eqp32}, nous obtenons
	\begin{equation*}
	F_{n}(123,132;q)=(1+q)F_{n-1}(123,132;q)+1-q.
	\end{equation*}
Puisque $F_n(123,213;q)=F_n(123,132;q)$, alors la proposition s'ensuit.
\end{proof} 
\bthm \label{thm423}
Pour tout $\tau \in \{132,213\}$, nous avons 
\begin{equation*}
F(123,\tau;q,z)=1+\frac{z(1-qz)}{(1-z)(1-(1+q)z)}.
\end{equation*}
\ethm
\begin{proof}
Pour $\tau \in \{132,213\}$, la f.g.o. associée à la récurrence \eqref{eq32} est  
	\begin{equation*}
	F(123,\tau;q,z)=1+z+(1+q)z(F(123,\tau;q)-1)+z\left( \frac{1}{1-z}-1-z\right).
	\end{equation*} 
Par conséquent,
	\begin{equation*}
	F(123,\tau;q,z)=1+\frac{z(1-qz)}{(1-z)(1-(1+q)z)}.
	\end{equation*}
\end{proof}
\brem
{\rm \label{remx} L'expression explicite de $F_n(123,\tau)$ suivante se déduit de la récurrence \eqref{eq32}: pour tout $\tau \in \{132,213\}$,
\begin{align}\label{eqok}
F_n(123,\tau)=\frac{(1+q)^{n-1}-1+q}{q}, \text{ pour tout } n\geq 1.
\end{align}}
\erem
Le résultat trouvé ici est une nouvelle interprétation du triangle \href{https://oeis.org/A299927}{A299927} de OEIS \cite{OEIS}, un triangle  que Bukata et al. ont récemment interprété en termes d'autres statistiques sur les permutations évitant une paire de motifs de longueur 3 \cite[Prop. 11]{BKLPRW}. Plus précisément, pour tout $\tau \in \{132,213\}$, pour tous entiers $n$ et $k\geq 0$,
\begin{align*}
|\{\sigma \in S_n(123,\tau)|\crs(\sigma)=k\}|=A299927(n,k).
\end{align*}
\bco \label{corok}
Nous avons $F_n^{n-1}(123,132;q)=F_{n,2}(123,213;q)=(1+q)^{n-2}$.
\eco
\bpr
Depuis la relation \eqref{eqp32}, on obtient en utilisant \eqref{eqok}
\begin{align*}
F_n^{n-1}(123,132;q)&=(1+q)^{n-2}.
\end{align*}
Puisque $S_{n,2}(123,213)=rci(S_n^{n-1}(123,132))$, alors 
\begin{align*}
F_{n,2}(123,213;q)=F_n^{n-1}(123,132;q)=(1+q)^{n-2}.
\end{align*}
\epr
On obtient également des nouvelles interprétations du triangle de Pascal \href{https://oeis.org/A007318}{A007318} car le Corollaire  \ref{corok} se traduit ainsi comme suit, pour tous entiers $n$ et $k\geq 0$,
\begin{align*}
|\{\sigma \in S_n^{n-1}(123,132)|\crs(\sigma)=k\}|=|\{\sigma \in S_{n,2}(123,213)|\crs(\sigma)=k\}|=\binom{n-2}{k}.
\end{align*}

 \subsection{Permutations $(321,132)$ et $(321,213)$-interdites}	
Comme  $\{321,132\}=rci(\{321,213\})$, alors $F_n(321,213;q)=F_n(321,132;q)$ pour tout $n\geq 0$.
\bthm \label{prop35} On note $[n]_q=1+q+\cdots+q^{n-1}$  pour tout entier $n\geq 1$. Pour tout $\tau \in\{132,213\}$, nous avons
\begin{equation*}
F_n(321,\tau;q)=1+\displaystyle\sum_{k=1}^{n-1}[n-k]_{q^k}.
\end{equation*} 
\ethm
\begin{proof} 
	Comme structure, nous avons $S_n(321,213)=S_n^1(321,213)\cup\{\alpha_2,\alpha_3,\ldots,\alpha_n\}$, où $\alpha_j=(n-j+2)\cdots (n-1)n12\cdots (n+1-j)$  pour tout $j\in [n]$.  Depuis cette structure,	
	\begin{eqnarray*}
		F_n(321,213;q) &=&F_{n}^1(321,213;q)+\sum_{j=2}^{n}q^{\crs(\alpha_j)}.
	\end{eqnarray*}
	Puisque $F_{n}^1(321,213;q)=F_{n-1}(321,213;q)$ et  $\crs(\alpha_j)=\inv(\alpha_j)-\exc(\alpha_j)=(j-1)(n-j)$ pour tout $j$, alors	
	\begin{eqnarray*}
		F_n(321,213;q)  &=& F_{n-1}(321,213;q)+\sum_{j=2}^{n}q^{(j-1)(n-j)}.
	\end{eqnarray*}
Quand on résout cette récurrence avec $F_1(321,213;q)=1$, nous obtenons $$F_n(321,213;q)=1+\sum_{k=1}^{n-1}\sum_{j=1}^{k}q^{j(k-j)}=1+\sum_{k=1}^{n-1}[n-k]_{q^k}.$$ 
Le théorème vient immédiatement de l'identité $F_n(321,213;q)=F_n(321,132;q)$.
\end{proof}

\subsection{Permutations $(123,312)$ et $(123,231)$-interdites}\label{subsec1}

Avant de calculer les expressions de  $F_n(123,312;q)$ et $F_{n}(123,231;q)$, nous allons prouver une identité qui les relie.
\bpro \label{pro441}
Pour tout $n\geq 2$, nous avons 
\begin{eqnarray}
F_n(123,312;q)&=&n-1+F_{n-1}(123,231;q); \label{eq33}
\end{eqnarray}
\epro
\begin{proof}
	Nous pouvons vérifier que $S_n(123,312)=\{\pi_1, \pi_2, \ldots, \pi_{n-1}\}\cup S_n^n(123,312)$, où $\pi_j=j\cdots 21n(n-1)\cdots(j+1)$ pour tout $j\in [n]$. Dans ce cas, nous avons
	\begin{equation*}
	F_n(123,312;q)=\sum_{j=1}^{n-1}q^{\crs(\pi_j)}+F_n^n(123,312;q).
	\end{equation*}
	Puisque $\{123,312\}^{-1}=\{123,231\}$, alors $F_n^n(123,312;q)=F_{n-1}(123,231;q)$ (d'après la propriété (iv) de la Proposition \ref{prop3x1}). De plus, comme $\crs(\pi_j)=0$ pour tout $j\in [n]$, alors la relation  \eqref{eq33} s'obtient immédiatement.
\end{proof} 
\bco \label{cor442}
On a la relation
\begin{equation*}
F(312,123; q,z) =1+\left(\frac{z}{1-z} \right)^2 +zF(231,123; q,z).
\end{equation*}
\eco
\begin{proof}
En utilisant la relation \eqref{eq33}, nous obtenons		
\begin{eqnarray*}
		F(123,312; q,z)&=& 1+z+ \sum_{n\geq 2}\left( n-1+F_{n-1}(123,231;q)\right) z^n		\\
		&=& 1+z+\left(\frac{z}{1-z}\right)^2 + z (F(123,231; q,z)-1).
\end{eqnarray*}
L'équation fonctionnelle  obtenue est équivalente à celle du corollaire.
\end{proof}

Pour tout entier  $k\in [n]$, notons d'abord
\begin{equation*}
 \sigma_{n,k,j}=\begin{cases}
(k+j)\cdots (k+2)(k+1)n(n-1)\cdots(k+j+1)k\cdots 21 & \text{ si } j<n-k;\\ 
n\cdots 2 1 & \text{ si } j=n-k.
 \end{cases}
\end{equation*}
Par convention, $\sigma_{n,0,j}=j\cdots 1n(n-1)\cdots(j+1)$ pour tout $j\in [n-1]$. 
\brem\label{rem51}
	{\rm Nous avons $\crs(\sigma_{n,k,j})=\crs(\sigma_{n,k,n-k-j})$ puisque $\sigma_{n,k,n-k-j}=rci(\sigma_{n,k,j})$. En particulier $\crs(\sigma_{n,0,j})=0$ pour $1\leq j<n$  et    $\crs(\sigma_{n,k,n-k})=0$  pour tout $k\in [n]$}. 
\erem
Le lemme suivant donne le nombre de croisements de  $\sigma_{n,k,j}$ pour tous entiers $k\in [n]$ et $j\in [n-k-1]$.
\blem \label{lem51}
	Nous avons les propriétés suivantes
	
	i) Si $k\geq \frac{n}{2}$, alors  $\crs(\sigma_{n,k,j})=j(n-k-j)$ pour tout $j\leq n-k$.
		
		ii) Si $\frac{n-1}{3}\leq k< \frac{n}{2}$, alors 
		\begin{equation*}
		\crs(\sigma_{n,k,j})=\begin{cases}
		\binom{j}{2}+j(k-j)	& \text{ si } \ \ \  j\leq \frac{n-1-k}{2};\\
		\binom{n-k-j}{2}+(n-k-j)(k-(n-k-j))	& \text{ si } \ \ \ \frac{n-1-k}{2}\leq j\leq n-k-1.\\
		\end{cases}
		\end{equation*}
		
		iii) Si $k< \frac{n-1}{3}$, alors 
		\begin{equation*}
		\crs(\sigma_{n,k,j})=\begin{cases}
		\binom{j}{2}+j(k-j)	& \text{ si } \ \ \ j\leq k;\\
		\binom{k}{2}	& \text{ si }  \ \ \ k< j\leq n-2k-1;\\
		\binom{n-k-j}{2}+(n-k-j)(k-(n-k-j))	& \text{ si }\ \ \  n-2k\leq j\leq n-k-1.\\
		\end{cases}
		\end{equation*}
\elem
\begin{proof} Selon les valeurs de $k$, nous utilisons simplement la structure de $\sigma_{n,k,j}$ pour connaître le nombre de croisements. 
	\begin{itemize}
		\item[i)] Le premier cas est clair parce que, si $k\geq \frac{n}{2}$, alors l'ensemble des croisements de $\sigma_{n,k,j}$ est $\{(x,y)| x\in [j] \text{ et }  y\in [j+1;n-k-j] \}$.
		\item[ii)] Assumons maintenant que $\frac{n-1}{3}\leq k< \frac{n}{2}$. Notons que $\frac{n-1}{3}\leq k\Leftrightarrow\frac{n-k-1}{2}\leq k$. Soit $j\leq n-k-1$.
		\begin{itemize}
			\item[-] Si $j\leq \frac{n-k-1}{2}$, l'ensemble des croisements de $\sigma_{n,k,j}$ est $\{(x,y)| x\in [j] \text{ et }  y\in [j+1;k] \}\cup\left(\cup_{i=1}^{j} \{(x,k+i)| x\in[i+1;j] \}\right)$. Ainsi, $\crs(\sigma_{n,k,j})=	j(k-j)+\binom{j}{2}$.
			\item[-] Si $j\geq \frac{n-k-1}{2}$, alors, d'après la Remarque \ref{rem51},  $\crs(\sigma_{n,k,j})=\crs(\sigma_{n,k,n-k-j})=\binom{n-k-j}{2}+(n-k-j)(k-(n-k-j))$.
		\end{itemize}
	\item[iii)] De façon similaire, sachant que $k<\frac{n-1}{3}\Leftrightarrow k<n-2k-1$, nous pouvons simplement distinguer les trois cas $j\leq k$, $k\leq j \leq n-2k-1$ et $n-2k-1\leq j <n-k-1$ pour prouver la dernière identité du lemme.		
	\end{itemize}
\end{proof}

\bthm \label{thm311} Pour tout entier $n\geq 2$, nous avons
	\begin{eqnarray*}
		F_n(123,312;q)&=&n+\sum_{k=1}^{\lfloor\frac{n-1}{3}\rfloor}\left(2\sum_{j=1}^{k}q^{\binom{j}{2}+j(k-j)}+(n-3k-1)q^{\binom{k}{2}}\right) \\
		&&+\sum_{k=\lfloor\frac{n-1}{3}\rfloor+1}^{\lfloor\frac{n-1}{2}\rfloor}\left( 2\sum_{j=1}^{\lceil\frac{n-k-1}{2}\rceil}q^{\binom{j}{2}+j(k-j)}-\gamma_{n,k}(q)\right) +\sum_{k=\lceil\frac{n}{2}\rceil}^{n-2}\sum_{j=1}^{n-k-1}q^{j(n-k-j)}.
	\end{eqnarray*}
	où $\gamma_{n,k}(q)=\begin{cases}
	0 & \text{ si } n-k-1 \text{ est pair};\\
	q^{\displaystyle \binom{\lceil\frac{n-k-1}{2}\rceil}{2}+\lceil\frac{n-k-1}{2}\rceil(n-k-\lceil\frac{n-k-1}{2}\rceil)} & \text{ si } n-k-1 \text{ est impair}. 
	\end{cases}$.
\ethm
\begin{proof}
	Posons  $\overline{\sigma_{n,k}}=\cup_{j=1}^{n-k-1}\{\sigma_{n,k,j}\}$ pour  tout $ k\leq n-1$. 
	En utilisant le Lemme \ref{lem51} et de la Remarque \ref{rem51},  on obtient
	\begin{equation}\label{eq51}
	\sum_{\pi\in \overline{\sigma_{n,k}}}q^{\crs(\pi)}=
	\begin{cases}
	1 & \text{ si } \ \ \ k=n-1;\\
	\displaystyle \sum_{j=1}^{n-k-1}q^{j(n-k-j)} & \text{ si } \ \ \ k\geq \frac{n}{2};\\
	\displaystyle 2\sum_{j=1}^{k}q^{\binom{j}{2}+j(k-j)}+(n-3k-1)q^{\binom{k}{2}} & \text{ si } \ \ \ \frac{n-1}{3}\leq k<\frac{n}{2};\\
	2\sum_{j=1}^{\lceil\frac{n-k-1}{2}\rceil}q^{\binom{j}{2}+j(k-j)}-\gamma_{n,k}(q)& \text{ si } \ \ \ 1 \leq k<\frac{n-1}{3};\\\
	n-1 & \text{ si } \ \ \ k=0.
	\end{cases}
	\end{equation}
	Puisque $ S_n(123,312)=\displaystyle\cup_{k=0}^{n-1}\overline{\sigma_{n,k}}$, alors $\displaystyle F_n(123,312;q)=\sum_{k=0}^{n-1}\sum_{\pi\in \overline{\sigma_{n,k}}}q^{\crs(\pi)}$. Nous obtenons l'expression souhaitée pour $F_n(123,312;q)$ de \eqref{eq51}.
\end{proof}

Grâce à la  relation \eqref{eq33} entre $F_n(123,231;q)$ et $F_n(123,312;q)$  de la Proposition \ref{pro441},  on peut aussi en déduire l'expression de $F_n(123,231;q)$.
\bco\label{coro4x1}
	Pour tout entier $n\geq 1$, nous avons
	\begin{eqnarray*}
		F_n(123,231;q)&=&1+\sum_{k=1}^{\lfloor\frac{n}{3}\rfloor}\left(2\sum_{j=1}^{k}q^{\binom{j}{2}+j(k-j)}+(n-3k)q^{\binom{k}{2}}\right) \\
		&&+\sum_{k=\lfloor\frac{n}{3}\rfloor+1}^{\lfloor\frac{n}{2}\rfloor}\left( 2\sum_{j=1}^{\lceil\frac{n-k}{2}\rceil}q^{\binom{j}{2}+j(k-j)}-\gamma_{n+1,k}(q)\right) +\sum_{k=\lceil\frac{n+1}{2}\rceil}^{n-1}\sum_{j=1}^{n-k}q^{j(n+1-k-j)}.
	\end{eqnarray*}
\eco

\subsection{Permutations $(213,132)$-interdites}\label{subsec2}
Fixons  $T=\{213,132\}$. En manipulant la structure de  $S_n(T)$,  nous pouvons trouver une relation de récurrence pour $F_{n,j}^k(T;q)$ qui nous permet d'obtenir $F_n(T;q)$, où $F_{n,j}^k(T;q)$ désigne le polynôme distributeur de la statistique $\crs$ sur l'ensemble $S_{n,j}^k(T):=S_{n}^k(T)\cap S_{n,j}(T)$ pour tous entiers $n$ et $i, j\in [n]$. 

Premièrement, pour tout $k\in [n]$, si $\sigma \in S_n^k(T)$, alors  $\sigma(k\cdots n)=12\cdots (n+1-k)$. 

Deuxièmement, pour tout $j \in [n]$, si $\sigma \in S_{n,j}(T)$ alors, $\sigma(1\cdots j)=(n+1-j)(n+2-j)\cdots n$. \\
En particulier, $S_{n}^1(T)=S_{n,n}(T)=\{12\cdots n\}$. De plus, si $k>1$ et $\sigma \in S_n^k(T)$, alors $j=\sigma^{-1}(n)<k$ et $S_{n,j}^k(T)\neq \emptyset\Leftrightarrow j<k$. Plus précisément,  $\sigma \in S_{n,j}^k \Leftrightarrow \sigma=(n+1-j)(n+2-j)\cdots n |\ \pi\ |12\cdots (n+1-k)$ où $\pi$ est une permutation $T$-interdite de $\{n+1-k,\ldots,n-j\}$. 

Alors, pour calculer $F_n(T;q)$, on utilise les identités évidentes suivantes:
\begin{eqnarray*}
	F_n^1(T;q)&=&F_{n,n}(T;q)=1,\\
	F_n^n(T;q)&=&F_{n-1}(T;q),\\
	F_n^k(T;q)&=&\sum_{j=1}^{k-1}F_{n,j}^{k}(T;q) \text{ pour } 1<k<n,\\
	\text{ et }   F_{n,j}(T;q)&=&\sum_{k=j+1}^{n}F_{n,j}^{k}(T;q) \text{ pour } 1\leq j<n.
\end{eqnarray*}
La récurrence pour $F_{n,j}^{k}(T;q)$ est donnée dans le théorème suivant.
\bthm\label{thm312}
Soit $n>1$. Pour tout entier $k$ satisfaisant  $1<k<n$, nous avons
\begin{equation*}
F_{n,j}^k(T;q)=\begin{cases}
q^{j(n-k)}F_{n-2j}^{k-j}(T;q) & \text{ si } j< n+1-k;\\
q^{j(j-1)}F_{n-2j}(T;q) & \text{ si } j=n+1-k;\\ 
q^{(n+1-k)(j-1)}F_{n-2(n+1-k),j-(n+1-k)}(T;q) & \text{ si } n+1-k<j<k-1;\\
q^{j(n-k)} & \text{ si } j= k-1;\\
0  & \text{ si } j\geq k.
\end{cases}
\end{equation*}
\ethm
\begin{proof} Pour  $1<k<n$, il est clair que nous avons  
	$$F_{n,j}^k(T;q)=0 \text{ si } j\geq k.$$  
	Le seul élément de $S_{n,k-1}^k(\T)$ est  $\pi=(n+2-k)\cdots (n-1)n12\cdots (n+1-k)$. De plus, de l'identité \eqref{eqarth}, puisque
	$\inv(\pi)=(k-1)(n+1-k)$, $\exc(\pi)=k-1$ et $\nes(\pi)=0$, alors $\crs(\pi)=(k-1)(n-k)$. Et par conséquent,
	\begin{equation*}
	F_{n,j}^k(T;q)=q^{j(n-k)} \text{ si } j=k-1.	
	\end{equation*}	
	Considérons maintenant une bijection de $S_{n,j}^k(\T)$ vers $S_{n-2,j-1}^{k-1}(\T)$ qui associe la permutation $\sigma=(n+1-j)(n+2-j)\cdots n |\ \pi\ |12\cdots (n+1-k)$ de $S_{n,j}^k(\T)$ à la permutation $\sigma'=(n-j)(n+1-j)\ldots (n-2) |\ \pi'\ |12\cdots (n-k)$ déduite de $\sigma$ en supprimant les lettres $n$ et $1$. Dans ce cas, 
	$$(\exc,\inv,\nes)(\sigma)=(1+\exc,n-2+k-j+\inv,k-1-j+\nes)(\sigma').$$ 
	Selon toujours l'identité  \eqref{eqarth}, on obtient $\crs(\sigma)=n-1+j-k+\crs(\sigma')$ et par conséquent
	\begin{equation}\label{eq45}
	F_{n,j}^k(213,132;q)=q^{n-1+j-k}F_{n-2,j-1}^{k-1}(213,132;q)
	\end{equation}
	Les relations souhaitées sur $F_{n,j}^k(T;q)$ pour $j<k-1$ se déduisent de la relation (\ref{eq45}). 
\end{proof}

\subsection{Permutations $(\tau_1,\tau_2)$-interdites, $(\tau_1,\tau_2)\in \{132,213\}\times \{231,312\}$}\label{subsec3}
Considérons d'abord le q-tableau   $(R_{n}^{k}(q))_{n,k}$ défini comme suit
\begin{equation} \label{q-tableau}	
\begin{cases}
R_{n}^n(q)=R_{n}^{n-1}(q)=1 & \\
R_{n}^{k}(q)=q^{\min\{k-1,n-1-k\}}R_{n-1}^{k}(q)+R_{n}^{k+1}(q) & \text{ si \ \ $1\leq k<n-1$} \\		
R_{n}^0(q)=R_{n-1}^{0}(q)+R_{n}^{1}(q) & 
\end{cases}
\end{equation}
Nous pouvons vérifier facilement que  $R_{n}^{k}(1)=2^{n-1-k}$ pour $0\leq k< n$, et  $\sum_{k=0}^{n}R_n^k(1)=2^{n}$. Ainsi, $(R_{n}^{k}(q))$ est un q-tableau des puissances de 2.  Nous présentons ici quelques valeurs de $(R_{n}^k(q))$ dans le Tableau \ref{table:qtable} pour des petites valeurs de $n$ et $k$.

\begin{table}[h]
	\centering
	
	\begin{tabular}{|c|llllll|}
		\hline 
		k& 0&1 & 2 & 3 & 4 \hspace{0.5cm} & 5 \\
		\hline
		0& $1$&&  & & &    \\
		1& $1$ &$1$&  & & &   \\
		2&$2$ &$1$& $1$&  & &\\
		3& $4$ &$2$& $1$& $1$ & & \\
		4& $7+q$ &$3+q$& $1+q$&  $1$ & &\\
		5& $11+4q+q^2$&$4+3q+q^2$& $1+2q+q^2$& $1+q$ & 1& 1\\
		\hline 
	\end{tabular}
	\label{table:qtable} 
	\caption{ Valeurs de $(R_{n}^k(q))$ pour $0\leq n\leq 6$ et $0\leq k\leq 3$} 
	
\end{table}

Nous allons montrer que ce q-tableau compte les permutations $T$-interdites selon le nombre de croisements, où $T$ est une des paires $\{213,312\},\{132,312\}, \{213,312\}$ et $\{213,312\}$ qui sont liées par les relations 
$\{132,312\}=rci(\{213,312\})$,  $\{132,231\}=rci(\{213,231\})$  et $\{213,312\}=\{213,231\}^{-1}$.

Pour toute permutation $\sigma \in S_n$ et pour tout entier $j\in [n]$, posons
\begin{itemize}
	\item $X_j(\sigma)=\{i<j/\sigma(i)\geq j\}$; 
	\item $Y_j(\sigma)=\{i+1<j/\sigma(i)\leq i \text{ et  } i+1\leq \sigma^{-1}(i+1)\}$; 
	\item $Z_j(\sigma)=\{(i,k)/i<k<\sigma(i)=k+1<\sigma(k) \text{ et } k+1\leq j\}$.
\end{itemize}
\blem\label{lem41} 	Pour toute permutation $\sigma $ et pour tout $1\leq j\leq |\sigma|+1$, nous avons $$\crs(\sigma^{(1,j)})=\crs(\sigma)+|X_j(\sigma)|+|Y_j(\sigma)|-|Z_j(\sigma)|.$$
\elem
\begin{proof}	Soit $\sigma$ une permutation et $1\leq j\leq |\sigma|+1$. Posons $\pi=\sigma^{(1,j)}$. Par définition,   
	\begin{equation*}
\pi(1)=j \text{  et	} \pi(1+i)=\begin{cases}
	\sigma(i)+1 & \text{ si } \sigma(i)\geq j;\\
	\sigma(i) & \text{ si } \sigma(i)< j.
	\end{cases}.
	\end{equation*}
	Les faits suivants viennent de cette définition.
	\begin{enumerate}
		\item $i\in X_j(\sigma)$ si et seulement si  $(1,1+i)$ est un croisement de $\pi$.	
		\item Si $i\in Y_j(\sigma)$ alors  $1+i\in L\!t(\pi)$, c'est-à-dire, $(i+1,\pi^{-1}(i+1))$ est un croisement de $\pi$.	
		\item Si $(i,k) \in Z_j(\sigma)$  alors $i+1<\pi(1+i)=1+k<\pi(1+k)$, c'est-à-dire, $i+1\in U\!t(\pi)$. 
		\item Si $(i,k)$ est un croisement de $\sigma$ qui n'est pas dans $Z_j(\sigma)$, alors $(i+1,k+1)$ est un croisement de $\pi$.
	\end{enumerate} 
	Par conséquent, $\crs(\sigma)-|Z_j(\sigma)|=\crs(\pi)-|X_j(\sigma)|-|Y_j(\sigma)|$. 
	Le lemme s'ensuit immédiatement.
\end{proof}

Soit $n$ un entier positif. Pour tout $k \in [n]$, nous dénotons par 
\begin{itemize}
	\item $G_n^{[k]}:=\{\sigma \in S_n|\sigma(n+1-i)=i, \ \ \forall\  i\  \in [k] \}$ (avec $G_n^{[0]}=S_n$ par convention);
	\item $H_n^{[k]}:=\{\sigma \in G_n^{[k]}|\sigma(n-k)\neq k+1\}$.
\end{itemize} 
Par ces notations, nous avons $G_{n}^{[n]}=G_{n}^{[n-1]}=\{n\ldots 21\}$ et  $G_{n}^{[k]}=H_n^{[k]}\cup G_n^{[k+1]}$ pour tout  $k<n-1$.
\blem\label{lcor41}
Pour tout $\sigma \in G_n^{[k]}$, nous avons $\crs(\sigma^{(1,k+1)})=\min\{k-1,n-k\}+\crs(\sigma)$.
\elem
\begin{proof}
	Puisque tout $\sigma \in G_{n}^{[k]}$ peut s'écrire sous la forme $\sigma =\pi|k\ldots 2 1$, où $\pi$ est une permutation  de $\{k+1,\ldots,n-1,n\}$, alors  $\sigma(i)\geq k+1$ pour tout $i\leq n-k$.
	Par conséquent, nous avons $Y_{k+1}(\sigma)=Z_{k+1}(\sigma)=\emptyset$ et
	\begin{equation*}
	X_{k+1}(\sigma)=\begin{cases}
		\{1,2,\ldots, k-1\} & \text{ si }\ \ \ k\leq \frac{n}{2};\\
		\{1,2,\ldots, n-1-k\} & \text{ si } \ \ \ k> \frac{n}{2}.
	\end{cases}
	\end{equation*}
Autrement dit, pour tout $k \in [n]$, $|Y_{k+1}(\sigma)|=|Z_{k+1}(\sigma)|=0$ et $|X_{k+1}(\sigma)|=\min\{k-1,n-1-k\}$.
Alors, notre lemme s'obtient du Lemme \ref{lem41}.
\end{proof}

Maintenant, nous notons  $G_{n}^k(T;q)=\sum_{\sigma \in G_{n}^{[k]}(T)} q^{\crs(\sigma)}$ pour tout $k\in [n]$ et pour tout ensemble de motifs $T$. Particulièrement, on a
\begin{equation*}
G_{n}^0(T;q)= F_{n}(T;q) \text{ et } G_{n}^1(T;q)= F_{n}^n(T;q).
\end{equation*}

\blem\label{lem46} Soit $T=\{213,312\}$.
	Pour tout $n\geq 1$ et pour tout $k<n-1$, 
	\begin{equation}\label{eq:lem46}
	G_{n}^k(T;q)=q^{\min\{k-1,n-1-k\}}G_{n-1}^k(T;q)+G_{n}^{k+1}(T;q).
	\end{equation}
\elem
\begin{proof}
	Premièrement, puisque $G_{n}^{[k]}(\T)=H_n^{[k]}(T)\cup G_n^{[k+1]}(T)$ pour $1\leq k<n-1$, nous obtenons 
	\begin{equation}\label{eqn3x}
	G_{n}^k(T;q)= \sum_{\sigma \in H_{n}^{[k]}(\T)}q^{\crs(\sigma)}+ G_{n}^{k+1}(T;q).
	\end{equation} 
	Deuxièmement, l'application de  $G_{n-1}^{[k]}(\T)$ vers $H_n^{[k]}(\T)$  qui associe $\sigma$ à $\sigma^{(1,k+1)}$ est bien définie et est bijective. Ainsi, en utilisant le Lemme \ref{lcor41}, nous obtenons
	\begin{equation}\label{eqn3y}
	\sum_{\sigma \in H_n^{[k]}(\T)}q^{\crs(\sigma)}=q^{\min\{k-1,n-1-k\}}\times G_{n-1}^k(T;q) .
	\end{equation}
	Notre lemme vient de la combinaison de (\ref{eqn3x}) et (\ref{eqn3y}).
\end{proof}

\bpro\label{pro44}	
	Soit $\tau \in \{213,132\}$. Pour tout  $n\geq 0$, nous avons
	\begin{equation*}\label{eq-main1}
	\sum_{\sigma \in G_n^{[k]}(312,\tau)}q^{\crs(\sigma)}=R_{n}^k(q) \text{ pour } 0 \leq k\leq n.
	\end{equation*}
\epro
\begin{proof} Soit $T=\{213,312\}$. Puisque $G_{n}^{[n]}(T)=G_{n}^{[n-1]}(T)=\{n\ldots 21\}$, nous avons 
	\begin{equation}\label{eq420}
	G_{n}^n(T;q)=G_{n}^{n-1}(T;q)=1.
	\end{equation} 
Rappelons que $S_n(T)=S_n^1(T)\cup S_n^n(T)$. Donc, 
\begin{equation*}
 F_{n}(T;q)=F_{n}^1(T;q)+F_{n}^{n}(T;q).
\end{equation*} 
Sachant que $F_{n}^1(T;q)=F_{n-1}(T;q)=G_{n-1}^0(T;q)$ et $F_{n}^{n}(T;q)=G_{n}^{1}(T;q)$, alors
\begin{equation}\label{eq:422}
G_{n}^0(T;q)=G_{n}^0(T;q)+G_{n}^{1}(T;q).
\end{equation}
Combinant les identités \eqref{eq:lem46}, \eqref{eq420} et \eqref{eq:422}, nous obtenons la relation suivante qui est équivalente à (\ref{q-tableau})
\begin{equation} \label{eq42}
\begin{cases}
G_{n}^{n-1}(T;q)&=G_{n}^{n}(T;q)=1.\\
G_{n}^{k}(T;q)&=q^{\min\{k-1,n-1-k\}} G_{n}^k(T;q)+G_{n}^{k+1}(T;q) \text{ pour tout $1\leq k<n-1$}.\\
G_{n}^{0}(T;q)&=G_{n-1}^{0}(T;q)+G_{n}^{1}(T;q).
\end{cases}
\end{equation}
Autrement dit, $G_{n}^{k}(213,312;q)=R_{n}^k(q)$ pour tout $k\geq 0$. La proposition s'ensuit dès que nous  utilisons
la relation $\{132,312\}=rci(\{213,312\})$.
\end{proof}
\bthm\label{thm46}
	Soit $\tau \in \{213,132\}$. Pour tout $n\geq 1$,	nous avons
	\begin{equation*}\label{eq3x}
		\sum_{\sigma \in S_n(312,\tau)}q^{\crs(\sigma)}=R_{n}^0(q) \text{ et } \sum_{\sigma \in S_n(231,\tau)}q^{\crs(\sigma)}=R_{n+1}^1(q) 
	\end{equation*}
\ethm
\begin{proof}
La première identité du théorème vient de la Proposition \ref{pro44}.	Puisque\\
 $f_{n+1}(S_{n}(213,231))=S_{n+1}^{n+1}(213,312)$, 
	\begin{equation*}
	\sum_{\sigma \in S_n(231,213)}q^{\crs(\sigma)}=F_{n+1}^{n+1}(312,213;q)=G_{n+1}^1(312,213;q)=R_{n+1}^1(q).
	\end{equation*}
	Comme $\{231,132\}=rci(\{231,213\})$ et $\{312,132\}=rci(\{312,213\})$, alors, on obtient la proposition.
\end{proof}

Nous allons conclure cette section par une relation entre  $F(312,\tau; q,z)$ et $F(231,\tau; \\ q,z)$, $\tau \in \{213,132\}$.
\bpro \label{pro4x}Pour tout $n\geq 2$, nous avons 
\begin{eqnarray} 
F_n(312,213;q)&=&F_{n-1}(312,213;q)+F_{n-1}(231,213;q). \label{eq34}
\end{eqnarray}
\epro
\begin{proof}
	Sachant que $S_n(312,213)=S_{n}^1(312,213)\cup S_n^{n}(312,213)$, nous avons
	\begin{eqnarray*}
		F_n(312,213;q)&=&F_{n}^{1}(312,213;q)+F_{n}^{n}(312,213;q).
	\end{eqnarray*}
Puisque $F_{n}^{1}(312,213;q)=F_{n-1}(312,213;q)$  et $F_{n}^{n}(312,213;q)=F_{n-1}(231,213;q)$ ( propriétés (i) et (iv) de la Proposition \ref{prop3x1}), alors nous obtenons la Proposition \ref{pro4x}.
\end{proof}
\bthm \label{thm4135}
Pour tout $(\tau,\tau') \in   \{132,213\}^2$, nous avons  la relation
\begin{equation*}
F(312,\tau; q,z)= 1+\left( \frac{z}{1-z}\right)F(231,\tau'; q,z). 	
\end{equation*} 
\ethm
\begin{proof}
	Depuis la récurrence \eqref{eq34}, nous obtenons l'équation 
	\begin{equation*}
	F(312,213;q,z)=1+z+z\left( F(312,213;q,z)+F(231,213;q,z)-2\right)
	\end{equation*}
	qui est équivalente à
	\begin{equation}\label{eqx3}
	F(312,213; q,z)= 1+\left( \frac{z}{1-z}\right)F(231,213; q,z). 	
	\end{equation} 
	On complète la preuve du théorème par $rci$.
\end{proof}

\section{Conclusion}
Nous avons énuméré les permutations interdisant deux motifs de $S_3$ selon le nombre de croisements. Pour chaque paire de motifs $T$ de $S_3$, nous avons manipulé la structure de $S_n(T)$ afin de trouver une relation de récurrence sur $F_n(T,q)$, la distribution de la statistique $\crs$ sur $S_n(T)$. A partir de la relation de récurrence ainsi trouvée, nous avons essayé d'en déduire soit l'expression explicite de $F_n(T,q)$, soit la forme close de $F(T;q,z)$.  Le Tableau \ref{tab32} suivant résume les classes des paires de motifs ainsi que les références des résultats obtenus.

\begin{table}[h]
	\centering
	
	\begin{tabular}{|l|l|l|}
		\hline
		\textbf{Classe}  & \textbf{Résultat} & \textbf{Référence}\\
		\hline
		$\{321, 231\}$ & Forme close de $F(321, 231;q,z)$ & Thm. \ref{thm413}\\
		\hline
$\{123, 213\}$, $\{123, 132\}$ & Forme close de $F(-;q,z)$ & Rem. \ref{remx}\\
\hline
$\{321, 213\}$, $\{321, 132\}$ & Expression de $F_n(-;q)$ & Thm. \ref{prop35} \\
\hline
$\{123, 312\}$  & Expression de $F_n(123, 312;q)$ & Thm. \ref{thm311}\\
\hline
$\{123, 231\}$  & Expression de $F_n(123, 231;q)$ & Cor. \ref{coro4x1}\\
\hline
$\{213, 132\}$ & Récurrence sur $F_n(213, 132;q)$ & Thm. \ref{thm312}\\
\hline
$\{213, 231\}$, $\{132, 231\}$ & q-tableau $(R_n(q))$& Thm. \ref{thm46} \\
\hline
$\{213, 312\}$, $\{132, 312\}$ & q-tableau $(R_n(q))$& Thm. \ref{thm46} \\
\hline
$\{312, 231\}$, $\{321, 312\}$ & Trivial (Sans croisement) & -\\
\hline
$\{123,321\}$ & Trivial &  -\\
\hline
	\end{tabular}
\caption{Classement des paires de motifs $T$ de $S_3$ selon la distribution de $\crs$.}\label{tab32}	
\end{table}

\addcontentsline{toc}{chapter}{Conclusion générale et perspectives}
	\chapter*{Conclusion générale et perspectives}
Dans cette thèse, nous avons essayé de répondre une partie des questions évoquées dans l'introduction. Plus précisément, en utilisant les fonctions génératrices et des bijections, nous avons effectué une étude combinatoire des permutations évitant un ou deux motifs de longueur 3 selon le nombre de croisements. Comme résultats, nous avons trouvé des équidistributions et des interprétations combinatoires.

Sur ce, nous avons étudié une bijection  $\Theta:S_n(321)\rightarrow S_n(132)$ qui est originalement construite par Elizalde et Pak \cite{ElizP}  et nous avons proposé une nouvelle formulation de $\Theta$. L'intérêt de notre formulation n'est pas seulement en terme de complexité, mais elle nous permet aussi de prouver que $\Theta$ conserve, non seulement le nombre de points fixes et le nombre d'excédances, mais aussi le nombre de croisements. Grâce à la liaison entre la bijection $\Theta$ de Elizalde et Pak et la bijection $\Gamma$ de Robertson \cite{ARob, ARob2} prouvée par Saracino \cite{Sarino}, nous pouvons en déduire que la bijection $\Gamma$  conserve également le nombre de points fixes.
 
Sur les permutations interdisant un motif de $S_3$, nous avons utilisé la bijection $\Theta$ et l'involution trivial miroir-complément-inverse $rci$  comme outils fondamentaux pour prouver l'équidistribution du nombre de croisements sur  $S_n(\tau)$, pour  $\tau\in \{132, 213, 321\}$.  
A travers le $q,p$-Catalan de Randrianarivony \cite{ARandr}, nous avons obtenu le développement en fraction continue suivant
\begin{equation*}
\sum_{\sigma\in S(321)} q^{\crs(\sigma)}z^{|\sigma|}=\sum_{\sigma\in S(132)} q^{\crs(\sigma)}z^{|\sigma|}=\sum_{\sigma\in S(213)} q^{\crs(\sigma)}z^{|\sigma|}=\frac{1}{1-\displaystyle\frac{z}{				
		1-\displaystyle\frac{z}{
			1-\displaystyle\frac{qz}{								
				1-\displaystyle\frac{qz}{
					1-\displaystyle\frac{q^2z}{
						1-\displaystyle\frac{q^2z}{				
							\ddots}
					}
			}}		
}}}. 
\end{equation*}  
 Nous avons également construit une bijection de $S_{n-1}(231)$ vers $S_n^n(312)$ qui conserve le nombre de croisements afin de trouver une relation entre les distributions du nombre de croisements sur $S_n(231)$ et $S_n(312)$ suivante:
 \begin{align}\label{relconc}
 	 \sum_{\sigma\in S(312)} q^{\crs(\sigma)}z^{|\sigma|}&= \frac{1}{1-z\displaystyle \sum_{\sigma\in S(231)} q^{\crs(\sigma)}z^{|\sigma|}}.	
 \end{align} 
 Il est à noter que, pour l'instant, nous n'avons trouvé aucune information sur la distribution de $\crs$ sur $S_n(123)$.
 
 Sur les permutations interdisant deux motifs de $S_3$, nous avons manipulé les structures de ces familles d'objets, puis utilisé des bijections  pour  ainsi établir des relations de récurrences sur les distributions du nombre de croisements. Dans plusieurs situations, la décomposition du nombre de croisements en termes du nombre d'excédances, du nombre d'inversions et du nombre d'imbrications, un résultat prouvé par Médicis et Viennot \cite{MedVienot} et Randrianarivony \cite{ARandr}, nous a grandement facilité les tâches. Comme résultats, nous avons trouvé des nouvelles interprétations combinatoires  des triangles \href{https://oeis.org/A076791}{A076791} et \href{https://oeis.org/A299927}{A299927} de OEIS \cite{OEIS}. 
 
La première perspective de cette thèse est de trouver les distributions du nombre de croisements sur $S_n(\tau)$ pour  $\tau\in \{231, 312, 123\}$ car elles restent encore ouvertes. La relation (\ref{relconc}) ci-dessus que nous avons trouvée  nous servira comme point de départ.
L'étude combinatoire des permutations restreintes selon le nombre d'imbrications sera aussi une autre perspective logique de notre recherche, une étude que nous avons déjà introduite dans \cite{Rakot2}. Sachant que le nombre de croisements et le nombre d'imbrications sont équidistribués sur $S_n$, il est important aussi de savoir comment ces deux statistiques seront liées sur $S_n(T)$, où $T$ est un ensemble de motifs quelconques. Pencher aussi vers les applications est intéressant. Dans ce cas, on peut se référer au travail de Corteel \cite{Cort} qui a mis en évidence son application en Physique.

\backmatter

\newpage
 \thispagestyle{empty} 
\underline{\textbf{Doctorant :}}

Paul Mazoto RAKOTOMAMONJY

rpaulmazoto@gmail.com

+26134 40 434 64.

\begin{center}
\textbf{\textcolor{blue}{COMBINATOIRE DES PERMUTATIONS RESTREINTES SELON LE NOMBRE DE CROISEMENTS}}
\end{center}
\textbf{Résumé}\\
Dans cette thèse, nous avons introduit et effectué une étude combinatoire des permutations interdisant un ou deux motifs de longueur 3 selon la statistique nombre de croisements.  Pour cela, nous avons manipulé une bijection d'Elizalde et Pak et construit d'autres bijections qui conservent le nombre de croisements. Comme résultats, nous avons trouvé, à travers ces bijections, diverses relations sur les distributions du nombre de croisements sur  les permutations restreintes ainsi que des interprétations combinatoires  en termes du nombre de croisements sur les permutations à motifs interdits de  certains triangles bien connus dans la littérature.

\begin{center}
\textbf{Mots clés :} Permutation restreinte, bijection, fonction génératrice, statistique nombre de croisements, interprétation combinatoire.
\end{center}

\vspace{0.5cm}
\begin{center}
\textbf{\textcolor{blue}{COMBINATORIAL OF RESTRICTED PERMUTATIONS ACCORDING TO THE NUMBER OF CROSSINGS}}
\end{center}

\textbf{Abstract}\\
In this thesis, we introduced and carried out a combinatorial study of permutations that avoid one or two patterns of length 3 according to the statistic number of crossings.  For this purpose, we manipulated a bijection of Elizalde and Pak and constructed other bijections that preserve the number of crossings. As results, we found, throughout these bijections, various relationships on the distributions of the number of crossings on restricted permutations as well as combinatorial interpretations in terms of the number of crossings on permutations with forbidden patterns of some well known triangles in the literature.

\begin{center}
\textbf{Keywords:} Restricted permutation, bijection, generating function, statistic number of crossings, combinatorial interpretation.
\end{center}

\underline{\textbf{Directeur de thèse:}}

Pr Arthur RANDRIANARIVONY

arthur.randrianarivony@gmail.com

+26134 48 997 37.

\end{document}